\documentclass[11pt,a4paper,twoside]{article}

\usepackage[notext]{kpfonts}
\usepackage{baskervald}
\usepackage[T1]{fontenc}
\usepackage[utf8]{inputenc}
\usepackage{lmodern}
\usepackage{microtype}
\usepackage[a4paper,margin=30mm,headheight=14pt]{geometry}
\usepackage{amsmath,amssymb,amsthm,mathtools}
\usepackage{tikz-cd}
\usepackage{verbatim}
\usepackage[all]{xy}
\usepackage{enumerate}
\usepackage{dsfont}
\usepackage{fancyhdr}
\usepackage{mathtools}
\usepackage{xcolor}

\usepackage[titles]{tocloft}
\usepackage[colorlinks=true,linkcolor=blue!55!black,citecolor=blue!55!black,
            urlcolor=blue!55!black, backref=page]{hyperref}

\newtheorem{thm2}{Theorem}
\newtheorem{cor2}{Corollary}
\newtheorem{conj2}{Conjecture}

\newtheorem{thm}{Theorem}[section]
\newtheorem{cor}[thm]{Corollary}
\newtheorem{lem}[thm]{Lemma}
\newtheorem{prop}[thm]{Proposition}

\theoremstyle{definition}
\newtheorem{defn}[thm]{Definition}

\newtheorem{ex}[thm]{Example}
\newtheorem{rmk}[thm]{Remark}

\DeclareMathOperator{\ch}{ch} \DeclareMathOperator{\Spec}{Spec}

 \DeclareMathOperator{\rk}{rk}
  
\DeclareMathOperator{\Sym}{Sym}

\DeclareMathOperator{\prim}{prim}

\newcommand{\C}{\ensuremath\mathds{C}}
\newcommand{\R}{\ensuremath\mathrm{R}}

\newcommand{\Z}{\ensuremath\mathds{Z}}
\newcommand{\Q}{\ensuremath\mathds{Q}}

\newcommand{\fa}{\ensuremath\mathfrak{a}}

\newcommand{\FF}{\ensuremath\mathcal{F}}
\newcommand{\h}{\ensuremath\mathfrak{h}}

\newcommand{\PP}{\ensuremath\mathds{P}}

\renewcommand{\D}{\ensuremath\mathrm{D}}
\newcommand{\HH}{\ensuremath\mathrm{H}}

\newcommand{\CH}{\ensuremath\mathrm{CH}}

\newcommand{\id}{\ensuremath\mathrm{id}}
\newcommand{\Gr}{\ensuremath\mathrm{Gr}}
\newcommand{\Stab}{\ensuremath\mathrm{Stab}}
\renewcommand{\top}{\ensuremath\mathrm{top}}

\newcommand{\QQ}{\mathbb{Q}}
\newcommand{\NN}{\mathbb{N}}
\newcommand{\LLL}{\mathbb{L}}

\newcommand{\UU}{\mathcal U}

\newcommand{\XX}{\mathcal X}
\newcommand{\YY}{\mathcal Y}
\newcommand{\MM}{\mathcal M}

\newcommand{\VV}{\mathcal V}
\newcommand{\OO}{\mathcal O}
\newcommand{\CC}{\mathcal C}
\newcommand{\Zz}{\mathcal Z}
\newcommand{\TT}{\mathcal T}
\newcommand{\WW}{\mathcal W}

\newcommand{\PPP}{\mathcal P}

\newcommand{\Sy}{\mathfrak S}
\newcommand{\EE}{\mathcal E}
\newcommand{\GDCH}{\operatorname{GDCH}}

\newcommand{\ima}{\hbox{Im}}

\renewcommand{\1}{\mathop{\mathds{1}}\nolimits} 

\renewcommand{\bar}{\overline}

\newcommand{\cart}{\ar@{}[dr]|\square} 
\renewcommand{\tilde}{\widetilde}

\newcommand{\B}{\ensuremath\mathrm{B}}
\newcommand{\GDB}{\operatorname{GDB}}
\newcommand{\Mrat}{\ensuremath\MM_{\mathrm{rat}}}

\newcommand{\Mhom}{\ensuremath\MM_{\mathrm{hom}}}
\newcommand{\CGr}{\ensuremath\mathrm{C}\mathrm{Gr}}
\newcommand{\Ku}{\ensuremath\mathcal{K}u}
\newcommand{\van}{\ensuremath\mathrm{van}}
\newcommand{\taut}{\ensuremath\mathrm{taut}}
\newcommand{\cl}{\ensuremath\mathrm{cl}}

\AtBeginDocument{%
	\def\MR#1{}
}

\newif\ifHideFoot
\HideFoottrue  
\HideFootfalse  

\ifHideFoot

\newcommand{\Lie}[1]{}
\newcommand{\Charles}[1]{}
\newcommand{\Robert}[1]{}
\newcommand{\Carlo}[1]{}

\else 

\newcommand{\marg}[1]{\normalsize{{
			\color{red}\footnote{{\color{blue}#1}}}{\marginpar[\vskip
			-.25cm{\color{red}\hfill$\Rightarrow$\tiny\thefootnote}]{\vskip
				-.2cm{\color{red}$\Leftarrow$\tiny\thefootnote}}}}}
\newcommand{\Lie}[1]{\marg{(Lie) #1}}
\newcommand{\Charles}[1]{\marg{(Charles) #1}}
\newcommand{\Robert}[1]{\marg{(Robert) #1}}
\newcommand{\Carlo}[1]{\marg{(Carlo) #1}}

\fi

\title{\textbf{The generalized Franchetta conjecture for some hyper-K\"ahler varieties, III. }}
\author{Lie Fu, Robert Laterveer, Carl Mazzanti, Charles Vial}
\date{ }

\begin{document}

\maketitle

\begin{abstract}
The generalized Franchetta conjecture states that for a locally complete family of polarized hyper-K\"ahler varieties, the generically defined cycles on any fiber of the universal family over the moduli stack inject into cohomology via the cycle class map. We prove the generalized Franchetta conjecture for the family of double EPW sextics discovered by O'Grady and for the family of double EPW cubes discovered by Iliev--Kapustka--Kapustka--Ranestad. We also prove a Franchetta-type result for the relative square of the universal family of double EPW sextics, which allows to relate double EPW sextics and Gushel--Mukai fourfolds on the level of motives.

The proof involves the geometry of Gushel--Mukai fourfolds, whose Kuznetsov component allows to realize double EPW sextics and cubes as moduli spaces of Bridgeland-stable objects, relatively over the parameter space of Gushel--Mukai fourfolds. A key step is that the canonical Chow--K\"unneth decomposition of  special Gushel--Mukai fourfolds is multiplicative; we use this to show that generically defined cycles on the third power of a special Gushel--Mukai fourfold inject into cohomology via the cycle class map. Along the way, we also show that the canonical Chow--K\"unneth decomposition of any Gushel--Mukai sixfold is multiplicative.
\end{abstract}

\renewcommand{\thefootnote}{}
\null
\renewcommand{\thefootnote}{\arabic{footnote}}

{\small\setlength{\parskip}{0pt}\tableofcontents}

\newpage
\section{Introduction}

\subsection{The Franchetta property}
Let $f\colon \XX\to S$ be a smooth projective (representable) morphism between smooth varieties (or more generally, Deligne--Mumford stacks) over the field of complex numbers, and let $X=X_s$ be the fibre of $f$ over a closed  point $s\in S(\C)$. We define the group of \textit{generically defined cycles} on $X$, with respect to the family $\XX\to S$, defined as the image of the Gysin restriction map from the Chow group of the family:
\[
\GDCH^*_S(X)\coloneqq\ima\bigl(\CH^*(\XX)\to \CH^*(X)\bigr).
\]
Here, as everywhere in this paper, Chow groups are taken with rational coefficients, and Chow groups of stacks as well as their functoriality are as in \cite{Vistoli-IntersectionTheoryStack, Kresch-ChowStack}. The family $f\colon \XX\to S$ is said to satisfy the \emph{Franchetta property in codimension $i$}, if the cycle class map
\[
\cl\colon \GDCH^i_S(X)\longrightarrow \HH^{2i}(X,\QQ)
\]
is injective for all fibres (or equivalently by Lemma~\ref{lem:spread}, for the very general fibres), and to satisfy the \emph{Franchetta property} if it does so in every codimension. 

We emphasize that the group $\GDCH^*_S(X)$ depends crucially on the base $S$ of the family (or rather, on the generic fiber $\XX_{\eta}$ with  $\eta=\Spec \C(S)$). For a dominant map  $S'\to S$, we have $\GDCH^*_S(X)\subseteq \GDCH^*_{S'}(X)$, and consequently, the Franchetta property for the dominant base-change family $\XX\times_S S'\to S'$ implies the Franchetta property for $\XX\to S$ (Lemma \ref{lem:base-change}). On the other hand, for a locally closed substack $S'\subseteq S$ containing the point $s$,  there is no implication in either direction between the Franchetta property of $\XX\to S$ and that of $\XX\times_{S}S'\to S'$.

\subsection{The Franchetta conjecture: from K3 surfaces to hyper-K\"ahler varieties}
Franchetta's work \cite{Franchetta}, together with its further developments \cite{Harer-H2MappingClassGroup, Mestrano-FranchettaForte, Kouvidakis-PicardOfUniversalCurve}, implies that the relative Picard group of the universal curve over the moduli stack of curves of genus $g\geq 2$ is generated by the canonical class. 
Inspired by this phenomenon for curves, O'Grady \cite{OGrady-K3} conjectured that the universal family of polarized K3 surfaces of genus $g\geq 2$ over its moduli space satisfies the Franchetta property, i.e.~the generically defined 0-cycles are generated by the so-called Beauville--Voisin class \cite{BV}. This was established by Pavic--Shen--Yin \cite{PSY} using Mukai models for $g\leq10$ and $g=12,13,16,18,20$. It has also been proven by Fu--Laterveer for $g=14$ \cite{FL}, by Lu for $g=11$ \cite{Lu} and for $g=22$ \cite{Lu2}, and by Meng for $g=21$ \cite{Meng}.

Compact hyper-K\"ahler varieties are the natural higher-dimensional generalizations of K3 surfaces (see \cite{beauville1983varietes, Huy99} for the basics). We asked in our previous work \cite{FLV1} whether O'Grady's conjecture can be extended to hyper-K\"ahler varieties:

\begin{conj2}[Generalized Franchetta conjecture \cite{FLV1}]
\label{conj:GFC}
Let $\MM$ be the moduli stack of a locally complete family of polarized hyper-K\"ahler varieties. Then the universal family $\XX\to\MM$ satisfies the Franchetta property, namely, the cycle class map restricted to generically defined cycles $\cl\colon \GDCH^i_{\MM}(X)\to \HH^{2i}(X,\QQ)$ is injective for any member $X$ in the family.
\end{conj2}

\paragraph{Known cases.} Besides K3 surfaces of low genus as mentioned above, Conjecture~\ref{conj:GFC} is established for the Fano varieties of lines on cubic fourfolds (and even for their relative squares) \cite{FLV1}, for the hyper-K\"ahler eightfolds of Lehn--Lehn--Sorger--van Straten \cite{FLV2}, and for Debarre--Voisin hyper-K\"ahler fourfolds \cite[Corollary 5.4]{Akazawa}. See also \cite{BeLi} for the study of a cohomological version of the conjecture.

\subsection{Main results}
There are very few available concrete geometric constructions of locally complete families of polarized hyper-K\"ahler varieties. Apart from the aforementioned families of Beauville--Donagi \cite{BD85}, Lehn--Lehn--Sorger--van Straten \cite{LLSvS} and Debarre--Voisin \cite{DebarreVoisin},
the so-called \textit{double EPW sextics} constructed by O'Grady \cite{OGrady06, OGrady13} and the so-called \textit{double EPW cubes} constructed by Iliev--Kapustka--Kapustka--Ranestad \cite{IKKR} seem to be the most interesting ones to study. 

Recall that these two families of hyper-K\"ahler varieties are both constructed from a generic Lagrangian subspace $A$ of the symplectic space $\bigwedge^3V_6$, where $V_6$ is a fixed 6-dimensional complex vector space. The Eisenbud--Popescu--Walter (EPW) sextic $Z_A\subseteq\PP(V_6)$ constructed in \cite{EPW} is the hypersurface
$$Z_A:= \bigg\{[V_1]\in \PP(V_6) ~\vert~ A\cap(V_1\wedge\bigwedge^2V_6)\neq0 \bigg\}.$$
For $A$ generic, $Z_A$ is an integral sextic hypersurface whose singular locus is a smooth surface, the double cover 
$$f\colon X_A\to Z_A$$
branched along that surface is a hyper-K\"ahler fourfold of $\mathrm{K}3^{[2]}$-type, called \textit{double EPW sextic} of O'Grady and $h_X\coloneqq f^*c_1(\OO(1))$ is a polarization of Beauville--Bogomolov square $2$ and divisibility 1; the resulting family is locally complete of dimension $20$ \cite{OGrady06,OGrady13}. 
Similarly, a generic Lagrangian subspace $A$ determines a degeneracy locus in $\Gr(3,V_6)$:
$$Z'_{A}:=\bigg\{[U_3]\in \Gr(3,V_6) ~\vert~ \dim \bigg(A\cap (\bigwedge^2U_3	\wedge V_6)\bigg)\geq 2\bigg\}.$$
The double cover $C_A$ of $Z'_{A}$ branched along its singular locus (which is a smooth threefold), is the \textit{double EPW cube} of Iliev--Kapustka--Kapustka--Ranestad \cite{IKKR}, which is a hyper-K\"ahler sixfold of $\mathrm{K}3^{[3]}$-type with polarization of Beauville--Bogomolov square 4 and divisibility 2, which likewise moves in a locally complete family of dimension $20$.

Let $\MM_{\widetilde{\mathrm{EPW}}}$ (resp.~$\MM_{\mathrm{IKKR}}$) denote the moduli stack of double EPW sextics (resp.~of double EPW cubes). Note that 
$\MM_{\widetilde{\mathrm{EPW}}}$ is an open substack of ${}^{[2]}\!\mathscr{M}^{(1)}_2$, which is the moduli stack of polarized hyper-K\"ahler varieties of $\mathrm{K}3^{[2]}$-type with polarization of Beauville--Bogomolov square 2 and divisibility 1, following the notation of \cite{DebarreMacri-PeriodMapHK}. Similarly, $\MM_{\mathrm{IKKR}}$ is an open substack of the moduli stack ${}^{[3]}\!\mathscr{M}^{(2)}_4$ of polarized hyper-K\"ahler varieties of $\mathrm{K}3^{[3]}$-type with polarization of Beauville--Bogomolov square 4 and divisibility 2.

 The main goal of this paper is to establish Conjecture~\ref{conj:GFC} for these two families, strengthening previous results modulo algebraic equivalence of Bolognesi--Laterveer \cite{BL}. In fact, for the first family, our result is slightly stronger in the sense that we can prove the Franchetta property even after the dominant base change $\MM_{\mathrm{EPW}}\to \MM_{\widetilde{\mathrm{EPW}}}$, where $\MM_{\mathrm{EPW}}$ is the moduli stack\footnote{$\MM_{\mathrm{EPW}}$ is generically an algebraic space, while $\MM_{\widetilde{\mathrm{EPW}}}$ is generically a (trivial) $\Z/2\Z$-gerbe, corresponding to the everywhere defined covering involution of double EPW sextics. It follows that $\MM_{\mathrm{EPW}}$ can be generically identified with the coarse moduli space of $\MM_{\widetilde{\mathrm{EPW}}}$.} of (sufficiently nice) Lagrangian spaces (or equivalently, the moduli stack of EPW sextic hypersurfaces).

\begin{thm2}\label{thm:main} \hfill
\begin{enumerate}[(i)]
    \item The universal family of double EPW sextics $\XX\to \MM_{\mathrm{EPW}}$ satisfies the Franchetta property. 
    \item The relative square $\XX^{2}_{/\MM_{\mathrm{EPW}}} \to \MM_{\mathrm{EPW}}$ satisfies the Franchetta property in codimension~$\le 4$, and the relative cube $\XX^{3}_{/\MM_{\mathrm{EPW}}} \to \MM_{\mathrm{EPW}}$ satisfies the Franchetta property in codimension~$\le 2$.
     \item The universal family of double EPW cubes over the moduli stack $\MM_{\mathrm{IKKR}}$ satisfies the Franchetta property.
\end{enumerate}
\end{thm2}

In our previous work \cite{FLV1, FLV2} proving Conjecture~\ref{conj:GFC} for the Beauville--Donagi family and the Lehn--Lehn--Sorger--van Straten family, the geometry of cubic fourfolds, which are \textit{Fano manifolds of cohomological K3 type} (see \cite{FLV-MCK}), plays an essential role. The basic proof strategy is to transfer the problem of algebraic cycles on families of hyper-K\"ahler varieties $\XX\to\MM$ to the study of algebraic cycles on families of cubic fourfolds $\YY\to B$, and eventually reduce to establishing a certain Franchetta property of the relative fourth power of the universal family of cubic fourfolds  $\YY^{4}_{/B}\to B$, where the existence of the so-called \textit{multiplicative Chow--K\"unneth decomposition} (as introduced in \cite{SV} and established in \cite{Diaz, FLV-MCK} for cubic hypersurfaces) is a crucial intermediate step.

In the case of double EPW sextics and double EPW cubes, there are also naturally associated Fano manifolds of cohomological K3 type, namely, the (even-dimensional) \textit{Gushel--Mukai} (GM) varieties. GM varieties are systematically studied in the series of works of Debarre--Kuznetsov \cite{DK1, DK2, DK3, DK4, DK5}; see also the survey \cite{Debarre}. Now we briefly recall their construction, present our results on GM varieties, and explain their relation to Theorem \ref{thm:main}.

\paragraph{Gushel--Mukai fourfolds.}

Let $V_5$ be a 5-dimensional vector space, $\CGr(2, V_5)$ the projective cone over the Grassmannian $\Gr(2, V_5)$, and $\OO(1)$ the Pl\"ucker polarization. By definition, a GM fourfold is a smooth transversal intersection of $\CGr(2, V_5)$ with two hyperplanes and a quadric:
\[Y=\CGr(2, V_5)\cap H_1 \cap H_2\cap Q.\] 
GM fourfolds come in two types: if the vertex $O$ of the projective cone $\CGr(2, V_5)$ lies in the two hyperplanes $H_1$ and $H_2$, the resulting GM fourfold is called \textit{special}; otherwise, if $O\notin H_1\cap H_2$, the corresponding GM fourfold is called \textit{ordinary}. Special GM fourfolds are indeed specializations of ordinary ones.

Denote $$\bar B\coloneqq\PP\HH^0(\CGr(2,V_5),\OO(1)^{\oplus2}\oplus\OO(2)),$$ 
and let $B\subseteq\bar B$ be the dense open subset parameterizing smooth GM fourfolds.
Let $\YY\to B$ be the universal family of GM fourfolds over $B$.
The proof of Theorem \ref{thm:main} uses geometric relations between double EPW sextics/cubes and GM fourfolds, and is ultimately reduced to the following results about \textit{special} GM fourfolds.

\begin{thm2}[Theorem~\ref{thm:MCK}]
\label{thm:MCK-intro}
The canonical Chow--K\"unneth decomposition of a \emph{special} GM fourfold, constructed in Definition~\ref{defn:CK}, is multiplicative.
\end{thm2}

\begin{thm2}[Proposition~\ref{prop:Y3}]
\label{thm:Y3-intro}
Let the universal family of GM fourfolds $\YY\to B$ be as above. 
Let $b\in B$ be such that $Y=Y_b$ is a \emph{special} GM fourfold. Then the cycle class map $$\cl \colon \GDCH^i_B(Y^m) \to \HH^{2i}(Y^m, \Q)$$ is injective for $m\leq 3$ (and $i$ arbitrary), and for $m\leq 6$ and $i\le 8$.
\end{thm2}

\begin{rmk}
\label{rmk:ConditionSpecialIsNecessary}
    The hypothesis of being special is essential in both Theorem~\ref{thm:MCK-intro} and Theorem~\ref{thm:Y3-intro}: it turns out that the canonical Chow--K\"unneth decomposition of a very general GM fourfold is \textit{not} multiplicative (Remark~\ref{rmk:noMCK}), hence the modified small diagonal provides a non-zero element in $\GDCH^*_B(Y^3)$ with trivial cohomology class.
    Let us also clarify that Theorem~\ref{thm:Y3-intro} does \textit{not} say that the relative cube of the universal family of special GM fourfolds satisfies the Franchetta property. Instead, Theorem~\ref{thm:Y3-intro} says that for a special GM fourfold $Y$, any homologically trivial cycle on $Y^3$ that is generically defined over $B$ (i.e.~for \textit{all} GM fourfolds, not just the special ones) is rationally equivalent to zero.
\end{rmk}

As a by-product of our results, we can establish the following analogue of Theorem~\ref{thm:MCK-intro} and Theorem~\ref{thm:Y3-intro} for Gushel--Mukai sixfolds, which are also of cohomological K3 type. By definition, GM sixfolds are smooth hypersurfaces in the linear system $|\OO(2)|$ on the Grassmannian cone $\CGr(2, V_5)$:
\[W=\CGr(2, V_5)\cap Q,\] 
or in other words, double covers of $\Gr(2, V_5)$ ramified along a smooth ordinary GM fivefold.
\begin{thm2}
\label{thm:MCK6-intro}
Let $S\subseteq |\OO_{\CGr(2, V_5)}(2)|$ be the open subset parameterizing smooth GM sixfolds. Let $\WW\to B$ be the universal family of GM 6-folds. 
\begin{enumerate}[(i)]
    \item The canonical Chow--K\"unneth decomposition (Definition~\ref{defn:CK6}) of a GM sixfold is multiplicative.
    \item For any $s\in S$, let $W=W_s$ be the fiber. Then the cycle class map $$\cl \colon \GDCH^*_S(W^m) \to \HH^*(W^m, \Q)$$ is injective for $m\leq 3$.
\end{enumerate}
\end{thm2}

\begin{rmk}
    Theorem~\ref{thm:MCK6-intro} will be proved in Section~\ref{S:sixfolds} (Theorem~\ref{thm:MCK6} and Corollary~\ref{cor:Y3-6}). Towards the final stage of preparation of our paper, we learned that a different proof of Theorem~\ref{thm:MCK6-intro}(i) has been found by Z. Xu \cite{XuZe-MotivicMultiplicativityCyclicCovers}.
\end{rmk}

\subsection{Applications}

The Franchetta property allows to lift homological relations among generically defined cycles to relations in Chow groups. 
\paragraph{Beauville--Voisin conjecture.}
A typical application in this direction is the so-called Beauville--Voisin conjecture  \cite{Beauville-SplittingBBFiltration, Voisin-ChowRingHK}, which says that for a hyper-K\"ahler variety $X$,  the cycle class map restricted to the sub-$\Q$-algebra of the Chow ring generated by divisors and the Chern classes of the tangent bundle is injective:
\[\cl\colon \langle \CH^1(X), c_i(T_X)\rangle_{\Q\mathrm{-alg}} \hookrightarrow \HH^*(X, \Q). \]

As observed in \cite[Proposition 2.5]{FLV1}, the Franchetta property for a family of hyper-K\"ahler varieties implies the Beauville--Voisin conjecture for a very general member of this family. Therefore, Theorem~\ref{thm:main} (iii) has the following immediate consequence:  
\begin{cor2}\label{cor:strongBV}
The Beauville--Voisin conjecture holds for a very general double EPW cube $X$. More strongly, the sub-$\Q$-algebra of  $\CH^*(X)$ generated by divisors,  Chern classes of  $T_X$ and the class of the fixed locus $W$ of the covering involution on $X$, is injective:
\[\cl\colon \langle \CH^1(X), c_i(T_X), [W] \rangle_{\Q\mathrm{-alg}} \hookrightarrow \HH^*(X, \Q). \]
\end{cor2}

\begin{rmk}
\begin{enumerate}[(i)]
\item The fact that one can include $W$ in Corollary \ref{cor:strongBV} is in line with the conjectures made in \cite{Voisin2016}, because $W$ is a constant cycle subvariety (in particular coisotropic) in $X$. Corollary \ref{cor:strongBV} implies in particular that $  W = {5\over 8}(3 h^3-hc_2)$ in $\CH^3(X)$, since this relation holds in cohomology \cite[Lemma 3.1]{Rizzo}.

\item For double EPW sextics, the Beauville--Voisin conjecture is known by work of Ferretti \cite{Ferretti} in the case of Picard number one and in general by Laterveer \cite{LatBV}.
   \end{enumerate} 
\end{rmk}

\paragraph{A motivic Beauville--Donagi-type relation.}
In our previous work \cite[Theorem 3]{FLV2}, we established a precise relation between the motive of a cubic hypersurface and the motive of its Fano variety of lines, lifting the Beauville--Donagi relation to the motivic level. 

In a similar fashion, we deduce from Theorem \ref{thm:main}(ii) such a motivic relation between a GM fourfold and the associated double EPW sextics:

\begin{cor2}[Corollary \ref{cor:motive-X}]\label{cor:relationXY} 
Let $X$ be the double EPW sextic associated to a GM fourfold~$Y$.
Then there is an isomorphism of (rational) Chow motives
$$ \h(X)\cong\Sym^2\bigl(\h^4_{\van}(Y)(1)\bigr)\oplus\bigoplus_{i=1}^{3}\h^4_{\van}(Y)(2-i)\oplus\bigoplus_{k=0}^4\1(-k),$$
where $\h^4_{\van}(Y)$ is the non-trivial part of the motive of $Y$ as in Definition \ref{defn:CK}. 
\end{cor2}

\paragraph{Orlov conjecture.}
In \cite{Orlov-DerivedCatAndMotives}, Orlov conjectured that derived equivalent smooth projective varieties have isomorphic (rational) Chow motives. By \cite{Kapustka-derived}, dual double EPW sextics are derived equivalent. Using Theorem \ref{thm:main}(ii) we can confirm Orlov's conjecture in this case: 

\begin{cor2}[Theorem~\ref{cor:dual-motives}]
\label{cor:DualdEPW6-intro}
    Let $X$ and $X^\vee$ be dual double EPW sextics. 
    Then there exists an isomorphism of (rational) Chow motives:
    $$ \h(X)\xrightarrow{\sim} \h(X^\vee).$$
\end{cor2}

In particular, Corollary~\ref{cor:relationXY} holds also for the dual double EPW sextic $X^\vee$ associated to~$Y$.
Note that dual double EPW sextics are in general \textit{not} birational, hence Corollary \ref{cor:DualdEPW6-intro} cannot be deduced from \cite{Rie14}.

\subsection{Sketch of proof of main result}
Theorem \ref{thm:main} will be proved via GM fourfolds, by reducing to Theorem \ref{thm:Y3-intro} (see \textbf{Step 1} below), while Theorem~\ref{thm:Y3-intro} (with Theorem~\ref{thm:MCK-intro} being a special case and crucial intermediate step) will be proved following the proof pattern developed in our previous work \cite{FLV-MCK} and \cite{FLV2} (\textbf{Step 2} below). 
Before sketching their proofs, we introduce some notation: $Y$ is a (smooth) GM fourfold, $X$ is the associated double EPW sextic and $C$ is the associated double EPW cube, via the moduli maps
\[
\begin{array}{c@{\qquad}c}
\begin{tikzcd}
& & \MM_{\widetilde{\mathrm{EPW}}} \\
\MM_{\mathrm{GM}} \arrow[r] & \MM_{\mathrm{EPW}} \arrow[ur] \arrow[dr] & \\
& & \MM_{\mathrm{IKKR}}
\end{tikzcd}
\, ;
&
\begin{tikzcd}
& & X = X_{A^\perp} \\
Y \arrow[r, mapsto] & A \arrow[ur, mapsto] \arrow[dr, mapsto] & \\
& & C = C_{A^\perp}
\end{tikzcd}
\end{array}
\]

\noindent\textbf{Step 1 (Reduce Theorem \ref{thm:main} to Theorem~\ref{thm:Y3-intro}).}  The main results are Corollary~\ref{cor:sextic-twists} and Theorem~\ref{thm:relative-transfer} where we show that the motives of $X$ and $C$ are controlled by the motives of $Y^2$ and $Y^3$ respectively, in a way that is generically defined over the parameter space $B$ of GM fourfolds. 

More precisely, the Kuznetsov component $\Ku(Y)$ of the derived category of $Y$, introduced in \cite{KP}, is a non-commutative K3 surface admitting Bridgeland stability conditions \cite{PPZ}. For $Y$ very general, $X$ (resp.~$C$) is realized in \cite{PPZ} (resp.~in \cite{KKM,FGLZ2}) as a moduli space of Bridgeland-stable objects in $\Ku(Y)$ with certain Mukai vector $v$ of square 2 (resp.~certain Mukai vector $w$ of square 4):
\[
X\simeq M_\sigma(\Ku(Y), v), \qquad C\simeq M_{\sigma'}(\Ku(Y), w).
\]
B\"ulles' argument \cite{Bulles}, in the refined form of our previous paper \cite{FLV2}, then exhibits the motives $\h(X)$ and $\h(C)$ as a direct summand of a finite sum of Tate twists of the motive $\h(Y^2)$, resp.\ of $\h(Y^3)$.
In this paper, we further refine B\"ulles' argument to control the Tate twists that appear. The algebraic cycles inducing the motivic embeddings and their left inverses, are generically defined over $B$. Consequently, we have inclusions that are compatible with cycle class maps:
\begin{equation}
\label{eqn:GDCHinclusion}
\GDCH^*_B(X)\subseteq \bigoplus \GDCH^*_B(Y^2), \qquad \GDCH^*_B(C)\subseteq \bigoplus \GDCH^*_B(Y^3),
\end{equation}
and so parts (i) and (iii) of Theorem~\ref{thm:main} would follow from the Franchetta property of~$Y^3$. Now the caveat here is that the Franchetta property does \textit{not} hold for $Y^3$ in general (Remark~\ref{rmk:ConditionSpecialIsNecessary}). However, a subtle but crucial point to make is the following:

\medskip
\noindent\textit{Observation}:
    $X$ and $C$ only depend on the Lagrangian datum $A$ determined by $Y$. Hence the inclusions \eqref{eqn:GDCHinclusion} hold for \textit{any} $Y$ with Lagrangian datum $A$ (such GM fourfolds move in a 4-dimensional family, called \textit{period partners}). Therefore, to prove the Franchetta property for $X$ and $C$, we have the flexibility of \textit{choosing} a suitable $Y$ with Lagrangian datum $A$, for which we can prove the Franchetta property for $Y^3$. The choice that works for our purpose is to take $Y$ to be a \textit{special} GM fourfold.

Since for a generic Lagrangian subspace $A$ in $\bigwedge^3V_6$, there always exists a special GM fourfold $Y$ with this Lagrangian datum, Theorem~\ref{thm:Y3-intro} is sufficient to conclude Theorem~\ref{thm:main}.
\medskip

\noindent\textbf{Step 2 (Cycles on $Y^3$).}  
Following the strategy of \cite{FLV-MCK} and \cite{FLV2}, we show Theorem~\ref{thm:Y3-intro} in the following substeps:

\noindent
    (i) \emph{(Generically defined cycles are tautological).} 
    For an arbitrary GM fourfold $Y$ and $m\in \NN$, the \textit{tautological ring} $\R^*(Y^m)$ is the subalgebra of $\CH^*(Y^m)$  generated by the Pl\"ucker divisor class $h$ and the Chern class $c_2$ of each factor and the partial diagonals. In Proposition~\ref{prop:generation}, we show that generically defined cycles are tautological in the following range:
    \[\GDCH^*_B(Y^m)=\R^*(Y^m)\quad\text{for }m\leq3, \qquad \GDCH^i_B(Y^m)=\R^i(Y^m)\quad\text{for }m\le 6\ \text{and}\ i\leq8 .\]
    
  \noindent
  (ii)  \emph{(Cohomological relation among tautological cycles).} Arguments similar to Tavakol \cite{Tavakol}, Yin \cite{Yin} and \cite{FLV-MCK} allow to classify cohomological relations among tautological cycles. Precisely, for $m\leq 45$, the cohomological relations in $\R^*(Y^m)$ are generated by pullbacks of three relations, one on $Y$, one on $Y^2$ and one on $Y^3$, which together express the multiplicativity of the canonical Chow--K\"unneth decomposition (Proposition~\ref{prop:GM-criteria}). 

    \smallskip\noindent
    (iii) \emph{(Establish relations modulo algebraic equivalence).} The relations on $Y$ and on $Y^2$ hold modulo algebraic equivalence for any $Y$ (Proposition~\ref{prop:alg-Y2}), and the relation on $Y^3$ holds modulo algebraic equivalence if $Y$ is special (Proposition~\ref{prop:MCK-alg}). Consequently, the following cycle class map is injective for $i\leq 8$ (Corollary~\ref{cor:Y4-alg})
    \[\GDCH^i_B(Y^4)/{\sim_{\mathrm{alg}}} \longrightarrow H^{2i}(Y^4, \Q).\]

   \noindent
  (iv) \emph{(Lifting to rational equivalence).}
By the Voisin--Voevodsky nilpotence theorem \cite{Voevodsky,VoisinNilp}, a homological splitting of $\h(Y^{(2)})$ into a direct sum of $\h(X)(-2)$, copies of $\h(Y)$ and Tate motives can therefore be lifted to a splitting in $\Mrat$ induced by generically defined correspondences (Corollary~\ref{cor:sym2-splitting}). The only obstruction to multiplicativity is then a single class $\Upsilon_Y=\tau\cdot p_1^*c_2\in\CH^6(Y\times Y)$. For special $Y$ it is symmetric and homologically trivial, hence comes from a generically defined class on $Y^{(2)}$, whose image under the above splitting consists of classes in $\GDCH^4_B(X)$, in $\GDCH^*_B(Y)$ and in Tate motives. All of these inject into cohomology, the first because $\GDCH^4_B(X)=\QQ\, h_X^4$ (Theorem~\ref{thm:zero-cycles-X}), which we prove by means of the Iliev--Manivel family of conics on the ordinary GM fourfolds, and then extend to special GM fourfolds by a spread argument. Therefore $\Upsilon_Y=0$, and all cohomological relations among tautological cycles on $Y^3$ hold in the Chow ring. Consequently, if $Y$ is a special GM fourfold we get the injectivity of the cycle class map on the tautological ring (Proposition \ref{prop:Y3}):
\begin{equation}
    \cl\colon \R^*(Y^m)\hookrightarrow \HH^*(Y^m, \QQ) \quad\text{for }m\leq 45.
\end{equation}

\subsubsection*{AI disclosure}
All mathematical results in this paper have been obtained by the authors without using~AI. We have used AI for proofreading and to generate polished notes from a preliminary version. During this writing and editing process, Claude Opus 5.5 discovered a subtle point that we had initially overlooked: our method does not allow to prove the stronger version of Theorem \ref{thm:main}(iii) for the universal family after the dominant base change $\MM_{\mathrm{EPW}}\to \MM_{\mathrm{IKKR}}$. This explains the current form of Theorem \ref{thm:main}, and the difference between Theorem \ref{thm:relative-transfer} (for double EPW cubes) and Corollary \ref{cor:sextic-twists} (for double EPW sextics); cf. also 
Remark \ref{rmk:DescendabilityKKMIsomorphim}.

Responsibility for all results in the paper remains with the authors.

\subsubsection*{Acknowledgments} R.L. would like to thank C.M. and C.V. for an invitation to Bielefeld in July 2026, where this project came into being.
L.F.\ is supported by the University of Strasbourg Institute for Advanced Study (USIAS), by the Institut Universitaire de France (IUF), by the Agence Nationale de la Recherche (ANR) under project ANR-24-CE40-4098, by the International Emerging Actions (IEA) project of CNRS.
The research of C.M.\ and C.V.\ was funded by the Deutsche Forschungsgemeinschaft
(DFG, German Research Foundation), Project-ID~491392403, TRR~358.

\subsubsection*{Conventions}
We work over the base field of complex numbers. All Chow groups $\CH^*(-)$ are with rational coefficients. $\CH^*_{\mathrm{hom}}(-)$ stands for the subgroup of homologically trivial cycles. $\B^*(-)$ is the group of algebraic cycles modulo \textit{algebraic equivalence}, again with rational coefficients. $\Mrat$ and $\Mhom$ are the categories of pure motives with rational coefficients modulo rational and homological equivalence respectively; see \cite{Andre}. We write $\h(-)$ for the contravariant functor that associates to a smooth projective variety $X$ its motive $\h(X)$, and $\1$ for the motive of a point.

\section{Preliminaries}

\subsection{Generically defined cycles and the Franchetta property}
We first recall the following notions central to the subject of the paper; see \cite{FLV1, FLV2}.
\begin{defn}\label{defn:gdch}
Let $\pi\colon\Zz\to S$ be a smooth projective morphism between smooth irreducible varieties (or more generally Deligne--Mumford stacks), and let $Z=\Zz_s$ be the fibre over a closed point $s\in S$. The subgroups of \emph{generically defined cycles} on $Z$ with respect to the family $\Zz\to S$ are
\[
\GDCH^i_S(Z)\coloneqq\ima\bigl(\CH^i(\Zz)\to\CH^i(Z)\bigr),
\qquad
\GDB^i_S(Z)\coloneqq\ima\bigl(\B^i(\Zz)\to\B^i(Z)\bigr),
\]
the maps being Gysin restriction to the fibre. We say that $\Zz\to S$ has the \emph{Franchetta property in codimension $i$} if the cycle class map $\cl\colon\GDCH^i_S(Z)\to\HH^{2i}(Z,\QQ)$ is injective for every closed point $s\in S$, and that it has the \emph{Franchetta property} if this holds in every codimension. Replacing $\GDCH$ by $\GDB$ throughout gives the Franchetta property \emph{modulo algebraic equivalence}.
\end{defn}

Clearly, $\GDCH^*_S(Z)$ is a $\Q$-subalgebra of $\CH^*(Z)$ containing all Chern classes of the tangent bundle of $Z$. We collect in the following two lemmas some basic behaviors of the Franchetta property. 

\begin{lem}[Spread]\label{lem:spread}
Let $\Zz\to S$ be as in Definition~\ref{defn:gdch} and let $\Gamma,\Gamma_1,\ldots,\Gamma_N\in\CH^i(\Zz)$.
\begin{enumerate}[(i)]
\item If $\Gamma|_{Z_s}=0$ in $\CH^i(Z_s)$ for the very general $s\in S$, then $\Gamma|_{Z_s}=0$ for every $s\in S$.
\item If $\Gamma|_{Z_s}$ lies in the $\QQ$-subspace generated by $\Gamma_1|_{Z_s},\ldots,\Gamma_N|_{Z_s}$ for the very general $s\in S$, then there are $\lambda_1,\ldots,\lambda_N\in\QQ$ such that $\Gamma|_{Z_s}=\sum_j\lambda_j\,\Gamma_j|_{Z_s}$ in $\CH^i(Z_s)$ for every $s\in S$.
\end{enumerate}
The same statements hold with $\CH^i$ replaced by $\B^i$.
\end{lem}

\begin{proof}
(i) is \cite[Lemma~3.1]{VoisinBook}; a proof can be found in \cite[Proposition~1.4]{VoisinUnir}. The key point is that a Hilbert scheme argument shows that the locus $S_\Gamma\coloneqq\{s\in S\mid \Gamma|_{Z_s}=0\}$ is a countable union of closed algebraic subsets of $S$. It must equal to $S$ if it contains a very general point. 

For (ii), put $\Gamma_\lambda\coloneqq\Gamma-\sum_j\lambda_j\Gamma_j$ for $\lambda\in\QQ^N$. By assumption the union $\bigcup_{\lambda\in\QQ^N}S_{\Gamma_\lambda}$ contains the very general point of $S$. It is a countable union of closed algebraic subsets,
so one of these closed subsets equals $S$; hence $S_{\Gamma_\lambda}=S$ for some $\lambda$ by (i). 

The proof modulo algebraic equivalence is identical, algebraic equivalence being likewise witnessed by countably many bounded families.
\end{proof}

For $\Gamma\in\CH^i(\Zz)$, the cohomology classes $[\Gamma|_{Z_s}]$ form a flat section of the local system $R^{2i}\pi_*\QQ$ on $S$; since $S$ is connected, $\Gamma|_{Z_s}$ is homologically trivial for one $s$ if and only if it is so for every $s$. Together with Lemma~\ref{lem:spread}(i) this shows that $\Zz\to S$ has the Franchetta property as soon as the cycle class map is injective on $\GDCH^*_S(Z_s)$ for the very general~$s\in S$. Recall the following basic result from \cite[Remark 2.6]{FLV1}:

\begin{lem}[Base change]\label{lem:base-change}
Let $\Zz\to S$ be as in Definition~\ref{defn:gdch} and let $Z=\Zz_s$.
\begin{enumerate}[(i)]
\item If $U\subseteq S$ is a dense open subset containing $s$, then $\GDCH^*_U(Z)=\GDCH^*_S(Z)$ and $\GDB^*_U(Z)=\GDB^*_S(Z)$.
\item Let $S'\to S$ be a morphism of smooth varieties (or stacks), 
let $\Zz'\coloneqq\Zz\times_SS'$ and let $s'\in S'$ be a preimage of $s$. Then, under the identification $\Zz'_{s'}=Z$,
\[
\GDCH^*_S(Z)\subseteq\GDCH^*_{S'}(Z),\qquad \GDB^*_S(Z)\subseteq\GDB^*_{S'}(Z).
\]
In particular, if $S'\to S$ is dominant and $\Zz'\to S'$ has the Franchetta property, then so does $\Zz\to S$.
\end{enumerate}
\end{lem}

\subsection{Gushel--Mukai fourfolds}\label{SS:GM4}

Let $V_5$ be a five-dimensional complex vector space, let $\Gr\coloneqq\Gr(2,V_5)\subseteq\PP(\bigwedge^2V_5)=\PP^9$ and let $\CGr\coloneqq\CGr(2,V_5)\subseteq\PP(\bigwedge^2V_5\oplus\C)=\PP^{10}$ be the projective cone over $\Gr$, with vertex $O$. Following \cite{DK1} (see also \cite{Debarre}), a \emph{Gushel--Mukai (GM) fourfold} is a smooth  transversal intersection
\[
Y=\CGr\cap H_1\cap H_2\cap Q\subseteq\PP^{10},
\]
where $H_1,H_2$ are hyperplanes and $Q$ is a quadric. Since $Y$ is smooth, $O\notin Y$. The linear projection away from $O$ gives rise to the \emph{Gushel morphism} $$\gamma\colon Y\to\Gr,$$ and we set
\[
h\coloneqq c_1(\OO_Y(1))=\gamma^*\sigma_1\in\CH^1(Y),\qquad
c_2\coloneqq c_2(\gamma^*\UU)=\gamma^*\sigma_{1,1}\in\CH^2(Y),
\]
where $\UU$ is the tautological rank-2 subbundle on $\Gr$. The GM fourfold $Y$ is called \emph{ordinary} if $O\notin H_1\cap H_2$ and \emph{special} otherwise. 
\begin{itemize}
    \item If $Y$ is ordinary, then the Gushel morphism $\gamma$ is a closed embedding and $\CGr\cap H_1\cap H_2$ is mapped isomorphically onto $Y'\coloneqq \Gr\cap H$ for some hyperplane $H\subseteq\PP(\bigwedge^2V_5)$ and $Y\simeq Y'\cap Q$. Note that $Y'$ is a degree-5 del Pezzo Fano fivefold.  
    \item  If $Y$ is special, then the Gushel morphism $\gamma$ factors as a double cover $p\colon Y\to Y_0\coloneqq\Gr\cap H_1'\cap H_2'$ onto a fourfold $Y_0$ (smooth by \cite[Proposition~2.22]{DK1}), branched along an ordinary GM threefold; we write $\iota_Y$ for the covering involution. Note that $Y_0$ is a degree-5 del Pezzo Fano fourfold. 
\end{itemize}

We first recall that the del Pezzo manifolds appearing above have "trivial" Chow group (see \cite[Definition 4.15]{VoisinBook} and \cite[Lemma~2.1]{FLV2}):
\begin{lem}[del Pezzo manifolds]
\label{lem:trivial-chow}
Let $H_1',H_2'\subseteq\PP(\bigwedge^2V_5)$ be hyperplanes.
\begin{enumerate}[(i)]
    \item If $Y'\coloneqq\Gr\cap H_1'$ is smooth of dimension $5$, then $Y'$ is of Tate Chow motive. 
    \item If $Y_0\coloneqq\Gr\cap H_1'\cap H_2'$ is smooth of dimension $4$, then $Y_0$ is of Tate Chow motive. 
\end{enumerate}
In particular, the cycle class maps for $Y'$ and $Y_0$ are isomorphisms.
Moreover, their Chow rings are generated by $h$ and $c_2$.
\end{lem}

\begin{proof}
The fact that Chow motives of $Y'$ and $Y_0$ are Tate is proved in \cite[Proposition~4.6]{LatPG}. Alternatively, the derived categories of $Y'$ and of $Y_0$ admit full exceptional collections by \cite[Section~6.1]{Kuznetsov}, so their Chow motives are Tate. Applying the realization functor, we see that the cycle class map is an isomorphism.

For the claim about ring generators, it suffices to check it for cohomology. Since $\HH^*(\Gr,\QQ)$ is generated by the Chern classes of the rank-2 tautological bundle, we only need to show that $\HH^*(\Gr,\QQ)\to\HH^*(Y',\QQ)$ and $\HH^*(\Gr,\QQ)\to\HH^*(Y_0,\QQ)$ are surjective. By the Lefschetz hyperplane and hard Lefschetz theorems, it remains to check the surjectivity for the middle degree. For $Y'$, as it has Tate motive, $\HH^5(Y',\QQ)=0$. For $Y_0$, $\HH^4(Y_0,\QQ)$ has dimension $2$, spanned by $c_2$ and $h^2$ (they are linearly independent since $\deg (h^4)=5$, $\deg (h^2c_2)=2$ and $\deg (c_2^2)=1$). 
\end{proof}

To parameterize GM fourfolds, we first consider
\[
\bar B\coloneqq\PP\HH^0\bigl(\CGr,\OO(1)^{\oplus2}\oplus\OO(2)\bigr).
\]
The \textit{parameter space} of GM fourfolds is defined to be the dense open subset $B\subseteq\bar B$ parameterizing those $b\in \bar{B}$ such that $Y_b$ is a smooth  fourfold. We denote by $\YY\to B$ the universal GM fourfold and by $\YY^m_{/B}\coloneqq\YY\times_B\cdots\times_B\YY$ its relative $m$-th power, with fibre $Y_b^m$ over a point $b\in B$.

\section{Motives of moduli spaces of Bridgeland-stable objects}
\label{S:motive-moduli}

Moduli spaces of Bridgeland-stable objects in 2-Calabi--Yau categories provide a rich source of projective hyper-K\"ahler varieties.
B\"ulles \cite{Bulles} showed that the Chow motive of a moduli space of stable sheaves on a K3 surface $S$ lies in the thick tensor subcategory of~$\Mrat$ generated by the motive $\h(S)$ of the surface.
Floccari--Fu--Zhang \cite{FFZ} extended his argument to the non-commutative setting, and \cite{FLV2} refined it so as to bound the precise \textit{power} of the underlying variety.
In this section, we refine the statement in the setting of moduli spaces of stable objects in the Kuznetsov components of GM fourfolds in order to control the \textit{Tate twists} and we specialize it to the double EPW sextics and double EPW cubes associated with GM fourfolds. The control on the Tate twists will be essential to our proof of Theorem \ref{thm:main}(ii).

\subsection{Refining B\"ulles' argument and GM fourfolds}\label{SS:Bulles}

Let $Y$ be a smooth projective variety and let $\mathcal A\subseteq\D^b(Y)$ be an admissible subcategory which is \emph{2-Calabi--Yau}, that is, whose Serre functor is the double shift $[2]$. 
Write $\Stab^\dagger(\mathcal A)$ for the distinguished connected component of the space of Bridgeland stability conditions \cite{Bridgeland} and assume that this is non-empty.
Equipping the topological K-theory of $Y$ with the Euler pairing $(\kappa,\kappa')\coloneqq-\chi(\kappa,\kappa')$, the \emph{Mukai lattice} of $\mathcal A$ is the torsion-free quotient of
\[
\widetilde{\HH}(\mathcal A)\coloneqq\{\kappa\in K_{\top}(Y)\ \vert\ (\kappa,[E])=0\ \text{ for all }E\in{}^{\perp}\mathcal A\},
\]
and the \emph{Mukai vector} of an object $E$ of $\mathcal A$ is $\ch(E)\cdot\sqrt{\mathrm{td}_Y}$, see \cite{AT14} or \cite[\S5]{FFZ}.
For an element $\lambda\in \widetilde{\HH}(\mathcal A)$,
call $\sigma\in\Stab^\dagger(\mathcal A)$ \emph{$\lambda$-generic} if stability and semistability coincide for objects of Mukai vector $\lambda$.
The moduli space of $\sigma$-stable objects of $\mathcal A$ with Mukai vector $\lambda$ is denoted by $M_\sigma(\mathcal A,\lambda)$.
Then, \cite[Thm. 0.1]{Bulles}, which was extended to the setting of 2-Calabi--Yau categories in \cite{FFZ}, exhibits $\h(M_\sigma(\mathcal A,\lambda))$ as a direct summand of a motive of the form $\bigoplus_{i}\h(Y^{r_i})(l_i)$ for some $r_i$ and $l_i$.

This was refined to control the powers $r_i$ that appear under the additional assumption that $H^{0,1}(M_\sigma(\mathcal A,\lambda))=0$ in \cite{FLV2}.
In that case, it was shown that $\h(M_\sigma(\mathcal{A},\lambda))$ is in fact a direct summand of a motive of the form $\bigoplus_i\h(Y^m)(l_i)$ for some $l_i$, where $m=\frac{1}{2}\dim M_\sigma(\mathcal A,\lambda)$.

We now extend this result as in the following theorem to control the Tate twists $l_i$ that appear, under the additional assumption of the triviality of the Chow groups of $Y$ in a certain range.

\begin{thm}\label{thm:twists}
    Let $Y$ be a smooth projective variety of dimension $2d$ such that rational and numerical equivalence agree in codimension $i\leq d$.
    Let $\mathcal{A}$ be a 2-Calabi--Yau admissible subcategory of $\D^b(Y)$.
    Let $\lambda\in \widetilde{\HH}(\mathcal A)$ be primitive, let $\sigma\in\Stab^\dagger(\mathcal A)$ be $\lambda$-generic, and assume that $M\coloneqq M_\sigma(\mathcal A,\lambda)$ is smooth projective of dimension $2m\coloneqq \lambda^2+2$ and satisfies $H^1(M,\OO_M)=0$.
    Then $\h(M)$ is a direct summand of $\h(Y^m)(m(d-1))^{\oplus r}$ for some $r$.
\end{thm}
The proof will be given in Section~\ref{SS:twists}.

\begin{rmk}
    Let us briefly discuss the assumption that rational and numerical equivalence agree in codimension $i\leq d$. 
The assumption will be used in the more explicit form of Lemma~\ref{lem:codim}, which is equivalent to it. 
In case the cycle class map is injective in codimension $i\leq d$, i.e. that rational and homological equivalence agree, the Lefschetz standard conjecture for $Y$ implies that homological and numerical equivalence agree.
More concretely, one can use the Bloch--Srinivas argument of \cite[Lemma~5.1]{FV2} to prove that it is sufficient to assume that:
\begin{itemize}
    \item $\CH_j(Y)$ is supported on a closed subvariety of dimension $j+1$ for $j\leq d-2$;
    \item and $H^{2k-1}(Y,\QQ)=0$ for all $2k-1<2d$. 
\end{itemize}
In particular, both of these assumptions apply to the following $Y$: K3 surfaces, cubic fourfolds, GM fourfolds \cite{DK2}, GM sixfolds \cite{DK2,BL}, and Debarre--Voisin 20-folds \cite{DebarreVoisin,Lat-DV}.
\end{rmk}

Let us now specialise to GM fourfolds $Y$.
Kuznetsov--Perry \cite{KP} established the semi-orthogonal decomposition
\[
\D^b(Y)=\langle\Ku(Y),\OO_Y,\gamma^*\UU^\vee,\OO_Y(1),\gamma^*\UU^\vee(1)\rangle,
\]
and showed that the Kuznetsov component $\Ku(Y)$ is 2-Calabi--Yau.
Perry--Pertusi--Zhao \cite{PPZ} constructed stability conditions on $\Ku(Y)$ and proved that $M_\sigma(\Ku(Y),\lambda)$, for $\lambda$ a primitive Mukai vector and $\sigma$ a $\lambda$-generic stability condition, is either empty or a smooth projective hyper-K\"ahler variety of dimension $\lambda^2+2$ \cite[Theorem~1.5]{PPZ}.

The numerical Grothendieck group of $\Ku(Y)$ contains a canonical rank-2 sublattice $\langle\Lambda_1,\Lambda_2\rangle$, namely the projection to $\Ku(Y)$ of the image of the pullback $K(\Gr)\to K(Y)$, with $\Lambda_1^2=\Lambda_2^2=2$ and $(\Lambda_1,\Lambda_2)=0$; see \cite{KP} and \cite[\S5.4]{PPZ}.
For each pair of coprime integers $(a,b)$, the moduli spaces $M_\sigma(\Ku(Y),a\Lambda_1+b\Lambda_2)$ form a smooth projective family over a dense open subset of the moduli space of ordinary GM fourfolds \cite[Theorem~5.15 and proof of Theorem~1.7]{PPZ}. 
The two cases we use are $(a,b)=(1,0)$, of dimension $4$, and $(a,b)=(1,-1)$, of dimension $6$; as in the introduction, we write $v\coloneqq\Lambda_1$ and $w\coloneqq\Lambda_1-\Lambda_2$, of squares $v^2=2$ and $w^2=4$, of divisibilities $\operatorname{div}(v)=1$ and $\operatorname{div}(w)=2$.

For a very general GM fourfold $Y$, Guo--Liu--Zhang \cite[Theorem~1.1]{GLZ} proved that $M_\sigma(\Ku(Y),v)$ is isomorphic to the double EPW sextic $X$ associated with $Y$.
In fact, $M_\sigma(\Ku(Y),v)$ is isomorphic to the double EPW sextic $X$ associated with $Y$ for \emph{general}~$Y$. 
The assumption that $Y$ is very general is only used in \cite{GLZ} to prove the stability of certain projection objects.
The extension of this result to general $Y$ is described in \cite[\S A.3]{FGLZ2}.
In particular, this means that the family of such moduli spaces and the family of double EPW sextics are isomorphic over the dense open subset
of $B$ where both are defined and smooth. 
Hence, we obtain the following corollary.
The construction of the quasi-universal object and hence of all the correspondences involved in the proof of Theorem~\ref{thm:twists} can be done over a family.

\begin{cor}[Double EPW sextics]
\label{cor:sextic-twists}
    Let $B^\circ\subseteq B$ be the open subset over which the associated double EPW sextic is smooth.
    Let $\XX\to B^\circ$ be the universal family of double EPW sextics.
    Then there exists a relative correspondence over $B^\circ$ inducing, for every $b\in B^\circ$, a split injection
    $$\h(X_b)\hookrightarrow \h(Y_b^2)(2)^{\oplus r}\ \ \hbox{in}\ \Mrat,$$
    whose left inverse is also induced by a relative correspondence over $B^\circ$.
\end{cor}

\begin{rmk}
    \label{rmk:DescendabilityKKMIsomorphim}
    Unfortunately, the same does not apply in the case of double EPW cubes.
Kapustka--Kapustka--Mongardi \cite[Theorem~1 and Corollary~1.3]{KKM} produced a Hodge isometry between the primitive second cohomology of the associated double EPW cube C and that of $M_\sigma(\Ku(Y),w)$.
Hence, by the Torelli theorem, they are isomorphic for very general $Y$.
While this can be extended to the existence of an isomorphism for general $Y$, meaning that Theorem~\ref{thm:twists} applies, it is not clear that this isomorphism extends to an isomorphism over an open subset of $B$ and hence there is a priori no relative correspondence inducing the split injection.
This will be addressed in Section~\ref{SS:Ku-GM}.
\end{rmk}

\subsection{Proof of Theorem \ref{thm:twists}}\label{SS:twists}

The proof consists of an analysis of the terms that appear in the proof of \cite[Theorem~0.1]{Bulles} in order to bound the twists that appear.
We remark that the following proof in fact recovers the refinement \cite[Theorem~1.1]{FLV2}.
We use the following version of the assumption on $Y$ that rational and numerical equivalence agree in codimension $i\leq d$.

\begin{lem}\label{lem:codim}
    Let $Y$ be a smooth projective variety of dimension $2d$. 
    Then rational and numerical equivalence agree in codimension $i\leq d$ if and only if $\h(Y)$ admits a decomposition 
    $$\h(Y)=\h_{\mathrm{tr}}(Y)\oplus\bigoplus_i\mathds{1}(-i)^{\oplus \rho_i}, $$
    where $\h_{\mathrm{tr}}(Y)$ satisfies $\CH^i(\h_{\mathrm{tr}}(Y))=0$ for all $i\leq d$, or equivalently, that $\h_{\mathrm{tr}}(Y)$ is cut out by a symmetric projector $\pi_{\mathrm{tr}}$ satisfying $\pi_{\mathrm{tr}}\circ\CH^i(V\times Y)=0$ for all smooth projective $V$ and all $i\leq d$.
\end{lem}
\begin{proof}
    Assume that rational and numerical equivalence agree in codimension $i\leq d$.
    Then $\CH^i(Y)$ is finite-dimensional for $i\leq d$ and the intersection pairing $\CH^i(Y)\times\CH^{2d-i}(Y)\to\QQ$ has trivial left kernel.
    Choose a basis $z_j^i$ of $\CH^i(Y)$ and choose $w_j^{2d-i}\in\CH^{2d-i}(Y)$ satisfying $\deg(z_j^i\cdot w_l^{2d-k})=\delta_{jl}\delta_{ik}$ for $i,k\leq d$.
    Then, define
    $$\pi_{\mathrm{tr}}\coloneqq \Delta_Y-\sum_{i\leq d}\sum_j w_j^{2d-i}\times z_j^i-\sum_{i<d}\sum_j z_j^i\times w_j^{2d-i}$$
    and denote the motive it cuts out by $\h_{\mathrm{tr}}(Y)$.
    By the choice of the $w_j^{2d-i}$, the subtracted summands are mutually orthogonal projectors, and the summand for $i=d$ is symmetric, as the intersection form on $\CH^d(Y)$ is; hence $\pi_{\mathrm{tr}}$ is a symmetric projector with $(\pi_{\mathrm{tr}})_\ast\CH^i(Y)=0$ for $i\leq d$.

    We now prove that $\pi_{\mathrm{tr}}\circ\Gamma=0$ for every $\Gamma\in\CH^i(V\times Y)$ for $i\leq d$.
    We argue by induction on $i$. 
    Assume that $V$ is connected and denote its generic point by $\eta$.
    Choose a countable field $k\subseteq \mathbb{C}$ over which $Y, V,\Gamma,$ and $\pi_{\mathrm{tr}}$ are defined. 
    Then, $\pi_{\mathrm{tr}}\circ\Gamma|_\eta=(\pi_{\mathrm{tr}})_\ast(\Gamma_\eta)$.
    Choose an embedding $k(V)\hookrightarrow \mathbb{C}$, and base change along this map, which is injective on Chow groups, so that $(\pi_{\mathrm{tr}})_\ast(\Gamma_\eta)=0$, as $(\pi_{\mathrm{tr}})_\ast\CH^i(Y)=0$.
    Hence, by localisation, $\pi_{\mathrm{tr}}\circ\Gamma$ is supported on $D\times Y$ for a proper closed $D\subseteq V$.
    Choosing a resolution of singularities $\tilde D$ of $D$, we see that $\pi_{\mathrm{tr}}\circ\Gamma=j_\ast(\pi_{\mathrm{tr}}\circ\delta),$ where $j\colon \tilde D\to V$ and $\delta\in\CH^{i-1}(\tilde D\times Y)$.
    Hence, we are done by induction. 
    The converse follows directly using $V=pt$.
\end{proof}

Assume now that $Y$ satisfies the conditions of Theorem~\ref{thm:twists}.
Moduli spaces of stable objects in admissible $2$-Calabi--Yau subcategories need not admit a universal family in general. 

What does exist is a \emph{quasi-universal} object $\EE_\rho$ on $M\times Y$ of some similitude $\rho$, i.e. an object whose restriction to $\{t\}\times Y$ is $E_t^{\oplus\rho}$ for every $t\in M$, where $E_t$ is the object parametrized by $t$: the moduli stack of $\sigma$-stable objects is a $\mathbb G_m$-gerbe over $M$ \cite[Theorem~21.24]{BLMNPS}, and Mukai's construction \cite[Theorem~A.5]{Mukai} applies verbatim.
Attached to $\EE_\rho$ there are a Brauer--Severi scheme $p\colon N\to M$, that is, a smooth projective morphism which is \'etale-locally a projective bundle of relative dimension $\rho-1$, a line bundle $\OO_N(\rho)$ restricting to $\OO(\rho)$ on the fibres of $p$, and an object $\tilde\EE$ of $\D^b(N\times Y)$ whose restriction to $\{x\}\times Y$ is $E_{p(x)}$ for every $x\in N$; see \cite[\S3]{Markman}.
The proof of Markman's formula only uses the fibres of $\tilde\EE$ and Serre duality in $\mathcal A$, so it applies on $N$ and gives
\[
c_{2m}(\tilde\PPP)=(p\times p)^*\Delta_M\ \ \text{ in }\ \CH^{2m}(N\times N),
\qquad
\tilde\PPP\coloneqq-Rp_{13*}\bigl(p_{12}^*\tilde\EE^\vee\otimes^{\LLL}p_{23}^*\tilde\EE\bigr);
\]
see \cite[\S3]{Markman}, \cite[Proposition~24]{Markman2007}, \cite{MarianZhao} and \cite[Proposition~5.2]{FFZ}.
Put $\xi\coloneqq\tfrac1\rho c_1(\OO_N(\rho))\in\CH^1(N)$, so that $p_*(\xi^{\rho-1})=1$.
Multiplying by $\xi^{\rho-1}\times\xi^{\rho-1}$ and pushing forward along $p\times p$ gives 
\begin{equation}\label{eq:NCK3-diagonal}
\Delta_M=(p\times p)_\ast(c_{2m}(\tilde\PPP)\cdot(\xi^{\rho-1}\times\xi^{\rho-1})) \ \text{ in }\ \CH^{2m}(M\times M).
\end{equation}
Let $\alpha\coloneqq \ch(\tilde\EE^\vee)\cdot p^\ast_Y\sqrt{\mathrm{td}_Y}$ and $\beta\coloneqq {}^t(\ch(\tilde\EE)\cdot p^\ast_Y\sqrt{\mathrm{td}_Y})$.
Grothendieck--Riemann--Roch applied to $p_{13}\colon N\times Y\times N\to N\times N$ yields $\ch(\tilde\PPP)=\beta\circ\alpha$.
Decompose $\alpha$ and $\beta$ by codimension into $\alpha^i\in\CH^i(N\times Y)$ and $\beta^j\in\CH^j(Y\times N)$, so that $\ch_\mu(\tilde\PPP)=\sum_{i+j=\mu+2d}\beta^j\circ\alpha^i$.
We use the following identity.

\begin{lem}\label{lem:logc}
Let $Z$ be a smooth projective variety and let $E\in\D^b(Z)$.
Then
\[
\log c(E)=\sum_{\mu\geq1}(-1)^{\mu-1}(\mu-1)!\,\ch_\mu(E),
\]
where $c(E)$ is the total Chern class and $\ch(E)$ is the Chern character of $E$.
\end{lem}

\begin{proof}
First, assume that $E$ is a vector bundle.
By the splitting principle, we may assume that $c(E)=\prod_{i=1}^{\rk E}(1+x_i)$, where the $x_i$ are the Chern roots of $E$.
Then $\ch_\mu(E)=\frac{1}{\mu!}\sum_ix_i^\mu$, and so
\[
\log c(E)=\sum_i\log(1+x_i)=\sum_{\mu\geq1}\frac{(-1)^{\mu-1}}{\mu}\sum_ix_i^\mu=\sum_{\mu\geq1}(-1)^{\mu-1}(\mu-1)!\,\ch_\mu(E).
\]
Since both the total Chern class and the Chern character factor through $K_0(Z)$, and every class in $K_0(Z)$ is a difference of classes of vector bundles, the identity holds for every $E\in\D^b(Z)$.
\end{proof}

Applying $\exp$ to the identity of Lemma~\ref{lem:logc}, we obtain
\[
c(\tilde\PPP)=\exp\bigl(c_1(\tilde\PPP)\bigr)\cdot\exp\Bigl(\sum_{\mu\geq2}(-1)^{\mu-1}(\mu-1)!\,\ch_\mu(\tilde\PPP)\Bigr).
\]
Hence $c_{2m}(\tilde\PPP)$ is a linear combination of monomials $c_1(\tilde\PPP)^{\tilde e}\ch_{\mu_1}(\tilde\PPP)\cdots\ch_{\mu_k}(\tilde\PPP)$ with $\mu_l\geq2$ and $\tilde e+\sum_l\mu_l=2m$; in particular, the number $k$ of factors of type $\ch_\mu(\tilde\PPP)$, counted with multiplicity, is at most $m$.
Since $Rp_\ast\OO_N=\OO_M$, the irregularity of $N$ vanishes, and we have $c_1(\tilde\PPP)=p_1^*\eta_1+p_2^*\eta_2$ for some $\eta_1,\eta_2\in\CH^1(N)$.
Lieberman's identities, as used in the proof of \cite[Theorem~0.1]{Bulles}, show that
\[
\prod_{l=1}^k\Gamma'_l\circ\Gamma_l={}^t\delta_N\circ(\Gamma'_1\otimes\cdots\otimes\Gamma'_k)\circ(\Gamma_1\otimes\cdots\otimes\Gamma_k)\circ\delta_N
\]
for $\Gamma_l\in\CH^*(N\times Y)$ and $\Gamma'_l\in\CH^*(Y\times N)$, where $\delta_N\in\CH^*(N\times N^k)$ is the correspondence induced by the diagonal embedding $N\hookrightarrow N^k$.
Combining this and expanding further shows that, after absorbing rational coefficients, $(p\times p)^*\Delta_M=c_{2m}(\tilde\PPP)=\sum_I\beta'_I\circ\alpha'_I$, where
\[
\alpha'_I\coloneqq p_N^*\eta_1^e\cdot\bigl((\alpha^{i_1}\otimes\cdots\otimes\alpha^{i_k})\circ\delta_N\bigr)\in\CH^{n_I}(N\times Y^k),
\qquad
n_I\coloneqq e+\textstyle\sum_li_l,
\]
and
\[
\beta'_I\coloneqq p_N^*\eta_2^f\cdot\bigl({}^t\delta_N\circ(\beta^{j_1}\otimes\cdots\otimes\beta^{j_k})\bigr)\in\CH^{2m+2dk-n_I}(Y^k\times N).
\]
A term $I$ comes from a monomial $c_1(\tilde\PPP)^{\tilde e}\ch_{\mu_1}(\tilde\PPP)\cdots\ch_{\mu_k}(\tilde\PPP)$ with $\tilde e+\sum_l\mu_l=2m$.
By splitting $c_1(\tilde\PPP)$ into $p_1^*\eta_1+p_2^*\eta_2$, we split $\tilde e=e+f$, and the indices above satisfy $i_l+j_l=\mu_l+2d$.
Pushing forward along $p$ as in \eqref{eq:NCK3-diagonal} we obtain $\Delta_M=\sum_I\beta_I\circ\alpha_I$, where
\begin{align*}
\alpha_I&\coloneqq(p\times\id)_*\bigl(\alpha'_I\cdot(\xi^{\rho-1}\times1)\bigr)\in\CH^{n_I}(M\times Y^k),\\
\beta_I&\coloneqq(\id\times p)_*\bigl(\beta'_I\cdot(1\times\xi^{\rho-1})\bigr)\in\CH^{2m+2dk-n_I}(Y^k\times M).
\end{align*}
As $\alpha_I\colon\h(M)\to\h(Y^k)(n_I-2m)$ and $\beta_I\colon\h(Y^k)(n_I-2m)\to\h(M)$, this exhibits $\h(M)$ as a direct summand of $\bigoplus_I\h(Y^{k_I})(n_I-2m)$, where $k_I$ is the number of factors of type $\ch_\mu(\tilde\PPP)$ of the term $I$.

\medskip

Now decompose the motive of $Y$ as in Lemma~\ref{lem:codim}.
Notice that a term with $k$ factors of type $\ch_\mu(\tilde\PPP)$ factors through a product of at most $k$ copies of $\h_{\mathrm{tr}}(Y)$ and Lefschetz motives.
Denote the summands of the decomposition of Lemma~\ref{lem:codim} by $\h_s$ and the corresponding projectors by $\pi_s$.
Every term $\beta_I\circ\alpha_I$ is the sum of the terms
\[
(\beta_I\circ\alpha_I)_{s_1,\ldots,s_k}\coloneqq\beta_I\circ(\pi_{s_1}\otimes\cdots\otimes\pi_{s_k})\circ\alpha_I,
\]
which factor through $\h_{s_1}\otimes\cdots\otimes\h_{s_k}(n_I-2m)$.
Putting
\[
(\beta'_I\circ\alpha'_I)_{s_1,\ldots,s_k}\coloneqq\beta'_I\circ(\pi_{s_1}\otimes\cdots\otimes\pi_{s_k})\circ\alpha'_I,
\]
the projection formula gives, as for \eqref{eq:NCK3-diagonal},
\[
(\beta_I\circ\alpha_I)_{s_1,\ldots,s_k}=(p\times p)_*\bigl((\beta'_I\circ\alpha'_I)_{s_1,\ldots,s_k}\cdot(\xi^{\rho-1}\times\xi^{\rho-1})\bigr),
\]
so that $(\beta_I\circ\alpha_I)_{s_1,\ldots,s_k}$ vanishes as soon as $(\beta'_I\circ\alpha'_I)_{s_1,\ldots,s_k}$ does.
We now show that many of these terms vanish.

\begin{lem}\label{lem:vanishing}
A term $(\beta'_I\circ\alpha'_I)_{s_1,\ldots,s_k}$ vanishes unless the indices are in the following ranges, the conditions applying to every $\h_{s_l}$.
\begin{enumerate}[(i)]
\item If $\h_{s_l}=\1(-a)$, the term vanishes unless $2d-j_l\leq a\leq i_l$.
\item If $\h_{s_l}=\h_{\mathrm{tr}}(Y)$, the term vanishes unless $i_l\geq d+1$ and $j_l\geq d+1$.
\end{enumerate}
\end{lem}

\begin{proof}
The term equals
\[
p_2^*\eta_2^f\cdot p_1^*\eta_1^e\cdot\Bigl({}^t\delta_N\circ\bigl((\beta^{j_1}\circ\pi_{s_1}\circ\alpha^{i_1})\otimes\cdots\otimes(\beta^{j_k}\circ\pi_{s_k}\circ\alpha^{i_k})\bigr)\circ\delta_N\Bigr),
\]
so it vanishes as soon as $\pi_{s_l}\circ\alpha^{i_l}=0$ or $\beta^{j_l}\circ\pi_{s_l}=0$ for some $l$.

(i) Suppose that $\h_{s_l}=\1(-a)$ is cut out by the projector $\pi_{s_l}=x\times x'$, with $x\in\CH^{2d-a}(Y)$ and $x'\in\CH^a(Y)$.
Then $\pi_{s_l}\circ\alpha^{i_l}=p_{N*}(\alpha^{i_l}\cdot p_Y^*x)\times x'$, and $p_{N*}(\alpha^{i_l}\cdot p_Y^*x)\in\CH^{i_l-a}(N)$ vanishes for codimension reasons unless $a\leq i_l$.
An analogous calculation, $\beta^{j_l}\circ\pi_{s_l}=x\times p_{N*}(p_Y^*x'\cdot\beta^{j_l})$ with $p_{N*}(p_Y^*x'\cdot\beta^{j_l})\in\CH^{j_l+a-2d}(N)$, shows the other inequality.

(ii) By Lemma~\ref{lem:codim} applied with $V=N$, we have $\pi_{\mathrm{tr}}\circ\alpha^{i_l}=0$ if $i_l\leq d$, and, $\pi_{\mathrm{tr}}$ being symmetric, $\beta^{j_l}\circ\pi_{\mathrm{tr}}={}^t\bigl(\pi_{\mathrm{tr}}\circ{}^t\beta^{j_l}\bigr)=0$ if $j_l\leq d$.
\end{proof}

Using Lemma~\ref{lem:vanishing}, we can now control the twists that appear.

\begin{proof}[Proof of Theorem~\ref{thm:twists}]
First, as observed above, the number $k$ of factors of type $\ch_\mu(\tilde\PPP)$ in the terms $\beta'_I\circ\alpha'_I$, counted with multiplicity, is at most $m$, so that $\alpha_I\in\CH^{n_I}(M\times Y^k)$ with $k\leq m$.
Thus $\h(M)$ is a direct summand of $\bigoplus_{k=0}^m\h(Y^k)(\ast)^{\oplus r_k}$, recovering the bound of \cite[Theorem~1.1]{FLV2}.
We now discuss the twists.

Consider a term $(\beta_I\circ\alpha_I)_{s_1,\ldots,s_k}$ such that $(\beta'_I\circ\alpha'_I)_{s_1,\ldots,s_k}\neq0$, and such that $u$ of the~$\h_{s_l}$ are equal to $\h_{\mathrm{tr}}(Y)$ and the remaining $k-u$, which we call the Lefschetz factors, are of the form $\1(-a_l)$.
Then the term factors through $\h_{\mathrm{tr}}(Y)^{\otimes u}(c)$, where
\[
c=n_I-2m-\sum_{\text{Lefschetz factors}}a_l.
\]
For a Lefschetz factor, define $t_l\coloneqq i_l-a_l$.
By Lemma~\ref{lem:vanishing}(i), $t_l\geq0$ and $a_l\geq2d-j_l$, which implies
\[
0\leq t_l=i_l-a_l\leq i_l-(2d-j_l)=i_l+j_l-2d=\mu_l.
\]
For the factors equal to $\h_{\mathrm{tr}}(Y)$, Lemma~\ref{lem:vanishing}(ii) gives $i_l\geq d+1$ and $j_l\geq d+1$, implying $d+1\leq i_l\leq\mu_l+d-1$ and $d+1\leq j_l\leq\mu_l+d-1$.
In this case, define $t'_l\coloneqq i_l-(d+1)$, so that $0\leq t'_l\leq\mu_l-2$.

Substituting $n_I=e+\sum_li_l$, we obtain
\[
c=e+\sum_li_l-2m-\sum_{\text{Lefschetz factors}}a_l=e+\sum_{\text{Lefschetz factors}}t_l+\sum_{\text{factors }\h_{\mathrm{tr}}(Y)}t'_l+(d+1)u-2m.
\]
Using $t_l\geq0$ and $t'_l\geq0$, we obtain $c\geq(d+1)u-2m$.
Furthermore,
\[
\sum_{\text{Lefschetz factors}}t_l+\sum_{\text{factors }\h_{\mathrm{tr}}(Y)}t'_l\leq\sum_{\text{Lefschetz factors}}\mu_l+\sum_{\text{factors }\h_{\mathrm{tr}}(Y)}(\mu_l-2)=\sum_l\mu_l-2u.
\]
Combining this with $e+\sum_l\mu_l\leq e+f+\sum_l\mu_l=2m$, we get
\[
(d+1)u-2m\leq c\leq2m-2u+(d+1)u-2m=(d-1)u.
\]
Since $\h(Y)$ contains $\1(b)$ for every integer $-2d\leq b\leq0$ by the construction in the proof of Lemma~\ref{lem:codim}, the motive $\h_{\mathrm{tr}}(Y)^{\otimes u}\bigl(m(d-1)+\sum_{l=1}^{m-u}b_l\bigr)$ is a direct summand of $\h(Y)^{\otimes m}(m(d-1))=\h(Y^m)(m(d-1))$ for every choice of integers $b_1,\ldots,b_{m-u}\in[-2d,0]$, and $\sum_lb_l$ takes every integer value in $[-2d(m-u),0]$.
Hence every motive $\h_{\mathrm{tr}}(Y)^{\otimes u}(c)$ with $c\in[2du-(d+1)m,(d-1)m]$ is a direct summand of $\h(Y^m)(m(d-1))$.
This applies to all the terms above, since $[(d+1)u-2m,(d-1)u]\subseteq[2du-(d+1)m,(d-1)m]$ as $u\leq m$.
Taking a large enough sum shows that $\h(M)$ is a direct summand of $\h(Y^m)(m(d-1))^{\oplus r}$ for some $r$.
\end{proof}

\subsection{Double EPW cubes as moduli spaces}\label{SS:Ku-GM}

As discussed after Corollary~\ref{cor:sextic-twists}, only the very general double EPW cube $C$ is isomorphic to $M_\sigma(\Ku(Y),w)$, where $w=\Lambda_1-\Lambda_2$.
Consider the open $B^\circ\subseteq B$ over which the associated double EPW cube is smooth and let $\CC_B\to B^\circ$ be the universal family.
Furthermore, let $U\subseteq B^\circ$ be an open on which we have smooth projective relative moduli spaces $M_U(w)\to U$.
Then $\CC_u\cong M_U(w)_u$ for very general $u\in U$ \cite[Corollary~1.3]{KKM}. 
Let 
$$J\coloneqq\mathrm{Isom}_U(\CC_B|_U,M_U(w)) $$
be the relative isomorphism scheme, an open subscheme of the relative Hilbert scheme of $\CC_B\times_U M_U(w)$.
Since its very general fibre is non-empty, an irreducible component $J_0$ of~$J$ dominates $U$.
In particular, $\CC_u\cong M_U(w)_u$ for the general $u\in U$.
The problem is that these isomorphisms a priori do not assemble to an isomorphism $\CC_B\to M_U(w)$ over $U$, so we do not obtain an analogue of Corollary~\ref{cor:sextic-twists}.
We can only prove the following weaker statement, which is still sufficient to prove the Franchetta property for $\CC_B\to B^\circ$.

\medskip

\begin{thm}[Double EPW cubes]
\label{thm:relative-transfer}
Let $B^\circ\subseteq B$ be the dense open subset over which the associated double EPW cube is smooth and let $\CC\to B^\circ$ be the universal family of double EPW cubes.
Let $\iota_C$ be the involution of $\CC$ over $B^\circ$ given by the covering involution of the double EPW cubes, and let
$\h(C_b)^+\coloneqq(C_b,\Delta^+|_{C_b\times C_b})$
be the direct summand of $\h(C_b)$ cut out by the $\iota_C$-invariant projector $\Delta^+\coloneqq\tfrac12(\Delta_{\CC_B/B^\circ}+\Gamma_{\iota_C})$.
Then there is a relative correspondence over $B^\circ$ inducing, for every $b\in B^\circ$, a split injection
\begin{equation}\label{eq:transfer-motives}
\h(C_b)^+\hookrightarrow\bigoplus_\ell\h(Y_b^3)(3)\ \ \hbox{in}\ \Mrat,
\end{equation}
whose left inverse is induced by a relative correspondence over $B^\circ$ as well.
\end{thm}

\begin{proof}
Let $M_U(w)\to U$, $J$, and $J_0$ be as above.
Since $U'\coloneqq J_0\to U$ is generically finite of degree $e\leq 2$, we can assume that  it is \'etale of degree $e\leq 2$ after possibly shrinking $U$.
Over $U'$ there is an isomorphism of $U'$-schemes 
$$\Phi\colon \CC\times_U U'\xrightarrow{\sim}M_U(w)\times_U U'.$$
When $e=1$, $\Phi$ is an isomorphism over $U$, and there is nothing to prove, as the correspondences that appear in the proof of Theorem~\ref{thm:twists} can be defined over $U$. 
When $e=2$ the two isomorphisms $\CC_u\to M_U(w)_u$ given by the preimage of $U'\to U$ differ by the covering involution $\iota_C$.

Now, consider the relative correspondences $\alpha_I$ and $\beta_I$ over $U$ that appear in the proof of Theorem~\ref{thm:twists}, see Section~\ref{SS:twists}. 
Let $\alpha'_I\coloneqq r_\ast\Phi^\ast r^\ast\alpha_I$ and $\beta'_I\coloneqq r_\ast\Phi^\ast r^\ast\beta_I$, where $r\colon U'\to U$ is the \'etale map. 
These are relative correspondences over $U$ between $C$ and $\YY^3_{/U}$.
Writing $\Gamma\coloneqq \sum_{g\in H}\Gamma_g$, where $H\coloneqq \{\id,\iota_C\}$, we obtain
$$\sum_I\beta'_I|_u\circ\alpha'_I|_u=2\Gamma\text{ and }\Delta^+\circ(\sum_I\beta'_I|_u\circ\alpha'_I|_u)\circ\Delta^+=4\Delta^+. $$
Since $\Delta^+$ is a relative cycle over $B^\circ$,
$$\tilde\alpha_I\coloneqq \alpha_I'\circ\Delta^+\text{ and }\tilde\beta_I\coloneqq \frac{1}{4}\Delta^+\circ\beta'_I $$
are relative correspondences over $U$ satisfying $\sum_I\tilde\beta_I\circ\tilde\alpha_I=\Delta^+$ on every fibre over $U$.

The localization sequence yields classes $\overline{\alpha}_I$ and $\overline{\beta}_I$, which restrict to $\tilde\alpha_I$ and $\tilde\beta_I$ over $U$.
The cycle $\sum_I\overline{\beta}_I\circ\overline{\alpha}_I-\Delta^+$ vanishes on $U$, so applying Lemma~\ref{lem:spread}(i) shows that it vanishes on $B^\circ$, completing the proof. 
\end{proof}


\section{Generically defined cycles and the tautological ring}
\label{S:tautological}

In this section, we prove that for every fibre $Y$ of the universal family of GM fourfolds $\YY\to B$ over the parameter space $B$, the generically defined cycles on $Y^m$ are exactly the tautological ones for $m\leq3$, and in codimension at most $8$ for $m=4$.

\subsection{Tautological cycles and stratified projective bundles}

We recall the notion of \textit{tautological rings} (see \cite[Definition~2.3]{FLV2}) and specify it for GM fourfolds:
\begin{defn}[Tautological ring]\label{def:tautological}
Let $P$ be a smooth projective variety with "trivial" Chow groups (i.e.~having Tate Chow motive). Let $Z\subseteq P$ be a smooth subvariety. Define
\begin{align*}
\R^*(Z) &\coloneqq \bigl\langle \ima(\CH^*(P)\to\CH^*(Z)),\, c_i(T_Z) \bigr\rangle \subseteq \CH^*(Z),\\
\R^*(Z^r) &\coloneqq \bigl\langle p_i^*\R^*(Z),\, p_{jk}^*\Delta_Z \bigr\rangle \subseteq \CH^*(Z^r).
\end{align*}
where $p_i$ and $p_{jk}$ are the natural projections of $Z^r$ onto its factors. 
\end{defn}
\begin{ex}
\label{ex:R*forGM4}
    For an ordinary GM fourfold $Y$, we always take $P=Y'$, the del Pezzo manifold containing $Y$, which has trivial Chow groups by Lemma~\ref{lem:trivial-chow}, so that  
    $$\R^*(Y)=\bigl\langle h, c_2\bigr\rangle \subseteq\CH^*(Y)$$ 
    thanks to Lemma~\ref{lem:trivial-chow} (the Chern classes of $T_Y$ come from $Y'$), and
\[
\R^*(Y^m)=\bigl\langle h_i,\,c_{2,i},\,\Delta_{ij}\bigr\rangle\subseteq\CH^*(Y^m),
\]
where $h_i\coloneqq p_i^*h, c_{2,i}\coloneqq p_i^*c_2,\Delta_{ij}\coloneqq p_{ij}^*\Delta_Y$.
For an arbitrary GM fourfold, we take the right-hand side as the definition of $\R^*(Y^m)$.
\end{ex}

Given a vector bundle $E$ on a variety $P$ and a subspace $V\subseteq\HH^0(P,E)$, we say that $(P,V)$ satisfies \emph{condition $(\star_r)$} if for all \textit{pairwise distinct} points $x_1,\ldots,x_r\in P$ the evaluation map $$V\to\bigoplus_{i=1}^rE(x_i)$$ is surjective, where $E(x)$ denotes the fibre of $E$ at $x$; see \cite[Definition~5.6]{FLV1} and \cite[Definition~2.5]{FLV2}. Clearly $(\star_r)$ implies $(\star_k)$ for all $k<r$, and $(\star_1)$ amounts to saying that $E$ is generated by sections in $V$.

The following result from \cite[Proposition~2.6]{FLV2} (see also \cite[Proposition 5.7]{FLV1}) was established by the so-called stratified projective bundle techniques:

\begin{prop}[Generic versus tautological]\label{prop:strat}
Let $P$ be a smooth projective variety with motive of Tate-type, let $E$ be a vector bundle on $P$ and let $V\subseteq\HH^0(P,E)$ be a subspace of sections generating $E$. Let $B_V\subseteq\PP(V)$ be the open subset parameterizing smooth zero loci of dimension $\dim P-\rk E$, let $\Zz\to B_V$ be the universal family and let $Z$ be the fibre over a closed point. Let $r\geq1$ be an integer, assume that $(P,V)$ satisfies $(\star_{r})$. Then
$$\GDCH^*_{B_V}(Z^r)=\R^*(Z^r).$$
\end{prop}

In this paper, we will need to apply Proposition \ref{prop:strat} in a more general situation where condition $(\star_r)$ is not necessarily satisfied. To this end, when condition $(\star_r)$ fails, we introduce the following \textit{degeneracy locus} of condition $(\star_r)$:
\[
T_r^\circ\coloneqq \left\{(x_1, \dots, x_r)\in P^r\setminus \bigcup_{j<k} \Delta_{j,k} ~\mid~ \text{ the evaluation map } V\to\bigoplus_{i=1}^rE(x_i) \text{ is not surjective}\right\}.
\]
Let $T_r\subset P^r$ be the Zariski closure of $T_r^{\circ}$. The following result used inductively and combined with Proposition~\ref{prop:strat} allows us to control the difference between generically defined cycles and tautological cycles in terms of the contributions from the degeneracy loci.

\begin{prop}
\label{prop:GenericVSTautWithFailingLoci}
Let $P$ be a smooth projective variety with motive of Tate-type, let $E$ be a vector bundle on $P$ and let $V\subseteq\HH^0(P,E)$ be a subspace of sections generating $E$. Let $B_V\subseteq\PP(V)$ be the open subset parameterizing smooth zero loci of dimension $\dim P-\rk E$, let $\Zz\to B_V$ be the universal family and let $Z$ be the fibre over a closed point. Let $r\geq1$ be an integer. Then we have
\[\GDCH^*_{B_V}(Z^r)\subseteq \R^*(Z^r)+\sum_{j<k}\delta_{j,k}\GDCH^*_{B_V}(Z^{r-1})+\ima \Bigl(\CH^*(T_r\cap Z^r)\to \CH^*(Z^r)\Bigr),\]
where $\delta_{j,k}\colon Z^{r-1}\hookrightarrow Z^r$ is the closed immersion of the big diagonal in $Z^r$ defined by $x_j=x_k$.
\end{prop}
\begin{proof}
Consider the $r$-th relative power of the universal family:
    \[
\Zz^r\coloneqq\Zz^{\times r}_{/B_V}\coloneqq\bigl\{(x_1,\ldots,x_r,[s])\in P^r\times\PP(V)\ \bigm|\ s(x_1)=\cdots=s(x_r)=0\bigr\}.
\]
Denote $(P^r)^\circ\coloneqq P^r\setminus  \bigcup_{j<k} \Delta_{j,k}$, the locus of pairwise distinct $r$-tuples in $P$. We have the following commutative diagram where all the squares are cartesian:
\[
\begin{tikzcd}
  T_r^\circ \cap Z^r \arrow[d, hook] \arrow[r, hook]& \mathcal{W}_r^\circ \arrow[r] \arrow[d, hook] & T_r^\circ \arrow[d, hook] \\
  U_r \arrow[r, hook] \arrow[d, hook] & \mathcal{U}_r \arrow[r, "q^\circ"] \arrow[d, hook] & (P^r)^\circ \arrow[d, hook] \\
  Z^r \arrow[r, hook] \arrow[d] & \Zz_r \arrow[r, "q"] \arrow[d, "\pi"] & P^r \\
  \{b\} \arrow[r, hook] & \PP(V) &
\end{tikzcd}
\]
where $Z$ is the fiber over $b\in B_V\subset \PP(V)$, $U_r=Z^r\setminus  \bigcup_{j<k} \Delta_{j,k}$, the maps from the second row to the third row are open immersions, $q$ is the evaluation map and $\pi$ is the natural projection. 

Let $\mathcal{W}_r$ be the Zariski closure of $\mathcal{W}_r^\circ$ in $\Zz_r$. We have the following commutative diagram:
\[
\begin{tikzcd}
  & \CH(\mathcal{W}_r)\oplus \bigoplus_{i=0}^{N}\CH(P^r) \arrow[r, two heads] \arrow[d, "f"] & \CH(\mathcal{W}_r^\circ)\oplus \bigoplus_{i=0}^{N}\CH((P^r)^\circ) \arrow[d, two heads, "f^\circ"]  &\\
  \bigoplus_{j<k}\CH(\Zz_{r-1}) \arrow[r, "\sum_{j<k}\delta_{j,k*}"] \arrow[d] & \CH(\Zz_r) \arrow[r, two heads] \arrow[d] & \CH(\mathcal{U}_r) \arrow[d] \arrow[r]& 0 \\
  \bigoplus_{j<k}\CH(Z^{r-1}) \arrow[r, "\sum_{j<k}\delta_{j,k*}"] & \CH(Z^r) \arrow[r, two heads] & \CH(U_r) \arrow[r] & 0
\end{tikzcd}
\]
where the second and the third rows are exact by the localization sequence for Chow groups, the map $f$ is defined as follows: $\CH(\mathcal{W}_r)\to \CH(\Zz_r)$ is the push-forward by inclusion, the morphism $\bigoplus_{i=0}^{N}\CH(P^r) \to \CH(\Zz_r)$ is given by $\sum_{i=0}^{N} q^*(-)\cdot \xi^i$, where $\xi=\pi^*c_1(\mathcal{O}_{\PP(V)}(1))$, the map $f^\circ$ is similarly defined. 

Now by the definition of the degeneracy locus $T_r^\circ$, the morphism $q^\circ\colon \mathcal{U}_r\to (P^r)^\circ$ is a $\PP^{N}$-bundle outside of the locus $T_r^\circ$, where $N=\dim(\PP(V))-r$. Consequently, the arrow~$f^\circ $  in the above diagram is surjective. A diagram chasing shows that
\[\CH^*(\Zz_r)=\ima (f) + \sum_{j<k}\delta_{j,k*}\CH^*(\Zz_{r-1}).\]
Composing with the Gysin restriction map to $\CH^*(Z^r)$ and using the commutativity of the above diagram, we obtain
\begin{align}
    \begin{split}
    &\GDCH_{B_V}(Z^r)=\\
    &\  \sum_{j<k}\delta_{j,k\ast}\GDCH^*_{B_V}(Z^{r-1})
     + \ima\left(\CH^*(\mathcal{W}_r)\oplus \bigoplus_{i=0}^{N}\CH^*(P^r) \xrightarrow{f} \CH^*(\Zz_r) \to \CH^*(Z^r)\right)\\
    =&\sum_{j<k}\delta_{j,k\ast}\GDCH^*_{B_V}(Z^{r-1})+\ima \Bigl(\CH^*(\mathcal{W}_r)\to \CH^*(Z^r)\Bigr)+\ima\Bigl(\CH^*(P^r)\to \CH^*(Z^r)\Bigr)\\
    =&\sum_{j<k}\delta_{j,k\ast}\GDCH^*_{B_V}(Z^{r-1})+\ima \Bigl(\CH^*(\mathcal{W}_r)\to \CH^*(Z^r)\Bigr)+\R^*(Z^r),
\end{split}
\end{align}
where in the last equality, we used the fact that $\xi$ restricts to 0 in $\CH^*(Z^r)$. Finally, since $\ima (\CH^*(\mathcal{W}_r)\to \CH^*(Z^r))\subseteq \ima (\CH^*(T_r\cap Z^r)\to \CH^*(Z^r))$, we can conclude.
\end{proof}

\subsection{Generically defined cycles on powers of a Gushel--Mukai fourfold}
Now we apply Proposition~\ref{prop:strat} and Proposition~\ref{prop:GenericVSTautWithFailingLoci} to GM fourfolds. Recall the definition of the tautological rings for GM fourfolds in Example~\ref{ex:R*forGM4}.

\begin{lem}\label{lem:GDCH-Y}
For every GM fourfold $Y=Y_b$, $b\in B$, one has $\GDCH^*_B(Y)=\R^*(Y)=\langle h, c_2\rangle$.
\end{lem}
\begin{proof}
The condition $(\star_1)$ is satisfied and the lemma follows from Proposition~\ref{prop:strat}.
\end{proof}

\begin{prop}\label{prop:generation}
Let $Y=Y_b$, $b\in B$, be any GM fourfold. Then
\begin{enumerate}[(i)]
    \item $\GDCH^*_B(Y^m)=\R^*(Y^m)\quad\text{for }m\leq3$;
    \item $\GDCH^i_B(Y^m)=\R^i(Y^m)\quad\text{for }i\leq8$ and $m\leq 6$.
\end{enumerate}
\end{prop}

\begin{proof}
The inclusion $\R^*(Y^m)\subseteq\GDCH^*_B(Y^m)$ is clear.
We prove the reverse inclusions.  
By Lemma~\ref{lem:spread}(ii), it suffices to prove the inclusion $\GDCH^i_B(Y_b^m)\subseteq\R^i(Y_b^m)$ for all $b$ in the dense open subset $B^{\mathrm{ord}}$ of $B$, where $B^{\mathrm{ord}}$ is the locus of ordinary GM fourfolds. Let $b\in B^{\mathrm{ord}}$ such that $Y=Y_b=\CGr\cap H_1\cap H_2\cap Q$.

\smallskip
\noindent\emph{Reduction to quadratic sections of the fivefold $Y'$.}
The del Pezzo fivefold $Y'\coloneqq\CGr\cap H_1\cap H_2$ is isomorphic to $\Gr\cap H$ for a hyperplane $H\subseteq\PP(\bigwedge^2V_5)$; we denote also $\Gr\cap H$ by $Y'$. By Lemma~\ref{lem:trivial-chow}, $Y'$ has trivial Chow groups and $\CH^*(Y')=\langle h,c_2\rangle$. 
Defining
\[
B'\coloneqq\bigl\{Q\in\PP\HH^0(\CGr,\OO(2))\ \bigm|\ Y'\cap Q \text{ is smooth of dimension }4\bigr\},
\]
an open subset of $|\OO(2)|$ containing $Q$, we have a natural map
\[
 B'\longrightarrow B
\]
sending $Q$ to $b=(Q, H_1, H_2)$.
By Lemma~\ref{lem:base-change}(ii),
\begin{equation}\label{eq:red-to-Yprime}
\GDCH^*_B(Y^m)\subseteq\GDCH^*_{B'}(Y^m).
\end{equation}

Let $V_0\coloneqq\ima\bigl(\HH^0(\CGr,\OO(2))\to\HH^0(Y',\OO_{Y'}(2))\bigr)$ and let $B''\subseteq\PP(V_0)$ be the open subset parameterizing smooth zero loci of dimension $4$; the subspace $V_0$ generates $\OO_{Y'}(2)$, since it contains the restrictions to $Y'$ of all quadrics of $\PP^{10}$.  
Choose a linear splitting $\lambda\colon V_0\to\HH^0(\CGr,\OO(2))$ of the restriction map, 
we have that by Lemma~\ref{lem:base-change} and \eqref{eq:red-to-Yprime},
\begin{equation}\label{eq:red-to-Bpp}
\GDCH^*_B(Y^m)\subseteq\GDCH^*_{B'}(Y^m)\subseteq\GDCH^*_{B''}(Y^m).
\end{equation}
We are thus reduced to the family of quadratic sections of the fixed del Pezzo fivefold $Y'$, to which Proposition~\ref{prop:strat} applies with $P=Y'$, $E=\OO_{Y'}(2)$ and $V=V_0$. Note that the tautological ring of $Y^m$ relative to the ambient variety $Y'$ is $\R^*(Y^m)$.

\smallskip
\noindent\emph{Independent conditions on lines and quadrics.}
Let $x_1,\ldots,x_m\in\PP^{10}$ be pairwise distinct points. Consider the evaluation map
\begin{equation}\label{eq:quadric-evaluation}
\operatorname{ev}_m\colon \HH^0\bigl(\PP^{10},\OO(2)\bigr)\longrightarrow\bigoplus_{i=1}^m\C_{x_i}.
\end{equation}
By \cite[Proposition~1]{EGH}, we have the following independence properties:
\begin{itemize}
    \item ($m=3$) it is clear that $\operatorname{ev}_3$ is surjective;
    \item ($m=4, 5$) the map $\operatorname{ev}_m$ is surjective if and only if no four points among $\{x_1, \dots, x_m\}$ are collinear;
    \item ($m=6$) the map $\operatorname{ev}_6$ is surjective if and only if $x_1, \dots, x_6$ do not lie on a conic and have no four collinear points.
\end{itemize}
We will need the fact that any conic $c$ containing $6$ points of $Y'$, of which no four are collinear, is contained in $Y'$. Indeed, assume first that $c$ is irreducible. 
    Any quadric containing $Y'$ either contains the entire plane spanned by $c$ or meets it in a conic meeting~$c$ in $6$ points, hence containing $c$.
    A hyperplane containing $Y'$ either contains the entire plane spanned by $c$ or meets it in a line, meeting $c$ in at most $2$ points. Hence $c\subset Y'$.
    Now assume $c=\ell_1\cup \ell_2$ is reducible. 
    Since no four of the points are collinear, each $\ell_i$ contains three of them and the same argument shows that each quadric and hyperplane containing~$Y'$ contains each $\ell_i$, hence $c$.
    Double lines cannot occur, since then all $6$ points would be collinear.

\smallskip
\noindent\emph{The case $m\leq3$.}
Since $V_0$ contains the restrictions to $Y'$ of all quadrics of $\PP^{10}$, the pair $(Y',V_0)$ satisfies  the condition $(\star_3)$.
By Proposition~\ref{prop:strat} and \eqref{eq:red-to-Bpp} we get \[ \GDCH^*_B(Y^m)\subseteq\GDCH^*_{B''}(Y^m)=\R^*(Y^m) \ \ \hbox{for}\ m\leq 3,\] 
and Lemma~\ref{lem:spread}(ii) extends the equality $\GDCH^*_B(Y_b^m)=\R^*(Y_b^m)$ from $b\in B^{\mathrm{ord}}$ to every $b\in B$. This proves (i).

\smallskip
\noindent\emph{The case $m>3$ (reduce to dimension estimate).}
The condition $(\star_3)$ is satisfied, but the condition $(\star_4)$ does not necessarily hold. Applying Proposition~\ref{prop:strat} and Proposition~\ref{prop:GenericVSTautWithFailingLoci} inductively, we obtain that
\begin{equation}
\label{eqn:GDCHYmWithCorrections}
    \GDCH_B^\ast(Y^m)=\R^\ast(Y^m)+\sum_{\substack{\delta\colon Y^r\hookrightarrow Y^m\\4\leq r\leq m}} \delta_\ast\Bigl(\ima\bigl(\CH_*(T_r\cap Y^r)\to\CH_\ast(Y^r)\bigr)\Bigr),
\end{equation}
where the sum runs through all partial diagonal embeddings $\delta\colon Y^r\hookrightarrow Y^m$.
Therefore, we need to study the degeneracy loci $T_4\cap Y^4$, $T_5\cap Y^5$, $T_6\cap Y^6$. To conclude part (ii) of the proposition, it suffices to show the following bounds for their dimensions:
\[\dim(T_4\cap Y^4)\leq 7, \qquad\dim(T_5\cap Y^5)\leq 11, \qquad\dim(T_6\cap Y^6)\leq 15. \]

\smallskip
\noindent\emph{The case $m=4$.}
We identify $Y'$ with $\Gr\cap H\subseteq\PP(\bigwedge^2V_5)$. It is clear that 
\[
T_4=\bigcup_{\ell\in F_1(Y')}\ell\times\ell\times\ell\times\ell,
\]
where $F_1(Y')$ denotes the Hilbert scheme of lines contained in $Y'$, and
\[
T_4\cap Y^4=\bigcup_{\ell\in F_1(Y')}(\ell\cap Y)^4 .
\]
If $\ell\subseteq Y'$ is a line with $\ell\subseteq Y$, its contribution to this union is $\ell^4$, of dimension $4$. If $\ell\not\subseteq Y$, then $\ell\cap Y=\ell\cap\{q=0\}$ is a subscheme of length $2$ of $\ell$, so that $(\ell\cap Y)^4$ is finite. 
Therefore,
\[
\dim\bigl(T_4\cap Y^4\bigr)\leq\max\{\dim(F_1(Y))+\dim(\ell^4),\,\dim(F_1(Y'))\}.
\]
Now $\dim F_1(Y)=3$ by \cite[Lemma~3.8]{DK5} (see also \cite[Theorem~4.7]{DK2}), while $\dim F_1(Y')=6$, hence we obtain that  $\dim\bigl(T_4\cap Y^4\bigr)\leq 7$.

\smallskip
\noindent\emph{The case $m=5$.} It is similar to the $m=4$ case. 
 We need to bound the dimension of $T_5\cap Y^5$, which is the locus of five points with at least four of them being collinear. If the line $\ell$ is contained in $Y$, then four points are constrained in $\ell$ while the remaining point can move freely on $Y$.
    When $\ell$ is not contained in $Y$, four points are constrained in $\ell\cap Y$ which consists of at most two points, and the remaining point can move freely on $Y$. Therefore,
    $$\dim(T_5\cap Y^5)\leq\max\{\dim F_1(Y)+\dim(\ell^4)+\dim(Y),\dim F_1(Y')+\dim(Y)\}.$$
    Using $\dim F_1(Y)=3$ and $\dim F_1(Y')=6$, we obtain $\dim(T_5\cap Y^5)\leq 11$.

\smallskip
\noindent\emph{The case $m=6$.} 
    The locus $ T_6\cap Y^6$ is the union of the locus $L\subset Y^6$ where at least 4 points are collinear and the locus $M\subset Y^6$ where the 6 points lie on a conic. 
    For the locus $L$, we reason as in the case $m=5$ (but now with two points moving freely on $Y$):
    \[\dim(L)\leq \max\{\dim F_1(Y)+\dim(\ell^4)+2\dim(Y),\dim F_1(Y')+2\dim(Y)\}=15.\]
    As for the locus $M$,
    let $(x_1, \dots, x_6)\in M$ with conic $c$.
    Denote the union of the irreducible components of $c$ contained in $Y$ by $c_Y$.
    Similarly to the fact that $\ell\cap Y$ consists of at most two points, one can check that $(c\cap Y)\setminus c_Y$ consists of at most $4$ points. 
    Let us denote the Hilbert schemes of conics on $Y,Y'$, $\Gr$ by $F(Y),F(Y'), F(\Gr)$ respectively. We have $\dim F(\Gr)=13$ and $\dim F(Y)=5$ for $Y$ general, cf. \cite{IM} and Section~\ref{S:MCK}.
   
    We now distinguish between $4$ cases:

    \noindent (1) If $c_Y=\varnothing$, $x$ lies in a finite set, and we use $\dim F(Y') \leq\dim F(\Gr)=13$.

    \noindent (2) If $c\subset Y$, using $\dim F(Y)=5$, these form a set of dimension at most $11$.

    \noindent (3) If $c$ is a double line, we have $x\in\ell^6\subset L$.

    \noindent (4) If $c=\ell_1\cup\ell_2$ with $\ell_1\subset Y$ and $\ell_2\nsubset Y$, we distinguish between four of the points lying on $\ell_1$, and at least three lying on $\ell_2\cap Y$, which is finite. 
    The former case we already considered, and the latter results in a set of dimension at most $\dim F_1(Y)+\dim F_1(Y')+3=12$.
    Hence $\dim(M)\leq 13$.
    
    We now conclude that
    \[\dim(T_6\cap Y^6)\leq\max\{15,13\}=15,\] 
   which completes the proof.
\end{proof}

\begin{rmk}\label{rmk:Y4}
For $b$ in the dense open subset $B^{\mathrm{ord}}$ of ordinary GM fourfolds, Proposition~\ref{prop:GenericVSTautWithFailingLoci} and Proposition \ref{prop:generation} give the more precise inclusion
\begin{equation}
\label{eqn:GDCHforY4}
    \GDCH^*_B(Y^4)\subseteq\R^*(Y^4)+\ima\Bigl(\CH_*(T_4\cap Y^4)\to\CH_*(Y^4)\Bigr),
\end{equation}
whose error term is concentrated in codimension at least $9$, where the component with largest dimension
is the image of $L\times_{F_1(Y)}L\times_{F_1(Y)}L\times_{F_1(Y)}L$ in $Y^4$, where $L\to F_1(Y)$ is
the universal line, so that the contribution of the second term of \eqref{eqn:GDCHforY4} is governed by the Chow groups of $F_1(Y)$.  The
latter is birational to a double cover of a hyperplane section of an EPW sextic
\cite[Theorem~4.7(c)]{DK2}, a threefold with ample canonical class.
\end{rmk}


\section{Homological relations among tautological cycles}\label{S:relations}

In this section, we consider a smooth projective variety $Y$ equipped with a subalgebra of tautological cycles satisfying the hypotheses (H1)--(H3) below; this includes GM fourfolds (Section~\ref{SS:GM-u}) and GM sixfolds (Section~\ref{S:sixfolds}).
We define the canonical Chow--K\"unneth decomposition of $Y$ and, following Tavakol \cite{Tavakol} and Yin \cite{Yin}, we give a presentation of the image in cohomology of the tautological ring of $Y^m$. This isolates classes on $Y$, $Y^2$ and $Y^3$ whose vanishing in $\CH^*$ (resp.\ in $\B^*$) implies that the cycle class map is injective on the tautological ring of $Y^m$ (resp.\ on its image in $\B^*(Y^m)$), for all $m$ up to a bound depending only on the cohomology of $Y$. For future reference, we work in greater generality than needed in this paper; in particular, we allow $Y$ to be of odd dimension.

\subsection{Setting}\label{SS:setting-u}

Let $Y$ be a smooth projective variety of dimension $n$, and let $\R^*(Y)\subseteq\CH^*(Y)$ be a graded $\QQ$-subalgebra, whose elements we call \emph{tautological}. 
Consider the restriction to $\R^*(Y)$ of the cycle class map $\cl\colon \R^*(Y) \to \HH^{*}(Y,\QQ)$.
We write
$$\HH^{2*}_{\taut}(Y)\coloneqq \operatorname{im} \bigl(\cl\colon \R^*(Y) \to \HH^{2*}(Y,\QQ)\bigr) \quad \text{and} \quad K\coloneqq\ker\bigl(\cl\colon\R^*(Y)\to\HH^{2*}(Y,\QQ)\bigr),$$ 
and we fix a graded $\QQ$-linear section $s\colon\HH^{2*}_{\taut}(Y)\to\R^*(Y)$ of the cycle class map, so that $\R^*(Y)=s(\HH^{2*}_{\taut}(Y))\oplus K$. 
We assume the following:

\begin{enumerate}[(H1)]
	\item $\HH^{2i+1}(Y,\QQ)=0$ for $2i+1\neq n$;
	\item $\HH^{2i}_{\taut}(Y)=\HH^{2i}(Y,\QQ)$ for $2i\neq n$;
	\item if $n$ is even, the intersection form is non-degenerate on $\HH^n_{\taut}(Y)$.
\end{enumerate}
Let $\HH^n_{\van}(Y)$ be the orthogonal complement of $\HH^n_{\taut}(Y)$ in $\HH^n(Y,\QQ)$ and denote by $d$ its dimension.
By (H3), $\HH^n(Y,\QQ)=\HH^n_{\taut}(Y)\oplus\HH^n_{\van}(Y)$
and we have
\begin{equation}\label{eq:van-kills-u}
\HH^n_{\van}(Y)\cdot \HH^{>0}_{\taut}(Y) = 0.
\end{equation}
By (H2), $\HH^{2n}(Y,\QQ)=\HH^{2n}_{\taut}(Y)$ is spanned by the class $[\mathrm{pt}]$ of a point, and we put $o\coloneqq s([\mathrm{pt}])\in\R^n(Y)$, a tautological zero-cycle of degree one.

For $m\geq1$, let
\[
\R^*(Y^m)\coloneqq\bigl\langle p_i^*\R^*(Y),\,p_{jk}^*\Delta_Y\bigr\rangle\subseteq\CH^*(Y^m)
\quad\text{and}\quad
\mathrm{D}^*(Y^m)\coloneqq\bigl\langle p_i^*\R^*(Y)\bigr\rangle\subseteq\R^*(Y^m)
\]
be the algebras of tautological and of \emph{decomposable} tautological classes.
Note that when $\R^*(Y)$ is the tautological ring of Section~\ref{S:tautological}, $\R^*(Y^m)$ coincides with the tautological ring of~$Y^m$ of Definition~\ref{def:tautological}.

\begin{lem}\label{lem:decomposable-u}
	$\ker \bigl(
	\mathrm{D}^*(Y^m) \to \HH^*(Y^m,\Q)\bigr)$ is  generated by $p_1^*K,\ldots,p_m^*K$.
\end{lem}
\begin{proof}
	This follows at once from the fact that the exterior product $\HH^*(Y,\QQ)^{\otimes m}\to\HH^*(Y^m,\QQ)$ is injective by the K\"unneth formula.
\end{proof}

\subsection{The canonical Chow--K\"unneth decomposition}\label{SS:CK-u}

\begin{defn}\label{defn:CK-u}
	For $2i\neq n$, let $(a_\alpha)$ be a basis of $\HH^{2i}(Y,\QQ)$ and $(a^\alpha)$ the basis of $\HH^{2n-2i}(Y,\QQ)$ such that $\deg(a^\alpha a_\beta)=\delta_{\alpha\beta}$; if $n$ is even, let $(u_\mu)$ be a basis of $\HH^n_{\taut}(Y)$ and $(u^\mu)$ the basis of $\HH^n_{\taut}(Y)$ such that $\deg(u^\mu u_\nu)=\delta_{\mu\nu}$, which exists by (H3). Set
	\[
	\pi^{2i}_Y\coloneqq\sum_\alpha s(a^\alpha)\times s(a_\alpha)\ \text{for}\ 2i\neq n,\quad
	\pi^n_{\taut}\coloneqq\sum_\mu s(u^\mu)\times s(u_\mu),\quad
	\tau\coloneqq\Delta_Y-\sum_{2i\neq n}\pi^{2i}_Y-\pi^n_{\taut},
	\]
	where $\pi^n_{\taut}$ is set to 0 if $n$ is odd, and put $\pi^n_Y\coloneqq\pi^n_{\taut}+\tau$ and $\pi^k_Y\coloneqq0$ for odd $k\neq n$. We call $\{\pi^k_Y\}_{0\leq k\leq2n}$ the \emph{canonical Chow--K\"unneth decomposition} of $(Y,\R^*(Y),s)$, and we write
	\[
	\h^k(Y)\coloneqq(Y,\pi^k_Y),\qquad \h^n_{\taut}(Y)\coloneqq(Y,\pi^n_{\taut}),\qquad \h^n_{\van}(Y)\coloneqq(Y,\tau).
	\]
	We further consider the classes
	\[
	\kappa\in K,\qquad
	\Upsilon_z\coloneqq\tau\cdot p_1^*z\in\R^*(Y\times Y)\ \text{for}\ z\in\R^{>0}(Y),\qquad
	\Theta_Y\coloneqq\tau_{12}\,\tau_{13}-\tau_{23}\,o_1\in\R^{2n}(Y^3),
	\]
	where $\tau_{jk}\coloneqq p_{jk}^*\tau$ and $o_i\coloneqq p_i^*o$.
\end{defn}

Let us justify the terminology \textit{canonical Chow--K\"unneth decomposition}. The sums defining $\pi^{2i}_Y$ and $\pi^n_{\taut}$ depend only on $s$, not on the chosen bases. The identity $(x\times x')\circ(y\times y')=\deg(x\cdot y')\,y\times x'$ with the convention that $\deg$ vanishes on cycles of non-zero dimension, together with $\deg(s(a)s(a'))=\deg(aa')$ for $a,a'\in\HH^*_{\taut}(Y)$, shows that the $\pi^k_Y$ are mutually orthogonal idempotents summing to $\Delta_Y$, with ${}^t\pi^k_Y=\pi^{2n-k}_Y$, and that $\pi^n_{\taut}$ and $\tau$ are symmetric orthogonal idempotents. In cohomology, $[\pi^{2i}_Y]$ is the K\"unneth projector onto $\HH^{2i}(Y,\QQ)$ for $2i\neq n$, $[\pi^n_{\taut}]$ is the orthogonal projector onto $\HH^n_{\taut}(Y)$ and $[\tau]$ is the orthogonal projector onto $\HH^n_{\van}(Y)$. 
By (H1) and (H2), $\{\pi^k_Y\}$ is thus a self-dual Chow--K\"unneth decomposition of $Y$.

\begin{prop}\label{prop:CK-basic-u} \hfill
	\begin{enumerate}[(i)]
		\item The action of the idempotent $\tau$ on $\R^*(Y)$ has image $K$ and kernel $s(\HH^{2*}_{\taut}(Y))$; in particular, $p_{1*}\tau=p_{2*}\tau=0$ in $\CH^0(Y)$.
		\item The classes $\tau\cdot o_1$, $\tau\cdot o_2$ and $\tau^2-(-1)^nd\,o_1o_2$ in $\CH^{2n}(Y\times Y)$ are homologically trivial; in particular, they vanish if $\CH_0(Y)=\QQ$.
		\item The classes $\kappa\in K$, $\Upsilon_z$ for $z\in\R^{>0}(Y)$, and $\Theta_Y$ are homologically trivial.
	\end{enumerate}
\end{prop}

\begin{proof} Item (i) is straightforward from the definition of $\tau$.
	Writing
	\begin{equation}\label{eq:tau-tensor-u}
	[\tau]=\sum_{\mu=1}^{d}e_\mu\otimes e^\mu\quad\in\HH^n_{\van}(Y)\otimes\HH^n_{\van}(Y)\subseteq\HH^{2n}(Y\times Y,\QQ),
	\end{equation}
	with $(e_\mu)$ a basis of $\HH^n_{\van}(Y)$ and $(e^\mu)$ the basis such that $\deg(e^\mu e_\nu)=\delta_{\mu\nu}$, and taking into account \eqref{eq:van-kills-u}, items (ii) and (iii) follow from straightforward cohomological computations.
\end{proof}

\subsection{Multiplicativity}
If $Y$ is smooth projective of dimension $n$, the small diagonal $\delta_Y\in\CH^{2n}(Y^3)$, viewed as a correspondence from $Y\times Y$ to $Y$, induces the intersection product. Following Shen--Vial~\cite{SV}, a Chow--K\"unneth decomposition $\{\pi^k_Y\}$ of $Y$ is \emph{multiplicative} (MCK decomposition) if
\begin{equation*}
\pi^k_Y\circ\delta_Y\circ(\pi^i_Y\otimes\pi^j_Y)=0\quad\text{in }\CH^{2n}(Y^3)\qquad\text{for all }k\neq i+j;
\end{equation*}
see \cite[Lemma~2.6]{FLV-MCK} for equivalent formulations.

For $\Gamma_1,\Gamma_2,\Gamma_3\in\CH^n(Y\times Y)$, we put
\[
M(\Gamma_1,\Gamma_2;\Gamma_3)\coloneqq\Gamma_3\circ\delta_Y\circ(\Gamma_1\otimes\Gamma_2)\in\CH^{2n}(Y^3),
\]
so that $\{\pi^k_Y\}$ is multiplicative if and only if $M(\pi^i_Y,\pi^j_Y;\pi^k_Y)=0$ for all $k\neq i+j$.

\begin{prop}[Multiplicativity criterion]\label{prop:MCK-criterion-u}
	Let $G\subseteq\R^{>0}(Y)$ be a set of algebra generators of $\R^*(Y)$.
	\begin{enumerate}[(i)]
		\item If $K=0$, $\Upsilon_g=0$ for all $g\in G$, and $\Theta_Y=0$, then the canonical Chow--K\"unneth decomposition of $Y$ is multiplicative. The same holds modulo algebraic equivalence.
		\item Conversely, assume that $K=0$. If the canonical Chow--K\"unneth decomposition of $Y$ is multiplicative, then $\Upsilon_z=0$ for all $z\in\R^{>0}(Y)$, and $\Theta_Y=0$.
	\end{enumerate}
\end{prop}

\begin{proof}
	(i) Since $K=0$, $\tau_*$ vanishes on $\R^*(Y)$ by Proposition~\ref{prop:CK-basic-u}(i). 
	On the other hand,
	the vanishing of the $\Upsilon_g$ for $g\in G$ gives
	\begin{equation}\label{eq:tau-contract-u}
	\tau\cdot p_1^*z=\tau\cdot p_2^*z=0\qquad\text{for all }z\in\R^{>0}(Y).
	\end{equation}
	The projectors $\pi^k_Y$, $k\neq n$, and $\pi^n_{\taut}$ are decomposable and it follows from the above that $M(\pi^i_Y,\pi^j_Y;\pi^k_Y)=0$ for $k\neq i+j$, except possibly for $(i,j,k)=(n,n,n)$, in which case we have $M(\pi^n_Y,\pi^n_Y;\pi^n_Y) = M(\tau,\tau;\tau)$.  Now
	\begin{equation}\label{eq:tau-small-diag-u}
	\tau_{13}\cdot\tau_{23}=\delta_Y\circ(\tau\otimes\tau)=\sum_kM(\tau,\tau;\pi^k_Y)=\tau_{12}\,o_3+M(\tau,\tau;\tau) \quad \text{in}\ \CH^{2n}(Y^3),
	\end{equation}
	where in the last equality, we have simply used that $M(\tau,\tau;\pi^{2n}_Y)=M(\tau,\tau;Y\times o)=\tau_{12}\,o_3$.
	Since $\Theta_Y=0$ is equivalent, by permuting the factors of $Y^3$, to $\tau_{13}\tau_{23}=\tau_{12}o_3$, we get $M(\tau,\tau;\tau)=0$, and the decomposition is multiplicative. The argument only uses products, pushforwards and compositions of correspondences, so it holds verbatim in $\B^*$.

	(ii) Assume that the decomposition is multiplicative, and let $\CH^i_{(t)}(-)\coloneqq(\pi^{2i-t})_*\CH^i(-)$ be the induced bigrading, for which $\CH^i_{(t)}\cdot\CH^{i'}_{(t')}\subseteq\CH^{i+i'}_{(t+t')}$ \cite[Definition~8.1]{SV}.
	As $K=0$, every $z\in\R^i(Y)$ satisfies $(\pi^{2i}_Y)_*z=z$, i.e., 
	$\R^*(Y)\subseteq\CH^*_{(0)}(Y)$.

	By \cite[Theorem~8.6]{SV} (see also \cite[Proposition~2.9]{FLV-MCK}), $Y\times Y$ carries the product multiplicative decomposition $\pi^k_{Y\times Y}=\sum_{l+l'=k}\pi^l_Y\otimes\pi^{l'}_Y$, with $(\pi^l_Y\otimes\pi^{l'}_Y)_*\Gamma=\pi^{l'}_Y\circ\Gamma\circ{}^t\pi^l_Y$ for $\Gamma\in\CH^*(Y\times Y)$. By orthogonality, $(\pi^l_Y\otimes\pi^{l'}_Y)_*\tau$ is $\tau$ for $l=l'=n$ and $0$ otherwise, so $\tau\in\CH^n_{(0)}(Y\times Y)$. 
	Fix now $z\in \R^i(Y)$ for some $i>0$.
	Since $p_1^*z=z\times Y$ lies in $\CH^i_{(0)}(Y\times Y)$, we have $\Upsilon_z\in\CH^{n+i}_{(0)}(Y\times Y)$, that is,
	\[
	\Upsilon_z=\sum_{l+l'=2n+2i}\pi^{l'}_Y\circ\Upsilon_z\circ\pi^{2n-l}_Y.
	\]
	Now $\Upsilon_z=\tau\circ\Gamma_z$, where $\Gamma_z\coloneqq\Delta_Y\cdot p_1^*z$ is the correspondence whose action induces multiplication by $z$. The terms with $l'\neq n$ vanish because $\pi^{l'}_Y\circ\tau=0$ and so $\Upsilon_z=\pi^{n}_Y\circ\tau\circ\Gamma_z\circ\pi^{n-2i}_Y$.
	Since $\pi^{n-2i}_Y$ is decomposable and since $\tau\circ\Gamma_z\circ(x\times x')=x\times\tau_*(z\cdot x')=0$ for  $x,x'\in\R^*(Y)$, we obtain $\Upsilon_z=0$. 
	Thus \eqref{eq:tau-contract-u} holds. 
	Arguing as in (i), the multiplicativity for $(i,j,k)=(n,n,n)$ gives $M(\tau,\tau;\tau)=0$, and then \eqref{eq:tau-small-diag-u} yields $\Theta_Y=0$.
\end{proof}

\subsection{The cohomological tautological ring}\label{SS:coh-taut-u}

Let $\bar\R^*(Y^m)\subseteq\HH^{2*}(Y^m,\QQ)$ be the image of $\R^*(Y^m)$ under the cycle class map. 
The algebra $\HH^{2*}_{\taut}(Y)=\bar\R^*(Y)$ is finite-dimensional, with $\HH^0_{\taut}(Y)=\QQ$ and $\HH^{2n}_{\taut}(Y)=\QQ\,[o]$, and the intersection pairing 
$\HH^{2*}_{\taut}(Y)\otimes\HH^{2*}_{\taut}(Y)\to\QQ$, $a\otimes a'\mapsto\deg(a\cdot a')$, where $\deg$ vanishes in degrees $\neq2n$, is perfect by (H2) and (H3). 
For $a\in\HH^*_{\taut}(Y)$ and $1\leq i\leq m$, write $a^{(i)}\coloneqq1\otimes\cdots\otimes a\otimes\cdots\otimes1\in\HH^*_{\taut}(Y)^{\otimes m}$, with $a$ in the $i$-th position, and let
\[
\TT_m\coloneqq\HH^{2*}_{\taut}(Y)^{\otimes m}[t_{jk}\mid 1\leq j<k\leq m],\qquad t_{jk}=t_{kj},
\]
be the graded polynomial algebra over $\HH^{2*}_{\taut}(Y)^{\otimes m}$ on formal variables $t_{jk}$ of degree $2n$. Since $\Delta_Y-\tau\in\mathrm{D}^n(Y\times Y)$, the algebra $\bar\R^*(Y^m)$ is generated by the classes $p_i^*a$, $a\in\HH^{2*}_{\taut}(Y)$, and $[\tau_{jk}]$, so that $a^{(i)}\mapsto p_i^*a$ and $t_{jk}\mapsto[\tau_{jk}]$ define a surjection $\TT_m\to\bar\R^*(Y^m)$. 
We put
\[
m_Y\coloneqq2d+1\ \text{ if $n$ is even},\qquad m_Y\coloneqq d+1\ \text{ if $n$ is odd},
\]
where $d\coloneqq\dim\HH^n_{\van}(Y)$ (which is even if $n$ is odd).

\begin{lem}\label{lem:relations-coh-u}
	The kernel of $\TT_m\to\bar\R^*(Y^m)$ is the ideal generated by the relations
	\begin{align}
	&t_{jk}\,a^{(j)}=0\ \text{for all}\ a\in\HH^{>0}_{\taut}(Y),\qquad t_{jk}^2=(-1)^nd\,[o]^{(j)}[o]^{(k)}\qquad\text{for }j\neq k,\label{eq:rel-Y2-u}\\
	&t_{ij}\,t_{ik}=t_{jk}\,[o]^{(i)}\qquad\text{for pairwise distinct }i,j,k,\label{eq:rel-Y3-u}
	\end{align}
	together with, in case $m>m_Y$, the relations obtained by substituting pairwise distinct elements of $\{1,\ldots,m\}$ for $1,\ldots,m_Y+1$ in
	\begin{equation}\label{eq:rel-kimura-u}
	\sum_{\sigma\in\Sy_{d+1}}\mathrm{sgn}(\sigma)\prod_{i=1}^{d+1}t_{i,\,d+1+\sigma(i)}=0\ \ \text{ if $n$ is even},\qquad
	\sum_{\varpi}\prod_{\{i,j\}\in\varpi}t_{ij}=0\ \ \text{ if $n$ is odd},
	\end{equation}
	where $\varpi$ runs through the partitions of $\{1,\ldots,d+2\}$ into pairs. In particular, for $m\leq m_Y$, the algebra $\bar\R^*(Y^m)$ is the quotient of $\TT_m$ by the relations \eqref{eq:rel-Y2-u} and \eqref{eq:rel-Y3-u}.
\end{lem}

\begin{proof} We first check that the relations \eqref{eq:rel-Y2-u}, \eqref{eq:rel-Y3-u} and \eqref{eq:rel-kimura-u}  do hold.
	The relations \eqref{eq:rel-Y2-u} and \eqref{eq:rel-Y3-u}  hold by \eqref{eq:van-kills-u} and by Proposition~\ref{prop:CK-basic-u}(ii) and (iii). If $n$ is even, \eqref{eq:rel-kimura-u} expresses the vanishing of $\bigwedge^{d+1}\HH^n_{\van}(Y)$. If $n$ is odd, the left-hand side of \eqref{eq:rel-kimura-u} maps to a class in $\HH^n_{\van}(Y)^{\otimes(d+2)}\subseteq\HH^*(Y^{d+2},\QQ)$ which is invariant under the permutations of the factors of~$Y^{d+2}$; by the Koszul sign rule, it is an element of $\bigwedge^{d+2}\HH^n_{\van}(Y)=0$.
	
	To see that these relations generate all relations, we follow \cite[\S3]{Yin} and \cite[Lemma~2.12]{FLV2}. 
	Choose a basis $\mathcal B$ of $\bigoplus_{0<j<n}\HH^{2j}_{\taut}(Y)$ consisting of homogeneous elements, such that every element of $\mathcal B$ has non-zero intersection with exactly one element of $\mathcal B$.
	Such a basis exists by (H2) and (H3), since the intersection pairing on $\HH^n_{\taut}(Y)$ diagonalizes if $n$ is even.
	Using \eqref{eq:rel-Y2-u} and \eqref{eq:rel-Y3-u}, every monomial of $\TT_m$ can be rewritten as a linear combination of monomials
	\begin{equation}\label{eq:tau-monomial-u}
	t_{I,\varpi}\cdot a_J\cdot o_L,\qquad
	t_{I,\varpi}\coloneqq\prod_{\{i,j\}\in\varpi}t_{ij},\quad a_J\coloneqq\prod_{j\in J}a_j^{(j)},\quad o_L\coloneqq\prod_{l\in L}[o]^{(l)},
	\end{equation}
	where $I,J,L\subseteq\{1,\ldots,m\}$ are pairwise disjoint, $\varpi$ is a partition of $I$ into pairs and $a_j\in\mathcal B$. 
	
	The algebra $\bar\R^*(Y^m)$ is graded by the cohomological degrees on the $m$ factors of $Y^m$; we call \emph{multidegree} the corresponding element of $\NN^m$, so that $a^{(i)}$ with $a\in\HH^j_{\taut}(Y)$ has multidegree $j$ at $i$, and $t_{ij}$ has multidegree $n$ at $i$ and at $j$. Note that $\bar\R^{nm}(Y^m)=\QQ\,o_1\cdots o_m$ has multidegree $(2n,\ldots,2n)$. Consider the pairing
	\begin{equation}\label{eq:taut-pairing-u}
	\bar\R^{i}(Y^m)\times\bar\R^{nm-i}(Y^m)\longrightarrow\bar\R^{nm}(Y^m)=\QQ\,o_1\cdots o_m.
	\end{equation}
	Let $t_{I,\varpi}a_Jo_L$ and $t_{I',\varpi'}a'_{J'}o_{L'}$ be two monomials \eqref{eq:tau-monomial-u} with non-zero product in $\bar\R^*(Y^m)$. Each index then receives total multidegree $2n$. An index $i\in I$ receives $n$ from the first monomial, and the missing $n$ can come neither from $[o]^{(i)}$, of multidegree $2n$, nor from a factor $a'^{(i)}_i$, which would be killed by \eqref{eq:rel-Y2-u}; hence $i\in I'$ and so $I=I'$. 
	It then follows that $J=J'$, $L'=(I\cup J\cup L)^\complement$, and, for each $j\in J$, that $a'_j$ is the unique element of $\mathcal B$ with $\deg(a_ja'_j)\neq0$.
	Therefore, after ordering the monomials \eqref{eq:tau-monomial-u} suitably, the matrix of \eqref{eq:taut-pairing-u} is block diagonal, the block attached to $(I,J,L,a_J)$ being a non-zero multiple of the matrix of the pairing $(t_{I,\varpi},t_{I,\varpi'})\mapsto t_{I,\varpi}\cdot t_{I,\varpi'}\cdot o_{I^\complement}$, where $\varpi,\varpi'$ run through the partitions of $I$ into pairs. By \eqref{eq:rel-Y2-u} and \eqref{eq:rel-Y3-u}, $t_{I,\varpi}t_{I,\varpi'}o_{I^\complement}=x^{c(\varpi,\varpi')}\,o_1\cdots o_m$, where $x\coloneqq(-1)^nd$ and $c(\varpi,\varpi')$ is the number of cycles of the graph on $I$ with edge set $\varpi\cup\varpi'$; this is the matrix $T_{|I|/2}(x)$ of Hanlon--Wales (see \cite[\S3]{Yin}). 
    Its eigen-spaces are indexed by the partitions $\lambda$ of $|I|$ into even parts, and for such a given partition $\lambda$, the corresponding eigenvalue is zero precisely when $\lambda$ has at least $d+1$ parts if $x=d$ \cite[Proposition~3.5]{Yin}, and precisely when $\lambda$ has a part of size at least $d+2$ if $x=-d$ \cite{Tavakol}. 
    Since such a partition $\lambda$ can only exist for $|I|>m_Y$, we see that for $m\leq m_Y$, the pairing \eqref{eq:taut-pairing-u} is perfect on the quotient of~$\TT_m$ by the ideal generated by \eqref{eq:rel-Y2-u} and \eqref{eq:rel-Y3-u}. The perfectness of the pairing excludes any further relation. For arbitrary $m$, the relations \eqref{eq:rel-kimura-u} is associated to exactly the degenerate eigenvalues, and we can conclude as in \cite[\S3]{Yin} and \cite{Tavakol}.
\end{proof}

\begin{prop}\label{prop:relations-generated-u}
	Let $G\subseteq\R^{>0}(Y)$ be a set of algebra generators of $\R^*(Y)$, and let $m\leq m_Y$. The kernel of the cycle class map on $\R^*(Y^m)$ is the ideal generated by the pullbacks, along the projections $Y^m\to Y$, $Y^m\to Y^2$ and $Y^m\to Y^3$, of the classes
	\[
	\kappa\in K,\qquad \Upsilon_g\ \ (g\in G),\qquad \tau^2-(-1)^nd\,o_1o_2,\qquad \Theta_Y,
	\]
	and of their transposes and permutations. 
\end{prop}

\begin{proof}
	Let $\TT'_m$ be the polynomial $\QQ$-algebra on symbols $g^{(i)}$ ($g\in G$, $1\leq i\leq m$) and $t_{jk}=t_{kj}$ ($j\neq k$), and let $\TT'_m\to\R^*(Y^m)$ be the homomorphism given by $g^{(i)}\mapsto p_i^*g$ and $t_{jk}\mapsto\tau_{jk}$, which is surjective
	since $p_{jk}^*\Delta_Y-\tau_{jk}=p_{jk}^*(\Delta_Y-\tau)$ is decomposable. Let $\fa$ be the kernel of the surjection $\QQ[G]\to\HH^*_{\taut}(Y)$; its image in $\R^*(Y)$ is $K$. 
	The composition $\TT'_m\to\R^*(Y^m)\to\bar\R^*(Y^m)$ factors through the surjection $\TT'_m\to\TT_m$, $g^{(i)}\mapsto[g]^{(i)}$, whose kernel is generated by the $\fa^{(i)}$. 
	Let $\omega\in\QQ[G]$ be a polynomial without constant term mapping to $o$. By Lemma~\ref{lem:relations-coh-u} and since $m\leq m_Y$, the kernel of $\TT'_m\to\bar\R^*(Y^m)$ is generated by the $\fa^{(i)}$, by the $t_{jk}g^{(j)}$ for $g\in G$
	(every element of $\HH^{>0}_{\taut}(Y)$ being the image of a polynomial without constant term in the elements of $G$), and by the $t_{jk}^2-(-1)^nd\,\omega^{(j)}\omega^{(k)}$ and $t_{ij}t_{ik}-t_{jk}\,\omega^{(i)}$. 
	The kernel of the cycle class map on $\R^*(Y^m)$ is the image of this ideal, hence is generated by the images of these generators, which are the classes listed.
\end{proof}

\begin{cor}\label{cor:relations-vanish-u}
	Assume that $K=0$, $\Upsilon_g=0$ for all $g\in G$, and $\Theta_Y=0$. Then the cycle class map is injective on $\R^*(Y^m)$ for all $m\leq m_Y$. If these classes vanish modulo algebraic equivalence, then the image of $\R^*(Y^m)$ in $\B^*(Y^m)$ injects into $\HH^{2*}(Y^m,\QQ)$ for all $m\leq m_Y$.
\end{cor}

\begin{proof}
	By Proposition~\ref{prop:relations-generated-u}, it suffices to show that  $\tau^2-(-1)^nd\,o_1o_2=0$.
	By Proposition~\ref{prop:MCK-criterion-u}(i), the canonical Chow--K\"unneth decomposition is multiplicative, and $\tau$ lies in $\CH^n_{(0)}(Y\times Y)$, as in the proof of Proposition~\ref{prop:MCK-criterion-u}(ii). Hence $\tau^2\in\CH^{2n}_{(0)}(Y\times Y)=(\pi^{2n}_Y\otimes\pi^{2n}_Y)_*\CH^{2n}(Y\times Y)=\QQ\,o_1o_2$, and $\tau^2=(-1)^nd\,o_1o_2$ because both sides have the same cohomology class. For the second assertion, $\B^*(Y^m)$ is a quotient ring of $\CH^*(Y^m)$ through which the cycle class map factors, so it suffices to check that $\tau^2-(-1)^nd\,o_1o_2$ is algebraically trivial, which is clear since a zero-cycle of degree zero on $Y\times Y$ is algebraically trivial.
\end{proof}

\begin{rmk}\label{rmk:kimura-u}
	For $m>m_Y$, one must moreover impose the relations \eqref{eq:rel-kimura-u}, with $\tau_{ij}$ in place of $t_{ij}$. If $n$ is even, they hold in $\CH^*(Y^{2d+2})$ if and only if $\bigwedge^{d+1}\h^n_{\van}(Y)=0$, that is, if and only if the motive of $Y$ is Kimura--O'Sullivan finite-dimensional \cite{Kimura}; compare \cite[Proposition~2.11]{FLV2}. If $n$ is odd, they hold in $\CH^*(Y^{d+2})$ as soon as $\Sym^{d+1}\h^n_{\van}(Y)=0$, that is, as soon as the motive of $Y$ is Kimura--O'Sullivan finite-dimensional. Indeed, the cycle $\sum_{\varpi}\prod_{\{i,j\}\in\varpi}\tau_{ij}\in\CH^*(Y^{d+2})$ is fixed by $\tau^{\otimes(d+2)}$ and invariant under the permutations of the factors of $Y^{d+2}$; it therefore defines a morphism from a Tate motive to $\Sym^{d+2}\h^n_{\van}(Y)$, which is a direct summand of $\Sym^{d+1}\h^n_{\van}(Y)\otimes\h^n_{\van}(Y)$.
\end{rmk}

\subsection{Double covers}\label{SS:double-u}
The following proposition will be applied to \emph{special} GM fourfolds, and to GM sixfolds.

\begin{prop}\label{prop:double-cover-u}
	Let $Y$ be a smooth projective variety of dimension $n$, let $p\colon Y\to P$ be a double cover of a smooth projective variety $P$ with trivial Chow groups, with covering involution $\iota$, and put $\R^*(Y)\coloneqq p^*\CH^*(P)$. Assume that the anti-invariant part of $\HH^*(Y,\QQ)$ is contained in $\HH^n(Y,\QQ)$. Then $(Y,\R^*(Y))$ satisfies (H1)--(H3), $\HH^*_{\taut}(Y)=\HH^*(Y,\QQ)^{\iota}$, and $\HH^n_{\van}(Y)$ is the anti-invariant part of $\HH^*(Y,\QQ)$. Moreover:
	\begin{enumerate}[(i)]
		\item $K=0$; thus $s$ is the inverse of $\cl\colon\R^*(Y)\xrightarrow{\sim}\HH^*_{\taut}(Y)$, and the canonical Chow--K\"unneth decomposition of $Y$ involves no choice.
		\item $\tau=\Delta^-_Y\coloneqq\tfrac12\bigl(\Delta_Y-\Gamma_\iota\bigr)$ in $\CH^n(Y\times Y)$.
		\item $M(\tau,\tau;\tau)=0$ in $\CH^{2n}(Y^3)$.
		\item	The canonical Chow--K\"unneth decomposition of $Y$ is multiplicative if and only if $\Upsilon_z=0$ for all $z\in\R^{>0}(Y)$. The same holds modulo algebraic equivalence.
	\end{enumerate}
\end{prop}

\begin{proof}
	First note that, since $P$ has trivial Chow groups,  it has zero odd cohomology. Hence $\HH^{2*}_{\taut}(Y)=p^*\HH^*(P,\QQ)=\HH^*(Y,\QQ)^\iota$. On the other hand, $\HH^j(Y,\QQ)=p^*\HH^j(P,\QQ)$ for $j\neq n$ by assumption. 
	This gives (H1) and (H2). 
	As $\iota^*$ preserves the intersection form, $\HH^n(Y,\QQ)$ is the orthogonal sum of its invariant part, which is $\HH^n_{\taut}(Y)$, and of its anti-invariant part. This gives (H3), and identifies $\HH^n_{\van}(Y)$ with the anti-invariant part of $\HH^n(Y,\QQ)$, hence of~$\HH^*(Y,\QQ)$.
	
	(i) Since $P$ has trivial Chow groups, this follows directly from the fact that $p_*p^*=2\,\id$.
	
	(ii) By assumption, the correspondence $\Delta^-_Y=\frac12(\id-\iota^*)$ is the projector onto $\HH^n_{\van}(Y)$, that is, $[\Delta^-_Y]=[\tau]$. The fibre product $Y\times_PY$ is the union of the graphs $\Delta_Y$ and $\Gamma_\iota$, each with multiplicity one, so that $\Delta^+_Y\coloneqq\frac12(\Delta_Y+\Gamma_\iota)=\frac12(p\times p)^*\Delta_P$ lies in $\mathrm{D}^n(Y\times Y)$. Therefore $\tau-\Delta^-_Y=\Delta^+_Y-(\Delta_Y-\tau)$ is a homologically trivial element of $\mathrm{D}^n(Y\times Y)$, and it vanishes by Lemma~\ref{lem:decomposable-u}, as $K=0$ by (i).
	
	(iii) For $\epsilon\in\{0,1\}^3$, the definition of $M$ and $\iota^{-1}=\iota$ give
	\[
	M(\Gamma_{\iota^{\epsilon_1}},\Gamma_{\iota^{\epsilon_2}};\Gamma_{\iota^{\epsilon_3}})=\bigl[\{(\iota^{\epsilon_1}y,\iota^{\epsilon_2}y,\iota^{\epsilon_3}y)\mid y\in Y\}\bigr],
	\]
	so that, $M$ being trilinear and $\tau=\frac12(\Delta_Y-\Gamma_\iota)$ by (ii),
	\[
	M(\tau,\tau;\tau)=\frac18\sum_{\epsilon\in\{0,1\}^3}(-1)^{\epsilon_1+\epsilon_2+\epsilon_3}\bigl[\{(\iota^{\epsilon_1}y,\iota^{\epsilon_2}y,\iota^{\epsilon_3}y)\mid y\in Y\}\bigr].
	\]
	Replacing $y$ by $\iota y$ shows that the cycle indexed by $\epsilon$ is the cycle indexed by $(1,1,1)-\epsilon$, while the two signs are opposite. Hence, the eight terms cancel in pairs and so $M(\tau,\tau;\tau)=0$.
	
	(iv) By (i), we have $K=0$.
	The vanishing of the $\Upsilon_z$ gives \eqref{eq:tau-contract-u}, hence \eqref{eq:tau-small-diag-u}, as in the proof of Proposition~\ref{prop:MCK-criterion-u}(i). By (iii), \eqref{eq:tau-small-diag-u} reads $\tau_{13}\tau_{23}=\tau_{12}o_3$, and permuting the factors gives $\Theta_Y=0$.
	We conclude with Proposition~\ref{prop:MCK-criterion-u}(i) that the canonical Chow--K\"unneth decomposition is multiplicative.  
	The converse is Proposition~\ref{prop:MCK-criterion-u}(ii).
	Finally, everything holds verbatim in $\B^*$.
\end{proof}

\subsection{Gushel--Mukai fourfolds}\label{SS:GM-u}

We begin with the following general lemma.

\begin{lem}\label{lem:divisor-u}
	Let $P$ be a smooth projective variety with trivial Chow groups and let $j\colon Y\hookrightarrow P$ be the zero locus of a section of a line bundle $L$ on~$P$. Assume that $Y$ has dimension $\dim P-1$ and that $j^*\CH^*(P)\subseteq\R^*(Y)$. Then $\Upsilon_{c_1(L|_Y)}$ lies in the ideal of $\R^*(Y\times Y)$ generated by $p_1^*K$ and $p_2^*K$. In particular, $\Upsilon_{c_1(L|_Y)}=0$ if $K=0$.
\end{lem}

\begin{proof}
	The square formed by the diagonal embeddings of $Y$ and $P$ and by $j$ and $j\times j$ is cartesian with excess bundle $N_{Y/P}=L|_Y$. 
	The excess intersection formula \cite[\S6.3]{Fulton} gives $(j\times j)^*\Delta_P=\Delta_Y\cdot p_1^*c_1(L|_Y)$. Since $P$ has trivial Chow groups, $\Delta_P\in\CH^*(P)\otimes\CH^*(P)$ \cite[Lemma~2.1]{FLV2}, hence $\Delta_Y\cdot p_1^*c_1(L|_Y)\in\mathrm{D}^{n+1}(Y\times Y)$ because $j^*\CH^*(P)\subseteq\R^*(Y)$. 
	Since $\Delta_Y-\tau$ is decomposable, the class $\Upsilon_{c_1(L|_Y)}=\Delta_Y\cdot p_1^*c_1(L|_Y)-(\Delta_Y-\tau)\cdot p_1^*c_1(L|_Y)$ is also decomposable. On the other hand, $\Upsilon_{c_1(L|_Y)}$ is homologically trivial by Proposition~\ref{prop:CK-basic-u}(iii). We conclude by Lemma~\ref{lem:decomposable-u}.
\end{proof}

Let $Y$ be a GM fourfold, and consider the tautological ring $\R^*(Y)=\langle h,c_2\rangle$ of Example~\ref{ex:R*forGM4}. 
By \cite[Proposition~3.1]{DK2}, the cohomology of $Y$ is concentrated in even degrees and the even Betti numbers of $Y$ are $1,1,24,1,1$. 
The class $\gamma_*[Y]$ is the class of a quadric section of a hyperplane section of $\Gr$ if $Y$ is ordinary, and twice the class of $Y_0$ if $Y$ is special.
In both cases $\gamma_*[Y]=2\sigma_1^2$, and the projection formula, together with $\sigma_1^6=5$, $\sigma_1^4\sigma_{1,1}=2$ and $\sigma_1^2\sigma_{1,1}^2=1$ on $\Gr$, gives $h^4=10$, $h^2c_2=4$ and $c_2^2=2$.
This shows that the pair $(Y,\R^*(Y))$ satisfies (H1)--(H3) with $d=22$, and that  $\R^1(Y)=\QQ\,h$ and $\R^2(Y)=\langle h^2,c_2\rangle$ inject into cohomology. Moreover, since $Y$ is Fano, $\R^4(Y) = \CH_0(Y)=\QQ$ injects in cohomology; in particular, $o$ is the class of a point.

\begin{defn}\label{defn:CK}
	The \emph{canonical Chow--K\"unneth decomposition} of a GM fourfold $Y$ is given by $o\coloneqq\frac1{10}h^4$ and
	\begin{equation*}
	\begin{gathered}
	\pi^0_Y\coloneqq o\times Y,\qquad \pi^2_Y\coloneqq\tfrac1{10}\,h^3\times h,\qquad
	\pi^6_Y\coloneqq\tfrac1{10}\,h\times h^3,\qquad \pi^8_Y\coloneqq Y\times o,\\
	\pi^4_{\taut}\coloneqq\tfrac12\,h^2\times h^2-h^2\times c_2-c_2\times h^2+\tfrac52\,c_2\times c_2,
	\end{gathered}
	\end{equation*}
	$\pi^k_Y\coloneqq0$ for odd $k$, $\tau\coloneqq\Delta_Y-\sum_{k\neq4}\pi^k_Y-\pi^4_{\taut}$ and $\pi^4_Y\coloneqq\pi^4_{\taut}+\tau$. We define $\h^k(Y)$, $\h^4_{\taut}(Y)$ and $\h^4_{\van}(Y)$ as in Definition~\ref{defn:CK-u}, and we put
	\[
	\Psi_Y\coloneqq h\,c_2-\tfrac25\,h^3\in\CH^3(Y),\qquad
	\Upsilon_Y\coloneqq \Upsilon_{c_2}=\tau\cdot p_1^*c_2\in\CH^6(Y\times Y).
	\]
	Note that $K=\QQ\,\Psi_Y$ and that
	Definition~\ref{defn:CK-u} with $s([h^3])=h^3$ gives the projectors of Definition~\ref{defn:CK}.
\end{defn}

\begin{prop}\label{prop:GM-instances}
	Let $Y$ be a GM fourfold.
	\begin{enumerate}[(i)]
		\item If $Y$ is special, then $K=0$.
		\item $\Upsilon_h=\tau\cdot h_1$ lies in the ideal generated by $p_1^*K$ and $p_2^*K$.
		In particular, if $Y$ is special, then $\Upsilon_h =0$.
	\end{enumerate}
	\end{prop}

\begin{proof}
	(i)  Looking at the Betti numbers of $Y$, we see that the anti-invariant part of $\HH^*(Y,\QQ)$ is contained in $\HH^4(Y,\QQ)$ and Lemma~\ref{lem:trivial-chow} gives that $\R^*(Y)=p^*\CH^*(Y_0)$, where $Y_0$ has trivial Chow groups. Therefore $K=0$ by Proposition~\ref{prop:double-cover-u}(i).
	
	(ii) We apply Lemma~\ref{lem:divisor-u}, with $c_1(L|_Y)=2h$. If $Y$ is ordinary, we take $P\coloneqq Y'=\Gr\cap H$, which is smooth by \cite[Proposition~2.22]{DK1}, has trivial Chow groups with $\CH^*(Y')=\langle h,c_2\rangle$ (Lemma~\ref{lem:trivial-chow}), and contains $Y=Y'\cap Q$ as the zero locus of a section of $L=\OO_{Y'}(2)$. If $Y$ is special, the fivefold $\CGr\cap H_1\cap H_2$ is the cone with vertex $O$ over $Y_0$, and we take $P$ to be its blow-up at $O$, which is the $\PP^1$-bundle $\PP_{Y_0}(\OO\oplus\OO(1))\to Y_0$, a smooth projective variety with trivial Chow groups. 
	Let $L_1$ be the pullback to $P$ of $\OO(1)$ from $\PP^{10}$. Since $O\notin Y$, we have $O\notin Q$, and $Y$ is the zero locus in $P$ of the pullback of the equation of~$Q$, which is a section of $L\coloneqq L_1^{\otimes2}$ that does not vanish along the exceptional divisor. By the projective bundle formula, $\CH^*(P)$ is generated as a module over the pullback of $\CH^*(Y_0)$ by $1$ and $c_1(L_1)$.
	Since $j^*L_1=\OO_Y(1)$ and the projection from $O$ restricts to $p$ on $Y$, we get $j^*\CH^*(P)=\langle h,c_2\rangle=\R^*(Y)$.
\end{proof}

\begin{prop}\label{prop:GM-criteria} Let $Y$ be a GM fourfold.
	\begin{enumerate}[(i)]
		\item The canonical Chow--K\"unneth decomposition of $Y$ is multiplicative if and only if $\Psi_Y=0$, $\Upsilon_Y=0$ and $\Theta_Y=0$.
		\item If $Y$ is special, it is multiplicative if and only if $\Upsilon_Y=0$. 
		\item For $m\leq45$, $\ker \bigl(\cl\colon \R^*(Y^m) \to \HH^*(Y^m,\QQ)\bigr)$ is the ideal generated by the pullbacks of $\Psi_Y$, $\Upsilon_Y$ and $\Theta_Y$ along all projections $Y^m\to Y^k$ composed with all permutations of the factors.
	\end{enumerate}
	Moreover, the statements hold verbatim modulo algebraic equivalence.
\end{prop}

\begin{proof}
	(i) If the decomposition is multiplicative, then, with the bigrading introduced in the proof of Proposition~\ref{prop:MCK-criterion-u}(ii), we have $h\in\CH^1_{(0)}(Y)$ and $c_2\in\CH^2_{(0)}(Y)$ by Definition~\ref{defn:CK}, and hence $hc_2\in\CH^3_{(0)}(Y)=(\pi^6_Y)_*\CH^3(Y)=\QQ\,h^3$. This gives $\Psi_Y=0$. 
	As $K=\QQ\,\Psi_Y$, Proposition~\ref{prop:MCK-criterion-u} with $G=\{h,c_2\}$ therefore shows that multiplicativity is equivalent to $\Psi_Y=0$, $\Upsilon_h=\Upsilon_{c_2}=0$ and $\Theta_Y=0$. On the other hand, $\Psi_Y=0$ implies $\Upsilon_h=0$ by Proposition~\ref{prop:GM-instances}(ii). 
	
	(ii) If $Y$ is special, Proposition~\ref{prop:double-cover-u}(iv), with $P=Y_0$, says that multiplicativity is equivalent to $\Upsilon_h=\Upsilon_{c_2}=0$, and $\Upsilon_h=0$ by Proposition~\ref{prop:GM-instances}(ii). 
	
	(iii) This is Proposition~\ref{prop:relations-generated-u}, with $G=\{h,c_2\}$ and $m_Y=2d+1=45$, together with $\CH_0(Y)=\QQ$, $K=\QQ \Psi_Y$ and Proposition~\ref{prop:GM-instances}(ii).
\end{proof}

\section{Lifting homological relations to algebraic equivalence}\label{S:algequiv}

Let $Y$ be a GM fourfold. By Proposition~\ref{prop:GM-criteria}(iii), the cycle class map is injective on the image of $\R^*(Y^m)$ in $\B^*(Y^m)$, for all $m\leq45$, as soon as the three classes
\[
\Psi_Y=h\,c_2-\tfrac25\,h^3\in\CH^3(Y),\quad \Upsilon_Y=\tau\cdot p_1^*c_2\in\CH^6(Y^2),\quad
\Theta_Y=\tau_{12}\tau_{13}-\tau_{23}\,o_1\in\CH^8(Y^3)
\]
vanish modulo algebraic equivalence. We prove the vanishing of the first two modulo algebraic equivalence for every GM fourfold, and of the third for the special ones.

\begin{lem}\label{lem:Psi-special}
	Let $Y$ be a special GM fourfold. Then the cycle class map is injective on $\GDCH^*_B(Y)=\R^*(Y)$. In particular, $\Psi_Y=0$ in $\CH^3(Y)$.
\end{lem}

\begin{proof}
	We have $\GDCH^*_B(Y)=\R^*(Y)$ by Lemma~\ref{lem:GDCH-Y}, and $K=0$ by Proposition~\ref{prop:GM-instances}(i). In particular, $\Psi_Y\in K$ vanishes.
\end{proof}

For $Y$ ordinary, the following proposition  is \cite[Lemmas~2.16 and~2.18]{LatBV}.
For a non-special GM fourfold, the vanishing of $\Psi_Y$ is not known, but it holds modulo algebraic equivalence:

\begin{prop}\label{prop:alg-Y2}
	Let $Y$ be a GM fourfold. Then $\B^3(Y)=\QQ\,h^3$, and
	\[
	\Psi_Y=0\ \text{ in }\ \B^3(Y),\qquad \Upsilon_Y=0\ \text{ in }\ \B^6(Y\times Y).
	\]
	Consequently, the cycle class map is injective on $\GDB^*_B(Y^m)$ for $m\leq2$.
\end{prop}

\begin{proof}
	Since $\CH_0(Y)=\QQ$, the effective motive $(Y,\Delta_Y-\pi^0_Y-\pi^8_Y)$ has trivial Chow group of zero-cycles.
	Hence \cite[Corollary~2.2]{VialCK}, which relies on the Bloch--Srinivas argument \cite{BlochSrinivas}, provides an isomorphism of Chow motives
	\[
	\h(Y)\ \cong\ \1\oplus (S,\pi,-1)\oplus\1(-4)\qquad\text{in }\Mrat,
	\]
	for some smooth projective surface $S$ and some idempotent correspondence $\pi \in \CH^2(S\times S)$. Hence $\B^3(Y) = \pi_*\B^2(S)$, and since the cycle class map is injective on $\B^2(S)$, it is also injective on $\B^3(Y)$. Using that $\HH^6(Y,\QQ)=\QQ\,h^3$ yields $\B^3(Y)=\QQ\,h^3$.  This establishes that $\Psi_Y$ vanishes in $\B^3(Y)$. 
	Now $\h(Y\times Y)=\h(Y)\otimes \h(Y)$ is isomorphic to a direct summand of a direct sum of $\h(S\times S)(-2)$ together with Tate twists of $\h(S)$ and of $\1$, so that $\B^6(Y\times Y)$ is motivically a direct summand of a direct sum of $\B^4(S\times S)$  with copies of $\B^*(S)$ and of~$\B^0(\mathrm{pt})$. 
Since the cycle class map is injective on $\B^4(S\times S)$, $\B^*(S)$ and $\B^0(\mathrm{pt})$ and since $\Upsilon_Y$ is homologically trivial (Proposition~\ref{prop:CK-basic-u}(iii)), we find that $\Upsilon_Y=0$ in $\B^6(Y\times Y)$.

	For the last assertion, Lemma~\ref{lem:GDCH-Y} and Proposition~\ref{prop:generation} give $\GDCH^*_B(Y^m)=\R^*(Y^m)$ for $m\leq2$, so that $\GDB^*_B(Y^m)$ is the image of $\R^*(Y^m)$ in $\B^*(Y^m)$. By Proposition~\ref{prop:GM-criteria}(iii), for $m\leq2$, the kernel of the cycle class map on $\R^*(Y^m)$ is the ideal generated by the pullbacks of $\Psi_Y$, of $\Upsilon_Y$ and of its transpose, and these vanish in $\B^*$.
\end{proof}

\begin{rmk}
	Alternatively, the vanishing of $\Upsilon_Y$ in $\B^6(Y\times Y)$ can be proven without the Bloch--Srinivas argument, by observing that $c_2$ is represented by a del Pezzo surface in $Y$.
\end{rmk}

\begin{prop}\label{prop:MCK-alg}
	Let $Y$ be a special GM fourfold, with covering involution $\iota_Y$. Then
	\[
	\tau=\tfrac12\bigl(\Delta_Y-\Gamma_{\iota_Y}\bigr)\ \text{ in }\CH^4(Y\times Y),
	\qquad\text{and}\qquad
	\tau_{12}\,\tau_{13}=\tau_{23}\,o_1\ \text{ in }\B^8(Y^3).
	\]
	In particular, the canonical Chow--K\"unneth decomposition of $Y$ is multiplicative modulo algebraic equivalence.
\end{prop}

\begin{proof}
	By the proof of Proposition~\ref{prop:GM-instances}(i), Proposition~\ref{prop:double-cover-u} applies to $p\colon Y\to Y_0$, with $\R^*(Y)=\langle h,c_2\rangle=p^*\CH^*(Y_0)$ and $G=\{h,c_2\}$.
	Specifically, Proposition~\ref{prop:double-cover-u}(ii) gives $\tau=\tfrac12\bigl(\Delta_Y-\Gamma_{\iota_Y}\bigr)$ in $\CH^4(Y\times Y)$.
	On the other hand, $\Upsilon_h=0$ by Proposition~\ref{prop:GM-instances}(ii), and $\Upsilon_{c_2}=\Upsilon_Y=0$ in $\B^6(Y\times Y)$ by Proposition~\ref{prop:alg-Y2}, so that Proposition~\ref{prop:double-cover-u}(iv), modulo algebraic equivalence, gives $\Theta_Y\coloneqq \tau_{12}\,\tau_{13}-\tau_{23}\,o_1 =0$ in $\B^8(Y^3)$ and the multiplicativity modulo algebraic equivalence.
\end{proof}

\begin{cor}\label{cor:Y4-alg}
	Let $Y$ be a special GM fourfold. Then, for every $m\leq45$, the image of $\R^*(Y^m)$ in $\B^*(Y^m)$ injects into $\HH^{2*}(Y^m,\QQ)$. In particular, the cycle class map is injective on $\GDB^*_B(Y^m)$ for $m\leq3$ and on $\GDB^i_B(Y^4)$ for $i\leq8$.
\end{cor}

\begin{proof}
	The classes $\Psi_Y$, $\Upsilon_Y$ and $\Theta_Y$ vanish modulo algebraic equivalence by Lemma~\ref{lem:Psi-special} and Propositions~\ref{prop:alg-Y2} and~\ref{prop:MCK-alg}; the first assertion thus follows from Proposition~\ref{prop:GM-criteria}(iii). By Proposition~\ref{prop:generation}, one has $\GDCH^*_B(Y^m)=\R^*(Y^m)$ for $m\leq3$ and $\GDCH^i_B(Y^4)=\R^i(Y^4)$ for $i\leq8$, so that $\GDB^*_B(Y^m)$, resp.\ $\GDB^i_B(Y^4)$, is the image of $\R^*(Y^m)$, resp.\ of $\R^i(Y^4)$, in $\B^*(Y^m)$.
\end{proof}

\section{The MCK decomposition for special Gushel--Mukai fourfolds}
\label{S:MCK}

In this section we prove Theorems \ref{thm:MCK-intro} and \ref{thm:Y3-intro} of the introduction. A main ingredient in proving Theorem \ref{thm:MCK-intro} is Theorem \ref{thm:zero-cycles-X}, stating that $\GDCH^4_B(X)=\QQ$ for any double EPW sextic $X$. Note that this result implies the codimension 4 part of Theorem \ref{thm:main}(i). To prove Theorem \ref{thm:zero-cycles-X}, we use the Iliev--Manivel description of the general $X$ in terms of conics in a GM fourfold.

\subsection{Zero-cycles on the associated double EPW sextic}

Let $Y$ be a GM fourfold,
and let
$
M_Y\coloneqq\CGr\cap H_1\cap H_2$
be its Grassmannian hull (for $Y$ ordinary, $M_Y$ is mapped
isomorphically onto the del Pezzo fivefold $Y'$ of Section~\ref{S:tautological} by the projection away from
the vertex of $\CGr)$.
We write $F(Y)$ for the Hilbert scheme of conics
contained in $Y$ and $\langle\Sigma\rangle\cong\PP^2$ for the linear span of a conic
$\Sigma\subseteq Y$. Following \cite[\S3.1]{IM} we set
\[
F_{\rho\sigma}(Y)\coloneqq\{[\Sigma]\in F(Y)\ \vert\ \langle\Sigma\rangle\subseteq M_Y\}=F_\rho(Y)\sqcup F_\sigma(Y),
\qquad F_\tau(Y)\coloneqq F(Y)\setminus F_{\rho\sigma}(Y),
\]
where $F_\tau(Y)$, $F_\rho(Y)$ and $F_\sigma(Y)$ correspond to $\tau$-conics resp. $\rho$-conics resp. $\sigma$-conics.
Likewise $F(\Gr)$ denotes the Hilbert scheme of conics in $\Gr$, a smooth projective spherical
variety of dimension $13$ \cite[\S3.1]{IM}, \cite[Lemma 3.12]{DK5}, with its loci
$F_{\rho\sigma}(\Gr)$ and~$F_\tau(\Gr)$.

Let $b\in B^\circ$, let $(V_6,V_5,A)$ be the Lagrangian datum of $Y_b$ and let
$p_b\coloneqq[V_5]\in\PP(V_6^\vee)$ be its Pl\"ucker point, so that $Z_b=Y_{A^\perp}$ and
$X_b=\widetilde Y_{A^\perp}$ in the notation of \cite[\S2.2]{DK2}. One has $p_b\in Y^\ell_{A^\perp}$ with $\ell=1$ if $Y_b$ is ordinary, $\ell=2$ if $Y_b$ is special
\cite[Proposition 3.13]{DK1}. 

We let $B^{\circ \circ}\subseteq B^\circ$ be the subset of those $b$ for which $Y_b$ is ordinary and
contains no plane. It is open (because the ordinary locus is open in $B$, and the locus of those $b\in B$
for which $Y_b$ contains a plane is closed, being the image of the incidence variety of pairs
$(\Pi,b)$ with $\Pi\subseteq Y_b$ a plane, which is proper over $B$) and dense (because the ordinary locus is dense and a general GM fourfold contains no plane \cite[Proposition 3.7]{DK5}).

\begin{thm}[Debarre--Kuznetsov]\label{thm:FX} Let $Y=Y_b$ and $X=X_b$, where $b\in  B^{\circ \circ}$.
The scheme $F(Y)$ is smooth irreducible of dimension~$5$, and 
there is a surjective morphism $\pi\colon F(Y)\to X$ which contracts $F_\rho(Y)$ and
$F_\sigma(Y)$ to two distinct points $x_\rho(b),x_\sigma(b)$, the two preimages under $f\colon X\to Z$ of the
Pl\"ucker point $p_b$. The morphism $\pi$ restricts to an \'etale-locally trivial $\PP^1$-fibration
\[
\pi_\tau\colon F_\tau(Y)\longrightarrow X^\tau\coloneqq X\setminus f^{-1}(p_b).
\]
\end{thm}

\begin{proof} This follows from \cite[Theorem 7.12]{DK5}, which makes more precise the generality conditions in the earlier result of Iliev--Manivel \cite{IM}.
The two points $x_\rho(b), x_\sigma(b)$ are distinct because $p_b\in Y^1_{A^\perp}$, hence $p_b\notin
Y^{\geq2}_{A^\perp}=\mathrm{Sing}(Z)$. 
\end{proof}

\begin{rmk}\label{rmk:IM-in-families}
The construction of Theorem \ref{thm:FX} makes sense relatively over $B^{\circ \circ}$. Indeed, the relative Hilbert scheme of conics $\FF\to B^{\circ \circ}$
is an open subscheme of the projective bundle~$\bar\FF$ of the proof of
Lemma~\ref{lem:conics-generation} below, hence is smooth of dimension $\dim B^{\circ \circ}+5$; as its
fibres are smooth of dimension~$5$, it is smooth over $B^{\circ \circ}$. The Pl\"ucker points $p_b$ form a
section of the family of EPW sextics $\XX_B/\iota\to B^\circ$, and we let $P\subseteq\XX_B$ be the
reduced preimage of that section under $f$ and $\XX^\tau_B\coloneqq\XX_B\smallsetminus P$; thus $P$ is
finite over $B^\circ$, of codimension $4$ in $\XX_B$, and $P\vert_{B^{\circ \circ}}\to B^{\circ \circ}$ is \'etale
of degree $2$, the double cover $f$ being \'etale away from $\mathrm{Sing}(Z_b)$. 

The morphism $\pi$
comes from the family of quadrics that the GM datum of $Y_b$ determines over the blow-up of
$\PP(V_6^\vee)$ at $p_b$ \cite[\S4]{DK5}: the Stein factorization of the associated relative
orthogonal Grassmannian of planes is the double covering $f\colon X_b\to Z_b$ blown up at the two
points over $p_b$ \cite[Lemma 5.5 and Proposition 5.6]{DK5}, and it descends to $\pi$ along the
birational morphisms of \cite[Theorem 7.12]{DK5}. 

Each of these steps commutes with base change, the
double covering attached to a Lagrangian datum being constructed over an arbitrary base
\cite[Theorems 4.2 and 5.2]{DKdouble} and agreeing there with the Stein factorization of the
corresponding relative orthogonal Grassmannian \cite[Proposition 3.10]{DKdouble}. Performing the construction for
$\YY\vert_{B^{\circ \circ}}\to B^{\circ \circ}$ therefore exhibits the morphism $\pi$ of Theorem \ref{thm:FX} as the fiberwise restriction of
a morphism of $B^{\circ \circ}$-schemes $\pi\colon\FF\to\XX_B\vert_{B^{\circ \circ}}$ contracting the relative
locus $\FF_{\rho\sigma}$ onto~$P\vert_{B^{\circ \circ}}$, and
$\pi_\tau\colon\FF_\tau\to\XX^\tau_B\vert_{B^{\circ \circ}}$ is an \'etale-locally trivial
$\PP^1$-fibration. 
\end{rmk}

\begin{lem}\label{lem:conics-Gr}
	The cycle class map $\CH^*(F(\Gr))\to\HH^{2*}(F(\Gr),\QQ)$ is an isomorphism, and there are classes
	$D,D'\in\CH^1(F(\Gr))$, $d_2\in\CH^2(F(\Gr))$ and $d_3\in\CH^3(F(\Gr))$ such that
	\[
	\CH^*(F(\Gr))=\langle D,D',d_2,d_3\rangle+\ima\Bigl(\CH_*(F_{\rho\sigma}(\Gr))\to\CH_*(F(\Gr))\Bigr),
	\]
	where $\langle\ \rangle$ denotes the $\QQ$-subalgebra generated by the indicated classes.
\end{lem}

\begin{proof}
	The variety $F(\Gr)$ is a smooth projective spherical variety \cite[\S3.1]{IM}, hence admits a
	decomposition into finitely many affine cells; this implies $F(\Gr)$ has trivial Chow groups. As shown in \cite[\S3.1]{IM}, there are birational maps
	\[
	S(\Gr)\ \xleftarrow{\ \phi_S\ }\ \widehat F(\Gr)\ \xrightarrow{\ \phi_F\ }\ F(\Gr),
	\]
	where $\widehat F(\Gr)$ is the nested Hilbert scheme of pairs $(\Sigma,V_4)$ with $\Sigma$ a conic
	contained in~$\Gr(2,V_4)$, the morphism $\phi_F$ is the blow-up of $F_\rho(\Gr)$, the variety
	$S(\Gr)$ parametrizes pairs $(P,V_4)$ with $P\subseteq\PP(\wedge^2V_4)$ a plane, and $\phi_S$
	contracts two disjoint divisors onto two subvarieties of codimension~$6$. The morphism
	$p\colon S(\Gr)\to\PP(V_5^\vee)$, $(P,V_4)\mapsto[V_4]$, is a $\Gr(3,\wedge^2V_4)$-bundle, so
	$S(\Gr)$ has trivial Chow groups; and since $\CH^*(\Gr(3,6))$ is generated by the Chern classes of
	the tautological subbundle \cite[Theorem 5.26]{EH}, the Chow ring of $S(\Gr)$ is generated by
	$p^*c_1(\OO(1))$ and by the Chern classes $d_1,d_2,d_3$ of the relative tautological subbundle
	of~$p$. As $\phi_S$ and $\phi_F$ restrict to isomorphisms on the complement of $F_{\rho\sigma}(\Gr)$, the images under $\phi_{F*}\phi_S^*$ of these generators, together with cycles supported on
	$F_{\rho\sigma}(\Gr)$, generate $\CH^*(F(\Gr))$. Taking for $D,D'$ the two divisor classes thus
	obtained proves the statement.
\end{proof}

An ordinary GM fourfold $Y_b$ is (isomorphic to) a codimension-2 section in the Grassmannian $\Gr$, meaning that there is a closed embedding 
$\gamma\colon F(Y_b)\hookrightarrow F(\Gr)$.

\begin{lem}\label{lem:conics-generation}
	Let $b\in B$ be such that $Y_b$ is an ordinary GM fourfold. Then
	$$\GDCH^*_B\bigl(F(Y_b)\bigr)=\ima\Bigl(\gamma^*\colon\CH^*(F(\Gr))\to\CH^*(F(Y_b))\Bigr).
	$$
\end{lem}

\begin{proof}
		Let $\bar\FF\coloneqq\{([\Sigma],b)\in F(\Gr)\times\bar B\ \vert\ \Sigma\subseteq Y_b\}$ and let
	$\bar p\colon\bar\FF\to F(\Gr)$ and $\bar q\colon\bar\FF\to\bar B$ be the projections. Since $\Sigma$ spans a plane,
	the restriction maps $\HH^0(\PP^{10},\OO(j))\to\HH^0(\Sigma,\OO_\Sigma(j))$ are surjective for
	$j=1,2$, with targets of dimension $3$ and~$5$, and $\HH^0(\PP^{10},\OO(j))\to\HH^0(F(\Gr),\OO(j))$ is
	surjective as well.
    The condition $\Sigma\subseteq H_1\cap Q$ therefore consists of $8$
	linear conditions of constant rank on $\bar B$, and $\bar p$ is a projective bundle. Hence
	$\bar\FF$ is smooth of dimension $\dim\bar B+5$. 
    Let $\xi\coloneqq c_1(\OO_{\bar p}(1))$, so that the projective bundle formula yields
	$$
	\CH^*(\bar\FF)=\bigoplus_{k\geq0}\bar p^*\CH^{*-k}(F(\Gr))\cdot\xi^k .
	$$
	Let $h_{\bar B}\in\CH^1(\bar B)$ be the hyperplane class and write $\bar q^*h_{\bar
		B}=\nu\xi+\bar p^*\beta$ with $\nu\in\QQ$ and $\beta\in\CH^1(F(\Gr))$. The scalar $\nu$ cannot vanish, as it would imply that $\bar q^*(h_{\bar B}^{\dim\bar B})=\bar p^*(\beta^{\dim\bar B})$, of which the right-hand side vanishes because
	$\dim\bar B>13=\dim F(\Gr)$, whereas the left-hand side is the class of a fibre of $\bar q$ and is
	non-zero. We can thus write $\xi=\nu^{-1}(\bar q^*h_{\bar B}-\bar p^*\beta)$. Now let
	$\alpha\in\CH^*(\bar\FF)$ and write $\alpha=\sum_k\bar p^*(z_k)\xi^k$ with
	$z_k\in\CH^*(F(\Gr))$. Restricting to the fibre $F(Y_b)$ of $\bar q$ sends $\bar q^*h_{\bar B}$ to zero. Hence
	$$
	\alpha\vert_{F(Y_b)}=\sum_k(-\nu^{-1})^k\,\bar p^*(z_k\beta^k)\vert_{F(Y_b)}
	$$
	is the restriction of a class pulled back from $F(\Gr)$. 
    This proves the inclusion ``$\subseteq$''.
    The reverse inclusion is clear.
\end{proof}

\begin{prop}\label{prop:conics-generators}
	Let $b\in B^{\circ \circ}$ and write $Y=Y_b$, $X=X_b$, $F_\tau=F_\tau(Y)$ and $\pi=\pi_\tau$. Let
	$\xi\in\CH^1(F_\tau)$ be the restriction of an ample class of $F(\Gr)$. Then $\xi$ has positive degree on the fibres of $\pi$, and there are classes $b_2,g_2\in\GDCH^2_B(X)$ and
	$b_3\in\GDCH^3_B(X)$ such that
	\[
	\GDCH^*_B(F_\tau)=\bigl\langle\,\pi^*h_X,\ \xi,\ \pi^*b_2,\ \pi^*g_2,\ \pi^*b_3\,\bigr\rangle .
	\]
\end{prop}

\begin{proof}
    Clearly the restriction of an ample class of $F(\Gr)$ to $F(Y)$ is ample and has positive degree on any fibre of $\pi$; the same then holds after restriction to the open $F_\tau=\pi^{-1}(X^\tau)\subset F(Y)$.
	
	We know from Lemmas~\ref{lem:conics-generation} and~\ref{lem:conics-Gr} (plus the trivial observation that classes supported on
	$F_{\rho\sigma}(\Gr)$ restrict to zero on $F_\tau$) that the ring
	$\GDCH^*_B(F_\tau)$ is generated by the restrictions of $D,D',d_2,d_3$. In particular
	$\GDCH^1_B(F_\tau)$ has dimension at most~$2$. It contains $\xi$ and~$\pi^*h_X$, which are linearly
	independent (indeed, $\pi^*h_X\neq0$ has degree $0$ on the fibres of $\pi$ while $\xi$ has positive
	degree), and so $\GDCH^1_B(F_\tau)=\langle\pi^*h_X,\xi\rangle$.
	
	Since $\pi$ is an \'etale-locally trivial $\PP^1$-fibration and $\xi$ has non-zero degree on its
	fibres, the projective bundle formula (with $\QQ$-coefficients) gives
	\begin{equation}\label{eq:pbf}
	\CH^*(F_\tau)=\pi^*\CH^*(X^\tau)\oplus\pi^*\CH^{*-1}(X^\tau)\cdot\xi,
	\end{equation}
and this holds relatively over $B^{\circ \circ}$, so that the two components of
	a generically defined class are generically defined. Writing
	$d_2\vert_{F_\tau}=\pi^*(b_2)+\pi^*(f_1)\,\xi$ and $d_3\vert_{F_\tau}=\pi^*(b_3)+\pi^*(g_2)\,\xi$ and
	using that $\CH^i(X)\to\CH^i(X^\tau)$ is bijective for $i\leq3$, we find classes $b_2,g_2\in\GDCH^2_B(X)$,
	$f_1\in\GDCH^1_B(X)$ and $b_3\in\GDCH^3_B(X)$. As the very general
	double EPW sextic has Picard number one, we can replace $f_1$ by $h_X$. This gives the asserted description.
\end{proof}

\begin{lem}\label{lem:conics-eigen}
	Let $b\in B^\circ$ be such that $Y_b$ is a special GM fourfold. Then $\GDCH^i_B(X_b)\subseteq\CH^i(X_b)^+$ for
	all $i$, where $\CH^*(X_b)^\pm$ are the eigenspaces of the covering involution $\iota$.
\end{lem}

\begin{proof}
	We know from \cite[Proposition 4.4]{BL2} that there exist correspondences inducing a split injection
	$\CH^*(X_b)^-\hookrightarrow\bigoplus\CH^*(Y_b)\oplus\QQ^r$, compatible with the cycle class maps.
	These correspondences, being defined relatively over a dense open subset of $B$, respect generically defined cycles (with
	respect to $B$), and so we obtain a split injection
	\[
	\GDCH^*_B(X_b)^-\cap\CH^*_{\mathrm{hom}}(X_b)\hookrightarrow
	\bigoplus\GDCH^*_B(Y_b)\cap\CH^*_{\mathrm{hom}}(Y_b).
	\]
	(Here the notation $\GDCH^*_B(X_b)^-$ makes sense because the involution $\iota$ acts on the family $\XX_B\to B^\circ$, and so
	 $\frac12(\id-\iota^*)$ preserves generically defined cycles.) 
    The GM fourfold $Y_b$ being special by assumption, Lemma~\ref{lem:Psi-special}
	implies that the right-hand side vanishes, and so
	\begin{equation}\label{eq:eigen-hom}
	\GDCH^*_B(X_b)\cap\CH^*_{\mathrm{hom}}(X_b)\subseteq\CH^*(X_b)^+ .
	\end{equation}
	
	Let us now consider a class $a\in\GDCH^i_B(X_b)$. Its cohomology class is invariant under the monodromy of the family
	$\XX_B\to B^\circ$, which is the monodromy of the locally complete family of double EPW sextics,
	since $B^\circ$ dominates $\MM_{\mathrm{EPW}}$ with connected fibres. The monodromy invariant
	subspace of $\HH^{2i}(X_b,\QQ)$ is spanned by the monomials of weight $i$ in $h_X$ and $c_2(X_b)$:
	indeed $\HH^4(X_b,\QQ)=\Sym^2\HH^2(X_b,\QQ)$ \cite{Verbitsky}, \cite{Bogomolov}, the monodromy
	group acts on $\HH^2(X_b,\QQ)$ through a subgroup of finite index of the orthogonal group of
	$h_X^\perp$ \cite{MarkmanSurvey}, so that its invariants in $\Sym^2\HH^2(X_b,\QQ)$ are spanned by
	$h_X^2$ and by the class of the Beauville--Bogomolov form, which is a combination of $h_X^2$ and
	$c_2(X_b)$; and $\HH^6$, $\HH^8$ are dual to $\HH^2$, $\HH^0$. These monodromy-invariant monomials are $\iota^*$-invariant, because $h_X$ and~$c_2(X_b)$ are. It follows that $\frac12(a-\iota^*a)$ is generically defined and homologically
	trivial, so it vanishes by \eqref{eq:eigen-hom} and hence $a\in\CH^i(X_b)^+$.
\end{proof}

\begin{lem}\label{lem:conics-Z}
	Let $Z\subseteq\PP(V_6^\vee)$ be an EPW sextic with only quotient singularities and let $h_Z$ be its
	hyperplane class. Then $\CH^i(Z)\cdot h_Z=\QQ\,h_Z^{i+1}$ for all $i$, and
	$\CH^i(Z)\cdot\CH^{4-i}(Z)=\QQ\,h_Z^4$ for $0<i<4$.
\end{lem}

\begin{proof}
	The first assertion is elementary and holds for any hypersurface with quotient singularities: letting
	$\varphi\colon Z\hookrightarrow\PP(V_6^\vee)$ denote the inclusion, we have
	\[ \CH^i(Z)\cdot h_Z=\varphi^*\varphi_*\CH^i(Z)\subseteq\varphi^*\CH^{i+1}(\PP(V_6^\vee))=\QQ h_Z^{i+1}.\]
	The second is \cite[Theorem 1.2]{LV}, which extends to Calabi--Yau hypersurfaces with quotient
	singularities the theorem of Voisin \cite[\S3]{VoisinCY} for smooth ones.
\end{proof}

\begin{lem}\label{lem:conics-S}
	Let $X$ be a double EPW sextic, with covering involution $\iota$ and fixed surface $S$.
	Then $\ima\bigl(\CH_0(S)\to\CH^4(X)\bigr)=\QQ\,h_X^4$.
\end{lem}

\begin{proof}
	The surface $S$ is a constant cycle surface by \cite[Theorem 1.8]{LYZ} (or, for $X$ general, \cite[Theorem 4.14]{Zhang}), so that all its points
	have the same class $c\in\CH^4(X)$. 
    Furthermore, it is the
	ramification locus of $f\colon X\to Z$, so that
	$f^*f_\ast[S]=2[S]$ in $\CH^2(X)$.
    Therefore
	$$
	[S]\cdot h_X^2=\tfrac12\,f^*\bigl(f_\ast[S]\cdot h_Z^2\bigr)\in\QQ\,f^*(h_Z^4)=\QQ\,h_X^4
	$$
	by Lemma~\ref{lem:conics-Z}. The class $[S]\cdot h_X^2$ is also represented by a zero-cycle
	supported on $S$ of degree $m>0$, so that $[S]\cdot h_X^2=mc$. Hence
	$c\in\QQ\,h_X^4$, completing the proof.
\end{proof}

We are now ready to prove the main result of this subsection.

\begin{thm}\label{thm:zero-cycles-X}
	Let $b\in B^\circ$ be such that $Y_b$ is a special GM fourfold, and let $X=X_b$ be the associated double EPW
	sextic. Then $\GDCH^4_B(X)=\QQ h_X^4$.
    In particular, the cycle class map is injective on
	$\GDCH^4_B(X)$.
\end{thm}

\begin{proof} We proceed in 3 steps.

	\emph{Step 1: relative classes and an identity over $B^{\circ \circ}$.}
	Let $H_X\in\CH^1(\XX_B)$ be the relative Pl\"ucker class. 
    We use the notation from 
	Proposition~\ref{prop:conics-generators}.
    For $b\in B^{\circ \circ}$, define $e\in\GDCH^2_B(X_b)$ and
	$f\in\GDCH^1_B(X_b)$ by 
    $$\xi^2=\pi^\ast(e|_{X_b^\tau})+\pi^\ast(f|_{X_b^\tau})\xi \text{ in }\CH^2(F_{\tau,b}).$$ Let $\eta_k\in\GDCH^{k-1}_B(X_b)$ be the unique class with
	$$\eta_k\vert_{X^\tau_b}=\pi_*(\xi^k)\text{ in } \CH^{k-1}(X_b^\tau),$$
    where $2\leq k\leq 4$.
    Such a class exists because $\CH^i(X_b)\to\CH^i(X^\tau_b)$ is an isomorphism
	for~$i\leq3$, since the complement of $X^\tau_b$ is finite. All these classes are defined relatively
	over~$B^{\circ \circ}$.
    Indeed, the morphism $\pi_\tau$ and the class $\xi$ are, and
	$\CH^i(\XX_B\vert_{B^{\circ \circ}})\to\CH^i(\XX^\tau_B\vert_{B^{\circ \circ}})$ is bijective for $i\leq3$ because
	$P$ has codimension~$4$. Taking closures and using the surjectivity of
	$\CH^*(\XX_B)\to\CH^*(\XX_B\vert_{B^{\circ \circ}})$, they are therefore the restrictions to the fibres over
	$B^{\circ \circ}$ of classes
	\[
	E_2\in\CH^1(\XX_B),\qquad B_2,G_2,E,E_3\in\CH^2(\XX_B),\qquad B_3,E_4\in\CH^3(\XX_B)
	\]
	restricting to $\eta_2$, to $b_2,g_2,e,\eta_3$ and to $b_3,\eta_4$ respectively. Let $L$ be the
	finite list of the products $U\cdot V$ of total codimension~$4$ of two monomials $U,V$ of positive
	codimension in $H_X,E_2,B_2,G_2,E,E_3,B_3,E_4$; note that $h_X^4=(H_X^2\cdot H_X^2)\vert_{X_b}$.
	
	We claim that for every $b\in B^{\circ \circ}$, every class $a\in\GDCH^4_B(X_b)$  lies in the $\QQ$-span of the restrictions
	to $X_b$ of the elements of $L$ and of the classes $[x_\rho(b)]$, $[x_\sigma(b)]$.
	
	Indeed, since $\pi_*\xi=d\,[X^\tau_b]$ for some non-zero $d\in\QQ$, the projection formula gives
	\[
	 a\vert_{X^\tau_b}={1\over d} \, \pi_*\bigl(\pi^*(a\vert_{X^\tau_b})\cdot\xi\bigr) \ \ \hbox{in}\ \CH^{4}(X^\tau_b).
	\]
	 The class $\pi^*(a\vert_{X^\tau_b})\cdot\xi$ lies in $\GDCH^5_B(F_\tau)$, hence by
	Proposition~\ref{prop:conics-generators} it is a $\QQ$-linear combination of monomials
	$\pi^*(u)\cdot\xi^k$, where $u$ is a monomial of degree $5-k$ in $h_X,b_2,g_2,b_3$ and $k\ge 1$ because $\CH^5(X^\tau_b)=0$. 
    
It follows that the restriction $a\vert_{X^\tau_b}$ lies in the $\QQ$-linear span of the finitely many elements
     \[ \pi_*(\pi^*(u)\cdot\xi^k)=u\cdot\pi_*(\xi^k)\ \in \CH^{4}(X^\tau_b).\] 

     Let us now assume that $1\leq k\leq4$.
     The cycle $u\cdot \pi_\ast(\xi^k) $ is the restriction to $X^\tau_b$ of
	$d\,u$ if $k=1$, and the restriction of $u\cdot\eta_k$ if $2\leq k\leq4$. In the first case, $u$ is a monomial of
	degree~$4$ in $h_X,b_2,g_2,b_3$, hence it is a product of two monomials of positive codimension. In the
	second case, $u$ has codimension $5-k>0$ and $\eta_k$ has codimension $k-1>0$. Either way we obtain that $u\cdot \pi_\ast(\xi^k)$ is
	the restriction of an element of $L$. 
    
    Next, let us assume $k=5$. Then one has $u=1$ and using \eqref{eq:pbf} we find
	\[
	\pi_*(\xi^5)=\pi_*\bigl(\pi^*(e)\,\xi^3+\pi^*(f)\,\xi^4\bigr)=e\cdot\eta_3+f\cdot\eta_4
	\quad\text{ in }\CH^4(X^\tau_b),
	\]
	where $f\in\GDCH^1_B(X_b)=\QQ\,h_X$. It follows that again $u\cdot \pi_\ast(\xi^k)$ is the restriction of a sum of elements
	of $L$. 
    
    We conclude that $ a\vert_{X^\tau_b}$ lies in the span of the restrictions to $X^\tau_b$ of
	the elements of $L$. This implies the claim because of the localization exact sequence for the open embedding $X^\tau_b\subset X_b$.
   	
	\smallskip
	\emph{Step 2: spread over $B^\circ$.}
	Let $\Xi_\rho,\Xi_\sigma\in\CH^4(\XX_B)$ be the classes of the closures in $\XX_B$ of the images of
	the two sections $x_\rho,x_\sigma\colon B^{\circ \circ}\to\XX_B\vert_{B^{\circ \circ}}$. For $b\in B^{\circ \circ}$ these
	closures meet the fibre $X_b$ transversally in the point $x_\rho(b)$ respectively
	$x_\sigma(b)$, so that $\Xi_\rho\vert_{X_b}=[x_\rho(b)]$ and $\Xi_\sigma\vert_{X_b}=[x_\sigma(b)]$.
	
	Write $L=\{U_j\cdot V_j\}_j$ and let $\Gamma\in\CH^4(\XX_B)$. By Step 1 applied to
	$a=\Gamma\vert_{X_b}$, the class $\Gamma\vert_{X_b}$ lies in the $\QQ$-span of the restrictions to $X_b$ of the finitely
	many classes $U_j\cdot V_j$, $\Xi_\rho$, $\Xi_\sigma$, for every $b\in B^{\circ\circ}$. Since $B^{\circ\circ}$ is a dense open in $B^{\circ}$, 
	Lemma~\ref{lem:spread}(ii) provides rational numbers $\lambda_j,\mu_\rho,\mu_\sigma$ such that
	\begin{equation}\label{eq:zc-spread}
	\Gamma\vert_{X_b}=\sum_j\lambda_j\,(U_j\cdot V_j)\vert_{X_b}
	+\mu_\rho\,\Xi_\rho\vert_{X_b}+\mu_\sigma\,\Xi_\sigma\vert_{X_b}
	\qquad\text{ in }\CH^4(X_b)
	\end{equation}
	for every $b\in B^\circ$.
	
	\smallskip
	\emph{Step 3: reduction to special $Y_b$.} The general double EPW sextic is associated to a special GM fourfold (this follows from Theorem \ref{thm:DKmap} below). In view of Lemma \ref{lem:spread}, to prove the theorem we may therefore reduce to the case  $X=X_b$
	where $b\in B^\circ$ is such that $Y_b$ is special. The closures of the two
	sections are contained in $P$, so that $\Xi_\rho\vert_X$ and $\Xi_\sigma\vert_X$ are represented by
	zero-cycles supported on $P\cap X=f^{-1}(p_b)$. As $Y_b$ is special, $p_b\in Y^2_{A^\perp}$ lies in
	$\mathrm{Sing}(Z)$, so that $f^{-1}(p_b)$ is a single point $s$ on the fixed surface $S$
	of $\iota$; hence $\Xi_\rho\vert_X$ and $\Xi_\sigma\vert_X$ lie in $\QQ\,h_X^4$ by
	Lemma~\ref{lem:conics-S}.
	
	Each of the classes $U_j\vert_X$ and $V_j\vert_X$ lies in $\GDCH^*_B(X)$ and has positive
	codimension. By Lemma~\ref{lem:conics-eigen}, $\GDCH^i_B(X)\subseteq\CH^i(X)^+=f^*\CH^i(Z)$; applying Lemma~\ref{lem:conics-Z} gives
	\[
	(U_j\cdot V_j)\vert_X\in f^*\bigl(\CH^i(Z)\cdot\CH^{4-i}(Z)\bigr)=f^*(\QQ\,h_Z^4)=\QQ\,h_X^4,
	\]
	where $i$ is the codimension of $U_j$, so that $0<i<4$. By \eqref{eq:zc-spread} we conclude that
	$\Gamma\vert_X\in\QQ\,h_X^4$. This proves the inclusion
	$\GDCH^4_B(X)\subseteq\QQ\,h_X^4$, and equality holds because $h_X^4$ is generically defined and
	non-zero.
\end{proof}

\subsection{A motivic splitting}

\begin{lem}\label{lem:conics-AJ}
There are relative correspondences over $B^\circ$ inducing, for every $b\in B^\circ$, mutually
inverse isomorphisms $\h^4_{\van}(Y_b)\cong\h^2_{\prim}(X_b)(-1)$ in $\Mhom$.
\end{lem}

\begin{proof}
By Corollary~\ref{cor:sextic-twists}, 
let $\alpha_\ell,\beta_\ell$ be relative correspondences over $B^\circ$  such that $\sum_\ell\beta_\ell\circ\alpha_\ell\vert_b=\Delta_{X_b}$ for
every $b\in B^\circ$. Composing the $\alpha_\ell$ with the projectors
$\tau\otimes\pi^{2j}_Y$ of Definition~\ref{defn:CK} and adding a correspondence inducing intersection with a
power of $h$, we obtain a finite number of relative correspondences $\theta_m$ whose fiberwise restrictions give maps
$$
(\theta_m\vert_b)_\ast\colon\HH^2(X_b,\QQ)\longrightarrow\HH^4_{\van}(Y_b,\QQ)(1).
$$
Consider $V\coloneqq\HH^4_{\van}(Y_b,\QQ)$. By the K\"unneth formula and Definition~\ref{defn:CK}, the Hodge
structure $\HH^*(Y_b^2,\QQ)$ is a direct sum of Tate structures, Tate twists of $V$, and Tate
twists of $V\otimes V$.
The $\theta_m\vert_b$ are the components of $\sum_\ell\alpha_{\ell}\vert_b$ along the
summands $V(k)$. 

Let us now assume that $b$ is very general, so that the Hodge structures $\HH^2_{\prim}(X_b,\QQ)$ and $V$ are
irreducible of rank $22$, and the Mumford--Tate group of $V$ is the orthogonal group of the
intersection form. It follows that $V\otimes V=\Sym^2V\oplus\wedge^2V$ contains no Hodge substructure
isomorphic to a Tate twist of $V$. Since $\sum_\ell\beta_\ell\circ\alpha_\ell\vert_b=\Delta_{X_b}$, the
Hodge structure $\HH^2_{\prim}(X_b,\QQ)$ is a direct summand of $\bigoplus_\ell\HH^*(Y_b^2,\QQ)(a_\ell)$,
and it is neither of Tate type nor a summand of a sum of twists of $V\otimes V$.
Thus, at least
one $\theta_m\vert_b$ is non-zero on $\HH^2_{\prim}(X_b,\QQ)$, and by irreducibility it must be an isomorphism.
Denote such a $\theta_m$ by~$\theta$. 
We then have $\theta^t\circ\theta\vert_b=\lambda\id$ on $\HH^2_{\prim}(X_b,\QQ)$ for some
$\lambda\in\QQ$
The scalar $\lambda$ is non-zero because $\langle\theta x,\theta y\rangle=\lambda\langle
x,y\rangle$ and both polarizations are non-degenerate. We thus obtain relative correspondences $\theta$ and $\lambda^{-1}\theta^t$ that induce
mutually inverse isomorphisms between $\h^2_{\prim}(X_b)(-1)$ and $\h^4_{\van}(Y_b)$ in $\Mhom$.

To extend to arbitrary $b\in B^\circ$, we note that
cohomology classes of relative
cycles are flat sections of the local systems $R^*\rho_*\QQ$ of the corresponding smooth projective families, so an identity between them
which holds at the very general point of $B^\circ$ holds at every point.
\end{proof}

\begin{cor}\label{cor:sym2-splitting}
There is a relative correspondence over $B^\circ$ inducing, for every $b\in B^\circ$ such that $Y_b$ is special,
a split injection in $\Mrat$
$$\Phi\colon \
\h(Y^{(2)})\hookrightarrow \h(X)(-2)\oplus\bigoplus_j\h(Y)(a_j)\oplus\bigoplus_k\1(b_k)
$$
for some $a_j,b_k$, such that the left inverse is also induced by a relative correspondence over~$B^\circ$.
\end{cor}

\begin{proof}
We begin by constructing the splitting on the level of $\Mhom$. By Definition~\ref{defn:CK} the canonical
Chow--K\"unneth projectors give $\h(Y)\cong\h^4_{\van}(Y)\oplus\bigoplus_k\1(b'_k)$ in $\Mhom$,
since $\h^{2i}(Y)$ is a sum of Tate motives for $i\neq2$ and $\h^4_{\taut}(Y)\cong\1(-2)^{\oplus2}$.
In view of the isomorphism of Lemma~\ref{lem:conics-AJ}, we thus obtain an isomorphism of homological motives
\begin{equation}\label{eqn:homY2}
\h(Y^{(2)})\cong\Sym^2\h^2_{\prim}(X)(-2)\oplus\bigoplus_j\h^4_{\van}(Y)(a_j)\oplus\bigoplus_k\1(b_k)\ \ \hbox{in}\ \Mhom .
\end{equation}
For any hyper-K\"ahler fourfold $X$ of $K3^{[2]}$-type, there is an isomorphism
induced by cup-product
   \[\Sym^2\HH^2(X,\QQ)\xrightarrow{\ \sim\ }\HH^4(X,\QQ)\] 
\cite{Verbitsky}, \cite{Bogomolov}. This implies there is an inclusion of homological motives
\[\Sym^2\h^2_{\prim}(X)\ \hookrightarrow\ \h^4(X)\ \hookrightarrow\\h(X)\ \ \hbox{in}\ \Mhom,\]
the inclusion and its left inverse being given by the small diagonal,
composed with the relevant projectors. Plugging this in the isomorphism \eqref{eqn:homY2} and replacing the summands $\h^4_{\van}(Y)(a_j)$ by
$\h(Y)(a_j)$, we obtain a split injection
\[
\Phi\colon\h(Y^{(2)})\ \hookrightarrow\ \h(X)(-2)\oplus\bigoplus_j\h(Y)(a_j)\oplus\bigoplus_k\1(b_k)
\qquad\text{ in }\Mhom,
\]
with left inverse $\Psi$. All correspondences used in this construction are restrictions of
relative correspondences over $B^\circ$: the projectors of Definition~\ref{defn:CK} and the
symmetrizing projector are polynomials in $h_i$, $c_{2,i}$ and $\Delta_Y$, the isomorphisms of
Lemma~\ref{lem:conics-AJ} are relative by construction, and the relative small diagonal of
$\XX_B\to B^\circ$, the class $h_X$ and the Chern classes of the relative tangent bundle are obviously
relative cycles. As in the proof of Lemma~\ref{lem:conics-AJ}, the identity
$\Psi\circ\Phi=\Delta_{Y^{(2)}}$ in $\Mhom$ holds at the very general point of $B^\circ$, hence at
every point of $B^\circ$.

It remains to lift $\Phi$ to $\Mrat$. Since $\h(Y^{(2)})$ is a direct summand of $\h(Y^2)$, we may
regard 
\[
N\coloneqq\Psi\circ\Phi-\Delta_{Y^{(2)}}\ \in \CH^8(Y^{(2)}\times Y^{(2)})
\]
as a self-correspondence of $Y^2$; by what we have said above this $N\in\GDCH^8_B(Y^4)$ is homologically trivial. But then Corollary~\ref{cor:Y4-alg} implies $N$ is
algebraically trivial, hence smash-nilpotent \cite{Voevodsky}, \cite{VoisinNilp}, and in particular
nilpotent for the composition of correspondences: $N^{\circ M}=0$ for some $M\in\NN$. Consequently
\[
(\Delta_{Y^{(2)}}+N)^{-1}=\Delta_{Y^{(2)}}-N+N^{\circ2}-\dots+(-1)^{M-1}N^{\circ(M-1)}
\]
is a well-defined self-correspondence of $Y^{(2)}$ modulo rational equivalence, and
$(\Delta_{Y^{(2)}}+N)^{-1}\circ\Psi$ is a left inverse of $\Phi$ in $\Mrat$. Being a polynomial in
generically defined correspondences, it is generically defined.
\end{proof}

\subsection{The MCK decomposition for special GM fourfolds}

Let $\bar B^{\mathrm{sp}}\subseteq\bar B$ be the locus that parametrized those $(H_1,H_2,Q)$
with $O\in H_1\cap H_2$, and let $B^{\mathrm{sp}}\coloneqq B\cap\bar B^{\mathrm{sp}}$, a dense open subset of
$\bar B^{\mathrm{sp}}$, parametrizing special GM fourfolds.

\begin{thm}\label{thm:MCK}
The canonical Chow--K\"unneth decomposition of a special GM fourfold is multiplicative.
\end{thm}

\begin{proof}
Let $Y=Y_b$ be a special GM fourfold. By Proposition~\ref{prop:GM-criteria}(i), it suffices to
prove that
$$
\Psi_Y=0\ \text{ in }\CH^3(Y),~ \Upsilon_Y=0\ \text{ in }\CH^6(Y^2),\text{ and }
\tau_{12}\tau_{13}=\tau_{23}o_1\ \text{ in }\CH^8(Y^3).
$$
The first is proven in Lemma~\ref{lem:Psi-special}.

We prove the second one. The class $\Upsilon_Y=\tau\cdot c_{2,1}$ is homologically trivial
by Proposition~\ref{prop:CK-basic-u}(iii) and generically defined with respect to $B$, since $\tau$ and $c_2$ are.
It is invariant under the exchange of the two factors of $Y^2$: indeed $\tau$ is symmetric, and,
since $\Delta_Y\cdot c_{2,1}=\Delta_Y\cdot c_{2,2}$, the class
\[
\tau\cdot(c_{2,1}-c_{2,2})=-\bigl(\pi^0_Y+\pi^2_Y+\pi^4_{\taut}+\pi^6_Y+\pi^8_Y\bigr)\cdot(c_{2,1}-c_{2,2})
\]
is a homologically trivial element of the decomposable subalgebra $\mathrm{D}^*(Y^2)$, hence lies in
the ideal generated by $p_1^*\Psi_Y$ and $p_2^*\Psi_Y$ by Lemma~\ref{lem:decomposable-u}, as $K=\QQ\Psi_Y$ (Section~\ref{SS:GM-u}), hence vanishes.
Consequently $\Upsilon_Y$ lies in the image of the symmetrizing projector, that is, it defines a
homologically trivial and generically defined class $\overline\Upsilon_Y$ in $\CH^6$ of the motive~$\h(Y^{(2)})$.

Let us assume that $b\in B^\circ$, and apply Corollary~\ref{cor:sym2-splitting}. Since
$\Phi$ and its left inverse are generically defined with respect to $B$, the components of
$\Phi_*\overline\Upsilon_Y$ are homologically trivial classes lying respectively in
$\GDCH^4_B(X_b)$, in the groups $\GDCH^{6+a_j}_B(Y)$ and in $\QQ$. The first group is $\QQ\,h_X^4$
by Theorem~\ref{thm:zero-cycles-X} and injects into cohomology, the second ones inject into
cohomology by Lemma~\ref{lem:Psi-special}, and the third ones inject trivially. So every component
of $\Phi_*\overline\Upsilon_Y$ vanishes and $\overline\Upsilon_Y=0$, because $\Phi$ is split
injective. Thus $\Upsilon_Y=0$ whenever $b\in B^\circ$ parametrizes a special GM fourfold. 
Since $B^\circ\cap B^{\mathrm{sp}}\subseteq B^{\mathrm{sp}}$ is a dense open subset, Lemma~\ref{lem:spread}(i) implies that $\Upsilon_Y$ vanishes for every $b\in B^{\mathrm{sp}}$.

The third vanishing follows from the second by Proposition~\ref{prop:double-cover-u}(iv), applied as in
the proof of Proposition~\ref{prop:MCK-alg}, since $\Upsilon_h=0$ by
Proposition~\ref{prop:GM-instances}(ii).
\end{proof}

\subsection{Injectivity of cycle class maps over the special locus}

\begin{prop}\label{prop:Y3}
Let $Y$ be a special GM fourfold. The cycle class map induces injections
\[ \begin{split}
& \R^*(Y^m)\hookrightarrow\HH^{2*}(Y^m,\QQ)\ \ \hbox{for}\ m\le 45,\\
&\GDCH^*_B(Y^m)\hookrightarrow\HH^{2*}(Y^m,\QQ) \ \ \hbox{for}\ m\le 3,\\
 & \GDCH_B^i(Y^m) \ \hookrightarrow\ \HH^{2i}(Y^m,\QQ)\ \ \hbox{for}\ i\le 8\  \hbox{and}\ m\leq 6.\\
 \end{split}\]
\end{prop}

\begin{proof}
By Theorem~\ref{thm:MCK} and Proposition~\ref{prop:GM-criteria}(i), the classes $\Psi_Y$, $\Upsilon_Y$
and $\Theta_Y$ vanish in the Chow ring, so that the cycle class map is injective on $\R^*(Y^m)$, for $m\leq 45$, by
Proposition~\ref{prop:GM-criteria}(iii).
We can conclude since $\GDCH^*_B(Y^m)=\R^*(Y^m)$ for $m\leq3$, and $\GDCH^i_B(Y^m)=\R^i(Y^m)$ for $i\le 8$ and $m\leq 6$, by Proposition~\ref{prop:generation}.
\end{proof}

\section{Proof of the main theorem}\label{S:main}

\subsection{The moduli stack}
We consider the quotient stack of Lagrangian spaces avoiding decomposable vectors:
  \[ \MM_{\mathrm{EPW}}:=[\mathrm{LGr}_{\mathrm{adv}}(\bigwedge^3V_6)/\mathrm{PGL}(V_6)],\] 
  where, following the notation in \cite{DK3}, $\mathrm{LGr}_{\mathrm{adv}}(\bigwedge^3V_6)$ is the open subset of $\mathrm{LGr}(\bigwedge^3V_6)$ parameterizing Lagrangian subspaces avoiding decomposable vectors.   We write
\[
  \XX\to\MM_{\mathrm{EPW}},\qquad \CC\to\MM_{\mathrm{EPW}}
\]
for the universal families of double EPW sextics and of double EPW cubes; they exist as families
over that stack because these two constructions are $\mathrm{PGL}(V_6)$-equivariant.  Chow groups of these stacks and
of the universal families over them are the $\mathrm{PGL}(V_6)$-equivariant Chow groups, with
rational coefficients, of the corresponding families over $\mathrm{LGr}_{\mathrm{adv}}(\wedge^3V_6)$
\cite{EdidinGraham}.

\subsection{The Debarre--Kuznetsov map}
Debarre--Kuznetsov attach to every GM fourfold $Y$ a Lagrangian
datum $(V_6,V_5,A)$, where $V_5\subseteq V_6$ is a hyperplane, and this datum determines $Y$
\cite[Theorem 3.6]{DK1}; the fourfold is ordinary if $\dim(A\cap\bigwedge^3V_5)=1$ and special if
$\dim(A\cap\bigwedge^3V_5)=2$.  The double EPW sextic and the double EPW cube associated with $Y$ in
the sense of Corollary~\ref{cor:sextic-twists} and Theorem~\ref{thm:relative-transfer} are the ones associated to the Lagrangian subspace
$A^\perp\subseteq\bigwedge^3V_6^\vee$ (or equivalently, for the double EPW cube, to $A$ itself), see \cite{IM,GLZ} and \cite{KKM,FGLZ2}. 

\begin{thm}[Debarre--Kuznetsov \cite{DK3}]\label{thm:DKmap} There is a  morphism of stacks
\[
  \nu\colon\MM_{\mathrm{GM}}\longrightarrow\MM_{\mathrm{EPW}},\qquad [Y]\longmapsto [A^\perp],
\]
that is dominant. On $S$-points, $\nu$ is the composite 
  \[ (\YY\to S) \ \mapsto\ (S,\VV_6,\VV_5,\mathcal A)\ \mapsto\ (S,\VV_6,\mathcal A) ,\]
  where the middle term is the family of Lagrangian data of \cite[Proposition 4.1]{DK3}. 

If $[A]\in\MM_{\mathrm{EPW}}$ parameterizes a smooth double EPW sextic $X_A$, the
fibre $\nu^{-1}([A])_{\mathrm{coarse}}$ of the map of coarse moduli spaces is isomorphic to $Z_{A^\perp}//{\mathrm{PGL}}(V_6)_A$, where ${\mathrm{PGL}}(V_6)_A$ is the stabilizer of $A$ in ${\mathrm{PGL}}(V_6)$ (the group ${\mathrm{PGL}}(V_6)_A$ is trivial for $A$ very general).
\end{thm}

\begin{proof} This is contained in \cite{DK3}. Note that in \cite[Section 6.1]{DK3} the morphism $\nu$ is constructed on the level of coarse moduli spaces. However,
closely reading \cite{DK3} one can lift this morphism to the level of moduli stacks, as we now show.
According to \cite[Proposition~4.1]{DK3} there is a morphism of stacks
  \[ a\colon \MM^{\mathrm{GM}}\ \to\ \MM^{\mathrm{Lag}}_9,\]
  defined as
  $ (S,W,\VV_6,\VV_5,u,q)\mapsto (S,\VV_6,\VV_5,\mathcal A)$,
  where ${\mathcal A}:=\ker(f_2)/\ima(f_1)$ is as in \cite[Equation (29)]{DK3}.
Furthermore, we know from \cite[Proposition 5.7]{DK3} that there is an isomorphism 
  \[ b\colon\MM^{\mathrm{Lag}}_9\cong [\bar{S}_4/\mathrm{PGL}(V_6)],\] 
  where
  \[ \bar{S}_4:=\Bigl\{ (A,V_5)\in \mathrm{LGr}_{\mathrm{adv}}(\wedge^3 V_6)\times\PP(V_6^\vee)\,\Bigl\vert\, \dim(A\cap \wedge^3 V_5)\in\{1,2\}\Bigr\} . \]
The projection $ \bar{S}_4\to  \mathrm{LGr}_{\mathrm {adv}}(\wedge^3 V_6)$ sending $(A,V_5)$ to $A$ is $\mathrm{PGL}(V_6)$-equivariant by construction, hence it induces a morphism of quotient stacks
  \[ p\colon [\bar{S}_4/ \mathrm{PGL}(V_6)]\ \to\ [\mathrm{LGr}_{\mathrm {adv}}(\wedge^3V_6)/\mathrm{PGL}(V_6)] =: \MM_{\mathrm{EPW}}.  \]
 The morphism $\nu$ is defined as the composite
$ \nu:=p\circ b\circ a $, 
  and it is readily seen this is a lift of the morphism of coarse moduli spaces constructed in \cite[Section 6.1]{DK3}.

The description of the fiber is \cite[Corollary 6.3]{DK3}, and the parenthesis is \cite[Proposition B.9]{DK1}.  
\end{proof}

\begin{rmk}
Thanks to Theorem \ref{thm:DKmap}, we can pull back generically defined cycles to our parameter space $B$.
Precisely, let $B^\circ\subseteq B$ be the dense open subset over which the double EPW sextic and double EPW cube are smooth,
and let $\XX_B\to B^\circ$ and $\CC_B\to B^\circ$ be the universal families, whose fibres $X_b$ and $C_b$ are the double EPW sextic and the
double EPW cube associated with $Y_b$.  The relative Lagrangian datum of $\YY\to B$ is a family of
Lagrangian subspaces of $\bigwedge^3\VV_6$ for a rank six vector bundle $\VV_6$ on $B$ \cite{DK3}, so that
 \[ \nu_B\colon  B^\circ \to \MM_{\mathrm{GM}}\to \MM_{\mathrm{EPW}}\]
 is a morphism to the stack $\MM_{\mathrm{EPW}}$,
and $\XX_B$ and $\CC_B$ are the pullbacks of $\XX$ and $\CC$ along it. Then pulling back
cycles along $\nu_B$ gives inclusions for any $m \geq 1$:
\begin{equation}\label{eq:basechange}
   \begin{split}
  \GDCH^*_{\MM_{\mathrm{EPW}}}(X^m)\ &\subseteq\ \GDCH^*_{B^\circ}(X^m),\\
  \GDCH^*_{\MM_{\mathrm{EPW}}}(C^m)\ &\subseteq\ \GDCH^*_{B^\circ}(C^m),\\
  \end{split}
\end{equation}
where we use the identifications $X=X_b=X_{\nu_B(b)}$ and
$C=C_b=C_{\nu_B(b)}$.
\end{rmk}

\subsection{A crucial reduction step}
Let $\MM_{\mathrm{IKKR}}$ denote the open substack of the moduli stack of
polarized hyper-K\"ahler sixfolds of $K3^{[3]}$-type parametrizing double EPW cubes, and let
$\CC_{\mathrm{IKKR}}\to\MM_{\mathrm{IKKR}}$ be its universal family.  The pullback of this family to $\MM_{\mathrm{EPW}}$ under the morphism
$\MM_{\mathrm{EPW}}\to\MM_{\mathrm{IKKR}}$, $A\mapsto (C_A,h_{C_A})$ is the universal family $\CC\to \MM_{\mathrm{EPW}}$.

The following lemma will be crucial to our argument proving Theorem \ref{thm:main}(iii):

\begin{lem}\label{lm:cube-stack} Let $C$ be a double EPW cube, and $\iota_C$ its covering involution. There is an inclusion
  \begin{equation}\label{eq:cube-stack} \GDCH^*_{\MM_{\mathrm{IKKR}}}(C)\ \subseteq\ \GDCH^*_{\MM_{\mathrm{EPW}}}(C)^{\iota_C}.\end{equation}
  \end{lem}

  \begin{proof} Let $[(C,h_C)]\in\MM_{\mathrm{IKKR}}$ be a point.
The automorphism group of $[(C,h_C)]$ is the group
$G:=\mathrm{Aut}(C,h_C)$ of automorphisms of $C$ preserving $h_C$. The group $G$ is finite and contains the covering involution $\iota_C$. In fact, if $[(C, h_C)]$ is very general in $\MM_{\mathrm{IKKR}}$, $G$ is of order~2 and generated by $\iota_C$.
The hyper-K\"ahler variety $C$ carries no non-zero vector field, and so the stack $\MM_{\mathrm{IKKR}}$ is
a Deligne--Mumford stack. 

Locally around the point $[(C,h_C)]$, there exists an affine scheme $V$ with a $G$-action, a fixed point $v$, and an \'etale representable morphism $[V/G]\to \MM_{\mathrm{IKKR}}$ sending $v$ to $[(C,h_C)]$, such that the universal family pulls back to $[\CC_V/G]$.
Since deformations of hyper-K\"ahler varieties are unobstructed, we may assume that $V$ is smooth of dimension $20$, so that $v\hookrightarrow V$ is a regular embedding and we can apply the refined Gysin pullback along it.
Furthermore, the inclusion $C\to\CC_V$ is $G$-equivariant and the composite $C\to\CC_V\to[\CC_V/G]$ is $G$-invariant.

The Chow groups with rational coefficients of a quotient stack are the equivariant Chow groups of
\cite{EdidinGraham}. Therefore, the restriction to the fibre $C$ of a cycle on $\CC_{\mathrm{IKKR}}$ factors through $[\CC_V/G]$, hence is invariant under $\mathrm{Aut}(C,h_C)$ and in particular under $\iota_C$. 
\end{proof}

\subsection{Proof of Theorem \ref{thm:main}}
We are now ready to prove Theorem~\ref{thm:main}, which asserts that the universal families
$\XX\to\MM_{\mathrm{EPW}}$ and $\CC_{\mathrm{}}\to\MM_{\mathrm{IKKR}}$ have the Franchetta
property. 

\begin{proof}[Proof of Theorem~\ref{thm:main}]
By Definition~\ref{defn:gdch} and Lemma~\ref{lem:spread}(i) it suffices to prove that the cycle
class map is injective on
$\GDCH^*_{\MM_{\mathrm{EPW}}}(X_A)$ and on $\GDCH^*_{\MM_{\mathrm{IKKR}}}(C_A)$ for 
a very general
$[A]\in\MM_{\mathrm{EPW}}$.  Fix such an $[A]$ and write $X\coloneqq X_A$ and $C\coloneqq C_A$.  Since $B\to
\MM_{\mathrm{GM}}$ is dominant and the fibre $\nu^{-1}([A])$ contains special GM fourfolds (Theorem \ref{thm:DKmap}), we may choose a point $b\in B^\circ$ with $\nu_B(b)=[A]$ such that the GM fourfold $Y\coloneqq Y_b$ is special.  Under
the identifications $X=X_b$ and $C=C_b$, combining the inclusions \eqref{eq:basechange} and
\eqref{eq:cube-stack} gives
\[
   \begin{split}
  \GDCH^*_{\MM_{\mathrm{EPW}}}(X_{})\ &\subseteq\ \GDCH^*_{B}(X),\\
   \GDCH^*_{\MM_{\mathrm{EPW}}}(X\times X)\ &\subseteq\ \GDCH^*_{B}(X\times X),\\
    \GDCH^*_{\MM_{\mathrm{EPW}}}(X^3)\ &\subseteq\ \GDCH^*_{B}(X^3),\\
  \GDCH^*_{\MM_{\mathrm{IKKR}}}(C_{})\ &\subseteq\ \GDCH^*_{B}(C)^{\iota_C},\\
  \end{split}
\]

Thanks to Corollary \ref{cor:sextic-twists} and Theorem \ref{thm:relative-transfer}, we can relate the right-hand side to generically defined cycles on powers of $Y$ and obtain inclusions
  \begin{equation}\label{eq:new-inclusion}
   \begin{split}
  \GDCH^*_{\MM_{\mathrm{EPW}}}(X_{})\ &\subseteq\ \bigoplus \GDCH^*_{B}(Y^2),\\
   \GDCH^i_{\MM_{\mathrm{EPW}}}(X\times X)\ &\subseteq\ \bigoplus \GDCH^{i+4}_{B}(Y^4),\\
   \GDCH^i_{\MM_{\mathrm{EPW}}}(X^3)\ &\subseteq\ \bigoplus \GDCH^{i+6}_{B}(Y^6),\\
  \GDCH^*_{\MM_{\mathrm{IKKR}}}(C_{})\ &\subseteq\ \bigoplus \GDCH^*_{B}(Y^3),\\
  \end{split}
\end{equation}

Since $Y$ is special, the cycle class map is injective on $\GDCH^*_B(Y^2)$ and
$\GDCH^*_B(Y^3)$ in all codimensions,  and on $\GDCH^j_B(Y^4)$ and $\GDCH^j_B(Y^6)$ for $j\le 8$, by Proposition~\ref{prop:Y3}. 
This yields the desired injectivity of the cycle class maps on
the left-hand side of \eqref{eq:new-inclusion} in the range stated in Theorem \ref{thm:main}.
\end{proof}

\begin{rmk}\label{rmk:Plucker}
One can show that the inclusions \eqref{eq:basechange} are strict in general, and that neither the family $\XX_B\to
B^\circ$ nor the relative square $\YY^2_{/B}\to B$ has the Franchetta property.

Indeed, let $b\in B^\circ$ and let $x_\rho(b),x_\sigma(b)\in X_b$ be the
two points to which the Iliev--Manivel morphism $\pi\colon F(Y_b)\to X_b$ contracts the loci
$F_\rho$ and $F_\sigma$.  They are the two points of $X_b$ lying above the Pl\"ucker point $p_b\in
Z_b$ of $Y_b$ \cite[Theorem 7.12]{DK5}.
As $b$ varies, the Pl\"ucker point defines a section of the family of EPW sextics
$\XX_B/\iota\to B^\circ$, and therefore $w(b):=x_\rho(b)+x_\sigma(b) \in\GDCH^4_{B^\circ}(X_b)$ for every
$b\in B^\circ$.  Assume that the cycle class
map were injective on $\GDCH^4_{B^\circ}(X_b)$ for very general $b$.  The class $6\,w(b)-h_X^4$ is
generically defined and homologically trivial, so it would vanish for very general $b$, and hence
for every $b\in B^\circ$ by Lemma~\ref{lem:spread}(i).  Fix a very general $[A]\in\MM_{\mathrm{EPW}}$
and put $X\coloneqq X_A$.  Letting $b$ run through the points of $\nu_B^{-1}([A])$ that parameterize ordinary
GM fourfolds, the corresponding Pl\"ucker points fill up $Z_A\smallsetminus\mathrm{Sing}\,Z_A$ (Theorem \ref{thm:DKmap}), and we would get $f^*(z)=\frac16h_X^4$ for every~$z$ in a dense subset of
$Z_A$, hence for every $z\in Z_A$ by Lemma \ref{lem:spread}(i). This would imply that the
correspondence 
\[\Lambda:=\Delta_X+\Gamma_\iota-\frac16(X\times h_X^4)\in\CH^4(X\times X)\] 
annihilates the Chow group of zero-cycles of $X$.  By the Bloch--Srinivas argument \cite{BlochSrinivas} (see also \cite[Chapter 3]{VoisinBook}), the
 correspondence $\Lambda$ would then be supported on $D\times X$ for some divisor $D\subseteq X$, and so act as zero on
$\HH^{4,0}(X)$.  But this is absurd: the involution $\iota$ is anti-symplectic, so $\iota^*$ is the
identity on $\HH^{4,0}(X)=\C\,\sigma^2$, where $\sigma$ is a symplectic form on~$X$, and $\Lambda$ acts there as multiplication by $2$.

Thus the cycle class map is not injective on $\GDCH^4_{B^\circ}(X_b)$ for very general $b$. 
Corollary~\ref{cor:sextic-twists} shows that the cycle class map is not
injective on $\GDCH^*_{B^\circ}(Y_b^2)$ either.  
\end{rmk}

\begin{rmk}\label{rmk:noMCK} Remark \ref{rmk:Plucker} implies in particular that
$\Psi_{Y_b}$ or $\Upsilon_{Y_b}$ is non-zero for very general~$b$, by
Proposition~\ref{prop:GM-criteria}(iii). Thus, it follows from Proposition~\ref{prop:GM-criteria}(i) that the
canonical Chow--K\"unneth decomposition of a very general GM fourfold is not multiplicative.

This has the following intriguing consequence: let $Y$ be a GM fourfold for which the canonical Chow--K\"unneth decomposition fails to be multiplicative, and let $Y^\sharp$ be a period partner of $Y$ that is special. Then $Y^\sharp$ has a multiplicative Chow--K\"unneth decomposition (Theorem~\ref{thm:MCK}), $Y$ and $Y^\sharp$ have equivalent Kuznetsov components \cite[Theorem 1.6]{KP2} and isomorphic Chow motives \cite[Theorem 8.2]{FM}. Thus, there is a graded isomorphism $\h^\ast(Y)\cong \h^\ast(Y^\sharp)$ where both sides are given the canonical Chow--K\"unneth decomposition.
However, this isomorphism is {\em not\/} an isomorphism of Frobenius algebra objects. This is in stark contrast to the case of cubic fourfolds: in \cite[Theorem 1]{FV2}, it is proven that cubic fourfolds with equivalent Kuznetsov component have isomorphic Chow motives as Frobenius algebra objects, and this isomorphism is compatible with the canonical Chow--K\"unneth decomposition.
\end{rmk}


\section{Applications to motives of double EPW sextics}\label{S:consequences}

In this section, we establish some motivic properties of double EPW sextics, which we deduce from Theorem~\ref{thm:main}(ii) and Proposition~\ref{prop:Y3}.

\subsection{Dual double EPW sextics}\label{SS:cons-dual}

Let $\mathrm{LGr}(\bigwedge^3 V_6)_{00}$ be the open subset of Lagrangian subspaces $A$ such that $A$ and $A^\perp$ contain no decomposable vectors and on which the associated double EPW sextics are smooth. 
O'Grady showed in \cite{OG-dual} that the rational cohomology groups $\HH^2_{\prim}(X_A,\QQ)$ and $\HH^2_{\prim}(X_{A^\perp},\QQ)$ are Hodge isometric, but that their period points differ, i.e. that this does not lift to an isomorphism of integral Hodge structures.
This means that although $X_A$ and $X_{A^\perp}$ are not birational in general, we would still expect their motives to be isomorphic. 
Unfortunately, the rational Hodge isometry constructed by O'Grady is a priori not induced by an algebraic correspondence, so we cannot lift it directly to motives.

O'Grady showed that the EPW sextics $Z_A$ and $Z_{A^\perp}$ are projectively dual \cite{OGrady06}, and hence we denote them by $Z$ and $Z^\vee$. 
The Gauss map $Z\dashrightarrow Z^\vee$ is birational and undefined on the surface along which $Z$ is singular. 
The closure of the graph of the Gauss map is both the blow-up of $Z$ in its singular locus, as well as the blow-up of $Z^\vee$ in its singular locus \cite{OG-dual}.

\begin{prop}\label{prop:inv-motives}
    Let $A\in\mathrm{LGr}(\bigwedge^3 V_6)_{00}$ and let $X_A$ and $X_{A^\perp}$ be the associated dual double EPW sextics. 
    The Gauss map on $Z$ and $Z^\vee$ induces an isomorphism
    $$ \h(X_A)^+\xrightarrow{\sim} \h(X_{A^\perp})^+ \text{ in }\Mrat,$$
    where $\h(X_A)^+$ is the invariant motive, cut out by the projector $\frac{1}{2}(\Delta_{X_A}+\Gamma_{\iota_{X_A}})$.
\end{prop}
\begin{proof}
    Denote $X\coloneqq X_A$ and $X^\vee\coloneqq X_{A^\perp}$.
    We prove that the Gauss map induces an isomorphism on cohomology and then use the Franchetta property to lift the isomorphism to motives. 
    In degrees $0$ and $8$, it is clearly an isomorphism.
    Now, use the blow-up formula to see that $H$ and $H^\vee$, the pullbacks of the polarizations of $Z$ and $Z^\vee$ to the blow-up, satisfy $H^\vee=5H-E$ on the blow-up, and we can compute $\deg(H\cdot H^{\vee 3})=30$, meaning that it is also an isomorphism in degrees $2$ and $6$.

    In degree $4$, the blow-up formula shows that 
    $$\HH^4(Z,\QQ)\oplus \HH^2(\mathrm{Sing}(Z),\QQ)(-1)\cong\HH^4(Z^\vee,\QQ)\oplus \HH^2(\mathrm{Sing}(Z^\vee),\QQ)(-1) $$
    and it suffices to prove that this isomorphism maps $\HH^4(Z,\QQ)$ isomorphically to $\HH^4(Z^\vee,\QQ)$.
    We have $H^4(Z,\QQ)=\QQ h^2_Z\oplus\QQ q_A\oplus W$, where $q_A$ is the class of the Beauville--Bogomolov form of $X$ and $W$ is an irreducible Hodge structure containing $\HH^{4,0}(Z)$.
    Since the Gauss map is birational, it is an isomorphism on $\HH^{4,0}(Z)$ and hence on $W$.
    The pullback of $\QQ h_Z^2\oplus \QQ q_A$ to the blow-up is spanned by $H^2$ and $4HH^\vee-H^{\vee 2}$ and symmetrically by $H^{\vee 2}$ and~$4HH^\vee-H^2$.
    The matrix representing the pairing between these two bases has determinant~53504, meaning that the Gauss map induces an isomorphism $\HH^\ast(Z,\QQ)\cong \HH^\ast(Z^\vee,\QQ)$.

    Hence, the Gauss map induces generically defined correspondences $\Gamma$ and ${}^t\Gamma$ which are isomorphisms on cohomology.
    Since they are generically defined, Theorem~\ref{thm:main}(ii) shows that the cycle class map is injective on the subalgebra generated by $\Gamma\circ{}^t\Gamma$ and ${}^t\Gamma\circ\Gamma$. The Cayley--Hamilton theorem then allows one to find a correspondence inverse to $\Gamma$, completing the proof.
\end{proof}

\begin{thm}\label{cor:dual-motives}
    Let $A\in\mathrm{LGr}(\bigwedge^3 V_6)_{00}$ and let $X_A$ and $X_{A^\perp}$ be dual double EPW sextics. 
    Then there exists an isomorphism 
    $$ \h(X_A)\xrightarrow{\sim} \h(X_{A^\perp}) \text{ in }\Mrat.$$
\end{thm}

\begin{proof}
    Denote $X\coloneqq X_A$ and $X^\vee\coloneqq X_{A^\perp}$.
    We decompose $\h(X)=\h(X)^+\oplus\h(X)^-$.
    Proposition~\ref{prop:inv-motives} yields an isomorphism $\h(X)^+\cong\h(X^\vee)^+$, so it remains to construct an isomorphism $\h(X)^-\cong\h(X^\vee)^-$.
    To this end, we use $\h(X)^-\cong\h^4_{\van}(Y)(1)\oplus\h^4_{\van}(Y)\oplus\h^4_{\van}(Y)(-1)$ \cite[Proposition~4.4]{BL2} and $\h^4(Y)\cong\h^4(Y^\vee)$ \cite[Theorem~8.2]{FM} for general GM fourfolds $Y$ and~$Y^\vee$ to which $X$ and $X^\vee$ are associated.
    We point out that choosing a special GM fourfold, the first isomorphism can be made to be generically defined via Lemma~\ref{lem:rat-conics-AJ}, but it is not clear that the construction of the isomorphism of \cite{FM} is generically defined. 
\end{proof}

\subsection{A motivic relation between double EPW sextics and Gushel--Mukai fourfolds}\label{SS:cons-sym2}

We begin with the following result, lifting Lemma~\ref{lem:conics-AJ} to $\Mrat$.

\begin{lem}\label{lem:rat-conics-AJ}
    There are relative correspondences over $B^\circ$ inducing, for every $b\in B^\circ$ such that $Y_b$ is special, mutually inverse isomorphisms 
    $\h^4_{\van}(Y_b)\cong\h^2_{\prim}(X_b)(-1)$ in $\Mrat$.
\end{lem}
\begin{proof}
    Lemma~\ref{lem:conics-AJ} constructs two relative correspondences over $B^\circ$ inducing mutually inverse isomorphisms $\h^4_{\van}(Y_b)\cong\h^2_{\prim}(X_b)(-1)$ in $\Mhom$ for every $b\in B^\circ$.
    Call these isomorphisms $\Gamma$ and $\Gamma^{-1}$. 
    Then, $\Gamma\circ\Gamma^{-1}=\id_{\h^2_{\prim}(X_b)(-1)}$ and $\Gamma^{-1}\circ\Gamma=\id_{\h^4_{\van}(Y_b)}$ and using Corollary~\ref{cor:sextic-twists} and Proposition~\ref{prop:Y3} this lifts to $\Mrat$.
\end{proof}

Theorem~\ref{thm:main}(ii) implies that $X$ admits a generically defined Chow--K\"unneth decomposition. We denote by $\h^{2i}(X)$ the weight-$2i$ direct summand of $\h(X)$ given by this decomposition.
The definition of $\h^2(X)$ agrees with the one in Lemma~\ref{lem:conics-AJ} and Lemma~\ref{lem:rat-conics-AJ}. 
Let $$\mu\colon \Sym^2\h^2(X)\to\h^4(X)$$ be the intersection-product morphism (i.e.~induced by the small diagonal).

\begin{thm}\label{thm:sym2}
    Let $X$ be a double EPW sextic.
The morphism $\mu\colon\Sym^2\h^2(X)\to\h^4(X)$ is an isomorphism in $\Mrat$.
\end{thm}

\begin{proof}
Since the cup product map $\Sym^2H^2(X)\to H^4(X)$ is an isomorphism for hyper-K\"ahler manifolds of $\mathrm{K}3^{[2]}$-type, $\mu$ induces an isomorphism on cohomology.
Analogously to Lemma~\ref{lem:rat-conics-AJ}, let $Y$ be a special GM fourfold such that the associated double EPW sextic is~$X$, then we have 
\begin{equation}\label{eqn:IsomH6Prim}
    \h^6_{\prim}(X)(1)\cong\h^4_{\van}(Y).
\end{equation}
Let $\zeta\colon\h^6(X)(2)\to\h^2(X)$ be the generically defined isomorphism given by $\frac{1}{12}h_X\times h_X$ between the Tate summand $\1(-1)$ of $\h^6(X)(2)$ spanned by $h_X^3$ and the Tate summand $\1(-1)$ of $\h^2(X)$ spanned by $h_X$, and by Lemma~\ref{lem:rat-conics-AJ} and \eqref{eqn:IsomH6Prim} between the primitive parts.
Let ${}^t\mu\colon \h^4(X)\to \Sym^2(\h^6(X)(2))$ be the transpose of $\mu$ combined with motivic Poincar\'e duality.
Then $\Sym^2\zeta\circ{}^t\mu$ induces an isomorphism on cohomology and using Corollary~\ref{cor:sextic-twists} and Proposition~\ref{prop:Y3} and the Cayley--Hamilton theorem, we obtain that $\mu$ is an isomorphism. 
\end{proof}

\begin{rmk}
    The proof of Proposition~\ref{prop:alg-Y2} shows that $\h^4(Y)$ is a summand in $\h(S)(-1)$ for some surface.
    Hence $\h^2(X)^-$ is too and this can be combined with Theorem~\ref{thm:sym2} to recover \cite[Theorem~A.5]{Mazzanti}, which uses the explicit geometry of double EPW sextics to establish that the intersection product $\CH^2(X)^-_{\hom}\otimes\CH^2(X)^-_{\hom}\to\CH^4(X)^+_{\hom}$ is surjective.
\end{rmk}

\begin{cor}\label{cor:motive-X}
Let $X$ be a double EPW sextic associated to a GM fourfold $Y$.
Then there is an isomorphism
\begin{equation}
\label{eqn:MotivicBDrelation}
    \h(X)\cong\Sym^2(\h^4_{\van}(Y)(1))\oplus\bigoplus_{i=1}^{3}\h^4_{\van}(Y)(2-i)\oplus\bigoplus_{k=0}^4\1(-k)\ \ \text{ in }\Mrat.
\end{equation}
The same holds for the dual double EPW sextic $X^\vee$ associated to $Y$.
\end{cor}

\begin{proof}
    Assume first that $Y$ is a special GM fourfold.
    Lemma~\ref{lem:rat-conics-AJ} shows that $\h^2(X)\cong\1(-1)\oplus\h^4_{\van}(Y)(1)$ and one can prove analogously that $\h^6(X)\cong\1(-3)\oplus\h^4_{\van}(Y)(-1)$.
    Decomposing into invariant and anti-invariant parts in Theorem~\ref{thm:sym2} one obtains $\h^4(X)\cong\1(-2)\oplus\h^4_{\van}(Y)\oplus\Sym^2(\h^4_{\van}(Y)(1))$, yielding the formula in the statement.
    To extend to arbitrary $Y$, we use the fact that any GM fourfold has a period partner that is special, and apply \cite[Theorem~8.2]{FM}, which says that GM fourfolds which are period partners have isomorphic motives. 
    To prove the last claim about $X^\vee$, recall that $X^\vee$ is the double EPW sextic associated to the dual GM fourfold $Y^\vee$, hence \eqref{eqn:MotivicBDrelation} expresses the motive of~$X^\vee$ in terms of the motive of $Y^\vee$. Applying again \cite[Theorem~8.2]{FM}, which says that $\h(Y)\simeq \h(Y^\vee)$, we obtain that $\h(X)\simeq \h(X^\vee)$.    
\end{proof}

\begin{rmk}
    Corollary \ref{cor:motive-X} gives an alternative proof of Theorem \ref{cor:dual-motives}. 
    However, both  proofs use \cite[Theorem 8.2]{FM}, which ultimately relies on the derived equivalence between Kuznetsov components of dual GM fourfolds established in \cite{KP2}. Consequently, the algebraic cycle inducing the isomorphism $\h(X)\simeq \h(X^\vee)$ is only implicit.
\end{rmk}

\subsection{Recovering Zhang's results}

Zhang proved in \cite{Zhang} that the universal conic obtained from the construction in Theorem~\ref{thm:FX} induces an isomorphism $\CH_0(X)^-\xrightarrow{\sim}\CH_1(Y)_{\mathrm{hom}}$ and that it vanishes on $\CH_0(X)^+_{\mathrm{hom}}$.
These results are obtained via explicit geometric methods and the aim of this section is to show how these can be recovered using the results in this paper. 

It is important to note that this requires $Y$ to be an ordinary GM fourfold which does not contain a plane, and hence we cannot apply the Franchetta property directly, see also Remark~\ref{rmk:Plucker}.

Consider the correspondences $ P\coloneqq (\pi\times\id)_\ast(\tilde P\cdot p_1^\ast\xi)\in\CH^3(X\times Y)$, where $\xi$ is the Plücker class on $F(Y)$, and $I^Y\coloneqq {}^tP\circ P\in\CH^2(X\times X)$, which depends on the choice of~$Y$.
We begin by proving the injectivity of $P$ on $\CH_1(Y)_{\mathrm{hom}}$, which uses the description of tautological cycles on GM fourfolds, even though the cycle class map is not injective on these. 

\begin{prop}\label{prop:zhang-injectivity}
    Let $Y$ be an ordinary GM fourfold containing no planes.
    Then the action of ${}^t P$ on $\CH_1(Y)_{\mathrm{hom}}$ is injective, and that of $P$ is surjective.
\end{prop}

\begin{proof}
    Consider the composition $\Gamma=P\circ(\delta_\ast h_X^2)\circ{}^t P\in\CH^4(Y\times Y)$.
    It is generically defined over $B$, so it is tautological by Proposition~\ref{prop:generation}.
    Its cohomology class is $[\Gamma]=\lambda[\tau]+[D]$, where $[D]$ is the cohomology class of some decomposable cycle $D$ and $\tau$ is as in Definition~\ref{defn:CK}. 
    We have $\Gamma=\lambda\tau+D+D_{\mathrm{hom}}$ for some homologically trivial tautological class $D_{\mathrm{hom}}$.
    By Proposition~\ref{prop:GM-criteria}, the homologically trivial part of $\R^4(Y^2)$ is spanned by $\Psi_Y\times h$ and $h\times \Psi_Y$.
    Decomposable classes act as zero on homologically trivial classes, and so the action of $\Gamma$ on $\CH_1(Y)_{\mathrm{hom}}$ agrees with that of $\lambda\tau$.
    Since $\tau$ acts as the identity and $\lambda$ can be checked to be non-zero in cohomology, we obtain the statement.
\end{proof}

We now continue by proving that the class $I^Y$ is in fact intrinsic to $X$, i.e. that it does not depend on the choice of $X$. 
To do this, we use the fact that a choice of GM fourfold corresponds to the choice of a point on the associated EPW sextic, and hence can be encoded in an additional copy of $X$.

\begin{prop}\label{prop:I-independence}
    The correspondence $I^Y\in\CH^2(X\times X)$ is independent of the choice of the GM fourfold $Y$, commutes with the involution on $X$, and $\Delta^+\circ I\circ\Delta^+$ is decomposable.
\end{prop}

\begin{proof}
    Let $\MM_{\mathrm{GM}}^{\circ\circ}$ be the image of $B^{\circ\circ}$ in $\MM_{\mathrm{GM}}$, i.e. the open moduli stack parametrizing ordinary GM fourfolds that do not contain a plane. 
    Furthermore, let $f\colon\XX\to\mathcal{Z}$ be the universal double cover over $\MM_{\mathrm{EPW}}$, where $\XX$ and $\mathcal{Z}$ are the universal families of double EPW sextics and EPW sextics, respectively. 
    There is a natural map $\mathcal{Z}\to\MM_{\mathrm{GM}}$ over $\MM_{\mathrm{EPW}}$ sending a point $[V_1]\in Z_A\subset\PP(V_6)$ to the GM fourfold associated to the Lagrangian data $(V_6^\vee,V_1^\perp, A^\perp)$. 
    Furthermore, denote by $\XX^{\circ\circ}$ the preimage of $\MM^{\circ\circ}_{\mathrm{GM}}$ under $\XX\to\mathcal{Z}\to\MM_{\mathrm{GM}}$.

    As explained in Remark~\ref{rmk:IM-in-families}, the correspondence $I^Y$ can be constructed in families, so we obtain a correspondence 
    $$ \mathcal{I}^\circ\in\CH^2(\XX^{\circ\circ}\times_{\MM_{\mathrm{EPW}}}\XX\times_{\MM_{\mathrm{EPW}}}\XX), $$
    such that $\mathcal{I}^\circ|_{\{x\}\times X_A\times X_A}=I^{Y_{f_A(x)}}$, where $Y_{f_A(x)}$ is the GM fourfold associated to the point~$f_A(x)\in Z_A$.
    Notice that $\mathcal{I}^\circ$ is pulled back from $\mathcal{Z}^{\circ\circ}\times_{\MM_{\mathrm{EPW}}}\XX\times_{\MM_{\mathrm{EPW}}}\XX$, where $\mathcal{Z}^{\circ\circ}$ is the preimage of $\MM_{\mathrm{GM}}^{\circ\circ}$ under $\mathcal{Z}\to\MM_{\mathrm{GM}}$.
    Choose a closure of this cycle and denote its pullback to $\XX_{/\MM_{\mathrm{EPW}}}^{3}$ by $\mathcal{I}$.
    
    Clearly, $\mathcal{I}$ is invariant under the involution $\iota$ in the first coordinate.   
    We can use this to compute that the class 
    $$\mathcal{S}\coloneqq \mathcal{I}-p_{23}^\ast\mathcal{T}-\mathcal{D} $$
    is homologically trivial on the fibres, where $\mathcal{T}\coloneqq\frac{1}{12}p_{23\ast}(\mathcal{I}\cdot p_1^\ast h_\XX^4)$ and $\mathcal{D}$ is a class whose cohomology class is monodromy invariant and lies in $\HH^2(X)^+\otimes \HH^2(X^2)$ and $\HH^4(X)^+\otimes \HH^0(X^2)$.
    In particular, $\mathcal{D}$ lies in the span of $p_1^\ast h_\XX\cdot p_2^\ast h_\XX$, $p_1^\ast h_\XX\cdot p_3^\ast h_\XX$, $p_1^\ast h_\XX^2$, and $p_1^\ast c_2$.
    Since $\mathcal{S}$ is generically defined over $\MM_{\mathrm{EPW}}$, we may use the Franchetta property for the relative cube of double EPW sextics established in Theorem~\ref{thm:main}(ii) to obtain that $\mathcal{S}=0$.
    Restricting it to $\{x\}\times X\times X$ gives $I^Y=T\coloneqq \mathcal{T}|_{X\times X}$, showing the first statement.

    Consider now $\Delta^-\circ I\circ \Delta^+$ and $\Delta^+\circ I\circ\Delta^-$. 
    They are generically defined over $\MM_{\mathrm{EPW}}$ and their Künneth components lie in $(h_X^\perp\otimes\QQ h_X)\oplus(h_X\cdot h_X^\perp)\otimes \HH^0(X)$, and these contain no monodromy invariants, so they vanish by Theorem~\ref{thm:main}(ii).
    The vanishing of these two components proves that $I$ commutes with the involution.
    An analogous argument shows that $\Delta^+\circ I\circ\Delta^+$ is decomposable.
\end{proof}

This result answers a question by Laterveer posed in \cite[Remark~2.26]{LatBV}, which asked whether the class $I^Y$ (which is written $I_\pi$ in loc. cit.) is intrinsic to $X$. 
We can now recover the rest of Zhang's results.

\begin{thm}[Zhang]\label{thm:zhang}
    Let $X$ be a double EPW sextic associated to an ordinary GM fourfold $Y$ that does not contain a plane. 
    Then the correspondence $P$ induced by the universal conic is zero on $\CH_0(X)^+_{\mathrm{hom}}$ and induces an isomorphism $\CH_0(X)^-\xrightarrow{\sim}\CH_1(Y)_{\mathrm{hom}}$.
    The inverse is given up to some scalar by the correspondence ${}^tP$ followed by multiplication by~$h_X^2$.
\end{thm}

\begin{proof}
    The fact that the correspondence $P$ is zero on $\CH_0(X)^+_{\mathrm{hom}}$ follows from Proposition~\ref{prop:zhang-injectivity}, together with the fact that $\Delta^+\circ I\circ\Delta^+$ is decomposable and hence acts trivially on homologically trivial cycles. 

    We now prove that there exists some $\lambda\in\QQ$ such that $h_X^2\cdot I^-_\ast (\alpha)=\lambda\alpha$ for all $\alpha\in\CH_0(Y)^-$.
    Consider the generically defined Chow--Künneth projector $\pi^{6,-}$ whose cohomology class is the projector onto $H^6(X)^-$\,; its existence is provided by Theorem~\ref{thm:main}(ii).
    The class $(p_2^\ast h_X^2\cdot I)\circ\pi^{6,-}-\lambda\pi^{6,-}\in\CH^4(X\times X)$ is generically defined over $\MM_{\mathrm{EPW}}$ and is homologically trivial for an appropriate choice of $\lambda$.
    Applying Theorem~\ref{thm:main}(ii) shows that it is zero. 
    To prove that $\lambda$ is non-zero, it suffices to check that $I^-$ acts non-trivially on cohomology. 
    This follows from the fact that $P$ sends the $(3,1)$-form on $Y$ to the symplectic form on $X$, see \cite[Section~4]{IM}.

    The fact that the inverse is given by some multiple of ${}^t P$ followed by multiplication by~$h_X^2$ can be seen from the above and the proof of Proposition~\ref{prop:zhang-injectivity}.
\end{proof}

\begin{rmk}
    We point out that this result is in fact slightly stronger than Zhang's, in two respects. 
    First, the generality assumptions on $Y$ are made explicit.
    Secondly, it can only be seen from \cite{Zhang} that the inverse is given by multiplication by the class $I_\ast(\frac{1}{12}h_X^4)$.
    Since it is generically defined, it is in fact a linear combination of $h_X^2$ and $c_2(X)$ and one can use \cite[Lemma~9.3]{Ferretti} to show that multiplication by $h_X^2$ and $c_2(X)$ agree up to a scalar on~$\CH^\ast(X)^-$. 
\end{rmk}

A further consequence regarding $I$ concerns the projector $\Delta_X^-:={1\over 2}(\Delta_X-\Gamma_\iota)$ defining the motive $h(X)^-$.

\begin{cor}\label{cor:anti} Let $X$ be a double EPW sextic. There is an equality
  \[ \Delta_X^-=I\cdot A + B\ \ \hbox{in}\ \CH^4(X\times X),  \]
  where $I:=I^Y \in\GDCH_{\MM_{\mathrm{EPW}}}^2(X\times X)$ is as in Proposition~\ref{prop:I-independence} and $A,B\in \bigl\langle (p_i)^\ast(h),(p_i)^\ast(c_2)\bigr\rangle$.
\end{cor}
 
\begin{proof} The equality $\Delta_X^-=I\cdot A +B$ is proven modulo algebraic equivalence in \cite[Proof of Proposition 3.2]{LatBV}, with all terms being generically defined with respect to $\MM_{\mathrm{EPW}}$
by Proposition~\ref{prop:I-independence}.
Applying Theorem~\ref{thm:main}(ii) concludes the proof.
\end{proof}

\begin{rmk} Corollary \ref{cor:anti} improves on the more complicated expression for $\Delta_X^-$ obtained in \cite[Proposition 3.2]{LatBV}.
The expression for $\Delta_X^-$ is analogous to the expression for $\Gamma_\varphi$ for Fano varieties of lines on cubic fourfolds obtained in \cite{Voisin-ChowRingHK} and \cite{SV}, where $\Gamma_\varphi$ is the closure of the graph of Voisin's rational self-map.
\end{rmk}


\section{Analogous results for Gushel--Mukai sixfolds}\label{S:sixfolds}

We show that the methods of Sections~\ref{S:tautological} to~\ref{S:MCK} also apply to GM sixfolds: their canonical Chow--K\"unneth decomposition is multiplicative (Theorem~\ref{thm:MCK6}), and generically defined cycles on their powers $Y^m$, $m\leq3$, inject into cohomology (Corollary~\ref{cor:Y3-6}).

\subsection{Gushel--Mukai sixfolds and their parameter space}

Let $V_5$, $\Gr=\Gr(2,V_5)$ and the cone $\CGr\subseteq\PP^{10}$ with vertex $O$ be as in Section~\ref{SS:GM4}. Following \cite{DK1} (see also \cite{Debarre}), a \emph{Gushel--Mukai sixfold} is a smooth dimensionally transverse intersection
\[
Y=\CGr\cap Q\subseteq\PP^{10},
\]
where $Q$ is a quadric. 	
The Gushel morphism $\gamma\colon Y\to\Gr$ identifies $Y$ with the double cover of $\Gr$ branched along an ordinary GM fivefold. We write $\iota_Y$ for the covering involution, and put
\[
h\coloneqq c_1(\OO_Y(1))=\gamma^*\sigma_1 \text{ and } c_2\coloneqq\gamma^*\sigma_{1,1}.
\]

Consider
\[
\bar B_6\coloneqq\PP\HH^0\bigl(\CGr,\OO(2)\bigr),
\]
and the dense open subset $B_6\subseteq\bar B_6$ parametrizing those $b$ whose zero locus $Y_b$ is a smooth sixfold.
This parametrizes GM sixfolds.
We denote by $\YY_6\to B_6$ the universal GM sixfold and by $\YY_6^m\coloneqq\YY_6\times_{B_6}\cdots\times_{B_6}\YY_6$ its relative $m$-th power.

We put $\R^*(Y)\coloneqq\langle h,c_2\rangle$ and, for $m\geq1$, $\R^*(Y^m)\coloneqq\langle h_i,c_{2,i},\Delta_{ij}\rangle\subseteq\CH^*(Y^m)$, with the notation of Section~\ref{S:tautological}. 

\subsection{Generically defined cycles on powers of a Gushel--Mukai sixfold}

\begin{lem}\label{lem:GDCH-Y6}
For every GM sixfold $Y=Y_b$, $b\in B_6$, one has $\GDCH^*_{B_6}(Y)=\langle h,c_2\rangle=\R^*(Y)$.
\end{lem}

\begin{proof}
This is the argument of Lemma~\ref{lem:GDCH-Y}, cf. also \cite[Proposition~4.4]{BL}.
\end{proof}

\begin{rmk}\label{rmk:strat6}
Proposition~\ref{prop:strat} and its proof remain valid for the smooth quasi-projective variety $P\coloneqq\CGr\smallsetminus\{O\}$, with $E\coloneqq\OO_P(2)$, with $V\coloneqq\HH^0(\CGr,\OO(2))$, which consists of the restrictions of the quadrics of $\PP^{10}$, $\CGr$ being projectively normal, and with $B_V\coloneqq B_6$. Indeed, of the hypothesis that $P$ be projective with trivial Chow groups, the proof only uses the surjectivity of $\CH^*(P)^{\otimes r}\to\CH^*(P^r)$, which holds here because $P^r\to\Gr^r$ is a vector bundle; $B_V$ only needs to be a dense open subset of $\PP(V)$ over which the zero loci are smooth and projective of dimension $\dim P-\rk E$; and \cite[Proposition~5.2]{FLV1} holds for any stratified projective bundle, its proof being localization and the projective bundle formula.
\end{rmk}

\begin{prop}\label{prop:generation6}
Let $Y=Y_b$, $b\in B_6$, be any GM sixfold. Then
\[
\GDCH^*_{B_6}(Y^m)=\R^*(Y^m)\quad\text{for }m\leq3,
\qquad\text{and}\qquad
\GDCH^i_{B_6}(Y^4)=\R^i(Y^4)\quad\text{for }i\leq 12 .
\]
\end{prop}

\begin{proof}
We apply Proposition~\ref{prop:strat} as in Remark~\ref{rmk:strat6}. Quadrics separate three distinct points of $\PP^{10}$ (see the proof of Proposition~\ref{prop:generation}), so $(P,V)$ satisfies $(\star_3)$, which gives the first equality. As in the proof of Proposition~\ref{prop:generation}, the evaluation map of condition $(\star_4)$ fails to be surjective at four pairwise distinct points of $P$ exactly when they lie on a line, and this line is contained in $\CGr$, which is cut out by the Pl\"ucker quadrics. Hence
\[
\GDCH^i_{B_6}(Y^4)\subseteq\R^i(Y^4)+\ima\bigl(\CH_*(T_4\cap Y^4)\to\CH_*(Y^4)\bigr),
\qquad
T_4\cap Y^4=\bigcup_{\ell\in F_1(\CGr)}(\ell\cap Y)^4 .
\]
A line $\ell\subseteq Y$ contributes $\ell^4$, and $\dim F_1(Y)=7$ \cite[Lemma~3.8(a)]{DK5}, whence a contribution of dimension at most $11$. A line $\ell\not\subseteq Y$, in particular every line through $O$, meets $Y$ in a subscheme of length $2$ and contributes finitely many points, and $\dim F_1(\CGr)=10$ \cite[Lemma~3.4]{DK5}. So $\dim(T_4\cap Y^4)\leq11$, its codimension in $Y^4$ is at least $13$, and $\GDCH^i_{B_6}(Y^4)=\R^i(Y^4)$ for $i\leq12$.
\end{proof}

\subsection{The canonical Chow--K\"unneth decomposition}

\begin{defn}\label{defn:CK6}
The \emph{canonical Chow--K\"unneth decomposition} of a GM sixfold $Y$ is given by $o\coloneqq\frac1{10}h^6$ and
\begin{equation}\label{eq:sixfold-projectors}
\begin{gathered}
\pi^0_Y\coloneqq o\times Y,\qquad \pi^2_Y\coloneqq\tfrac1{10}\,h^5\times h,\qquad \pi^{10}_Y\coloneqq\tfrac1{10}\,h\times h^5,\qquad \pi^{12}_Y\coloneqq Y\times o,\\
\pi^4_Y\coloneqq\tfrac12\,h^4\times h^2-h^2c_2\times h^2-h^4\times c_2+\tfrac52\,h^2c_2\times c_2,\qquad \pi^8_Y\coloneqq{}^t\pi^4_Y,\\
\pi^6_{\taut}\coloneqq\tfrac12\,h^3\times h^3-h^3\times hc_2-hc_2\times h^3+\tfrac52\,hc_2\times hc_2,
\end{gathered}
\end{equation}
$\pi^k_Y\coloneqq0$ for odd $k$, $\tau\coloneqq\Delta_Y-\sum_{k\neq6}\pi^k_Y-\pi^6_{\taut}$ and $\pi^6_Y\coloneqq\pi^6_{\taut}+\tau$. We further put $\Upsilon_Y\coloneqq \Upsilon_{c_2}= \tau\cdot p_1^*c_2\in\CH^8(Y\times Y)$.
\end{defn}

\begin{prop}\label{prop:GM6-instance}
Let $Y$ be a GM sixfold.
\begin{enumerate}[(i)]
\item $(Y,\R^*(Y))$ satisfies (H1)--(H3) of Section~\ref{S:relations}, with $d=22$, and $\CH_0(Y)=\QQ$; moreover $\HH^6_{\taut}(Y)=\langle h^3,hc_2\rangle$ has Gram matrix $\left(\begin{smallmatrix}10&4\\4&2\end{smallmatrix}\right)$.
\item $K=0$, and Definition~\ref{defn:CK-u} gives the projectors of Definition~\ref{defn:CK6}.
\item $\Upsilon_h=\tau\cdot h_1=0$ in $\CH^7(Y\times Y)$.
\item For any $m\leq45$, the kernel of the cycle class map on $\R^*(Y^m)$ consists of the ideal generated by the pullbacks to $Y^m$ of $\Upsilon_Y$ and $\Theta_Y$ and their transposes. The kernel of the cycle class map on the image of $\R^*(Y^m)$ in $\B^*(Y^m)$ is generated by the images of these classes.
\end{enumerate}
\end{prop}

\begin{proof}
(i) The variety $Y$ is Fano hence rationally connected, and so $\CH_0(Y)=\QQ$. The Hodge diamond  of $Y$ is known \cite[Proposition~3.1]{DK2}: all odd Betti numbers are 0, the even Betti numbers are $1,1,2,24,2,1,1$. 
We have equality $\R^*(Y)=\gamma^*\CH^*(\Gr)$ because $\CH^*(\Gr)=\QQ[\sigma_1,\sigma_{1,1}]$ \cite[Theorem~5.26]{EH}, and $\HH^*_{\taut}(Y)=\gamma^*\HH^*(\Gr,\QQ)$ with $\gamma^*$ injective; since the Betti numbers of $Y$ and of $\Gr$ are the same in degree $\not=6$, property (H2) holds. One has $\HH^6_{\taut}(Y)=\gamma^*\HH^6(\Gr,\QQ)=\langle h^3,hc_2\rangle$, which has the asserted Gram matrix of determinant $4$, which gives (H3). It follows that $d=24-2=22$.

\smallskip\noindent
(ii) Proposition~\ref{prop:double-cover-u}(i) applied to $\gamma\colon Y\to\Gr$ gives that $K=0$. It follows that $s$ is the inverse of the cycle class map on $\R^*(Y)$. The basis dual to $(h)$ in $\HH^{10}(Y,\QQ)$ is $(h^5/10)$, and the basis dual to $(h^2,c_2)$ in $\HH^8(Y,\QQ)$ is $(\tfrac12h^4-h^2c_2,\,-h^4+\tfrac52h^2c_2)$. The Gram matrix of $(h^3,hc_2)$ in $\HH^6(Y,\QQ)$ has inverse $\left(\begin{smallmatrix}1/2&-1\\-1&5/2\end{smallmatrix}\right)$. This gives projectors as in \eqref{eq:sixfold-projectors}.

\smallskip\noindent
(iii) We apply Lemma~\ref{lem:divisor-u}, with $c_1(L|_Y)=2h$, and $P$ the blow-up of $\CGr$ at $O$. Note that $P$ is the $\PP^1$-bundle $\PP_{\Gr}(\OO\oplus\OO(1))\to\Gr$, hence $P$ is a smooth projective variety with trivial Chow groups. Let $L_1$ be the pullback to $P$ of $\OO(1)$ from $\PP^{10}$. Since $O\notin Y$, we have $O\notin Q$, and the sixfold $Y$ is the zero locus in $P$ of the pullback of the equation of~$Q$, which is a section of $L\coloneqq L_1^{\otimes2}$ that does not vanish along the exceptional divisor. By the projective bundle formula, $\CH^*(P)$ is generated as a module over the pullback of $\CH^*(\Gr)$ by $1$ and $c_1(L_1)$.
Since $j^*L_1=\OO_Y(1)$ and the projection from $O$ restricts to $\gamma$ on $Y$, we get $j^*\CH^*(P)=\langle h,c_2\rangle=\R^*(Y)$. As $K=0$, Lemma~\ref{lem:divisor-u} gives $\Upsilon_{2h}=2\,\Upsilon_h=0$.

\smallskip\noindent
(iv) This is Proposition~\ref{prop:relations-generated-u} with $G=\{h,c_2\}$ and $m_Y=2d+1=45$, together with $\CH_0(Y)=\QQ$, $K=0$ and (iii). The kernel of the cycle class map on the image of $\R^*(Y^m)$ in $\B^*(Y^m)$ is the image of its kernel on $\R^*(Y^m)$.
\end{proof}

\begin{prop}\label{prop:sixfold-tau-u}
Let $Y$ be a GM sixfold, with covering involution $\iota_Y$. Then the cycle class map is injective on $\R^*(Y)=\langle h,c_2\rangle=\GDCH^*_{B_6}(Y)$, and
\[
\tau=\tfrac12\bigl(\Delta_Y-\Gamma_{\iota_Y}\bigr)\ \hbox{ in }\CH^6(Y\times Y),\ 
\tau\cdot h_1=0\ \hbox{ in }\CH^7(Y\times Y),\ 
M(\tau,\tau;\tau)=0\ \hbox{ in }\CH^{12}(Y^3).
\]
Moreover, the canonical Chow--K\"unneth decomposition of $Y$ is multiplicative if and only if $\Upsilon_Y=0$ in $\CH^8(Y\times Y)$, and it is multiplicative modulo algebraic equivalence as soon as $\Upsilon_Y=0$ in $\B^8(Y\times Y)$.
\end{prop}

\begin{proof}
This follows from Proposition~\ref{prop:double-cover-u} with $P=\Gr$, together with Proposition~\ref{prop:GM6-instance}(i)--(iii). The equality $\R^*(Y)=\GDCH^*_{B_6}(Y)$ follows from Lemma~\ref{lem:GDCH-Y6}.
\end{proof}

\begin{cor}\label{cor:sixfold-CH-u}
Let $Y$ be a GM sixfold such that $\Upsilon_Y=0$ in $\CH^8(Y\times Y)$. Then the canonical Chow--K\"unneth decomposition of $Y$ is multiplicative. 

Consequently, the cycle class map is injective on $\R^*(Y^m)$ for all $m\leq45$. In particular, the cycle class map is injective on $\GDCH^*_{B_6}(Y^m)$ for $m\leq3$, and on $\GDCH^i_{B_6}(Y^4)$ for $i\leq12$.
\end{cor}

\begin{proof}
Proposition~\ref{prop:sixfold-tau-u} proves multiplicativity of the decomposition, and implies the vanishing $K=0$, $\tau\cdot h_1=0$. In view of Proposition~\ref{prop:double-cover-u}(iv), we also have the vanishing $\Theta_Y=0$ in $\CH^{12}(Y^3)$. Applying Corollary~\ref{cor:relations-vanish-u}, we then find the cycle class map is injective on $\R^*(Y^m)$ for $m\leq45$. 
Lemma~\ref{lem:GDCH-Y6} and Proposition~\ref{prop:generation6} identify $\GDCH^*_{B_6}(Y^m)$ with $\R^*(Y^m)$ for $m\leq3$ and $\GDCH^i_{B_6}(Y^4)$ with $\R^i(Y^4)$ for $i\leq12$, proving the last statement.
\end{proof}

\subsection{Relations modulo algebraic equivalence}

\begin{prop}\label{prop:sixfold-alg-u}
Let $Y$ be a GM sixfold. Then $\Upsilon_Y=0$ in $\B^8(Y\times Y)$. Consequently the canonical Chow--K\"unneth decomposition of $Y$ is multiplicative modulo algebraic equivalence.

For every $m\leq45$ the image of $\R^*(Y^m)$ in $\B^*(Y^m)$ injects into $\HH^{2*}(Y^m,\QQ)$. In particular, the cycle class map is injective on $\GDB^*_{B_6}(Y^m)$ for $m\leq3$, and on $\GDB^i_{B_6}(Y^4)$ for $i\leq12$.
\end{prop}

\begin{proof}
Let $\Omega$ be a universal domain containing $\C$, abstractly isomorphic to $\C$. We know that $\CH_i^{\mathrm{hom}}(Y_\Omega)=0$ for $i\le 1$ \cite[Theorem~5.1]{BL}, and so $\CH_0(Y_\Omega)=\QQ\,o$ and $\CH_1(Y_\Omega)=\QQ\,h^5$. Hence the effective motive $M:=(Y,\Delta_Y-\pi^0_Y-\pi^2_Y-\pi^{10}_Y-\pi^{12}_Y)$ satisfies $\CH_0(M_\Omega)=\CH_1(M_\Omega)=0$, and the same holds for the dual of $M$.
 Thanks to \cite[Corollary~2.2]{VialCK} (applied with $m=n=2$), $M$ is then a direct summand of $\h(S)(-2)$ for some smooth projective surface $S$. 
Hence, we obtain a split injection
$$
\h(Y)\ \hookrightarrow\ \h(S)(-2)\oplus\bigoplus_k\1(b_k)\qquad\text{ in }\Mrat.
$$
Tensoring this injection with itself shows that $\h(Y\times Y)$ is a direct summand of the direct sum of $\h(S\times S)(-4)$ with Tate twists of $\h(S)$ and of $\1$. Hence $\B^8(Y\times Y)$ is motivically a direct sum of $\B^4(S\times S)$ and copies of $\B^*(S)$ and $\QQ$. As $\Upsilon_Y$ is homologically trivial by Proposition~\ref{prop:CK-basic-u}(iii), it follows from the fact that  algebraic and homological equivalence agree for cycles of dimension 0 and for cycles of codimension at most 1 that $\Upsilon_Y=0$ in $\B^8(Y\times Y)$.
The canonical Chow--K\"unneth decomposition is then multiplicative modulo algebraic equivalence by Proposition~\ref{prop:sixfold-tau-u}. 

We know that $\Theta_Y=0$ in $\B^{12}(Y^3)$ by Proposition~\ref{prop:double-cover-u}(iv); since moreover $K=0$ and $\tau\cdot h_1=0$, Corollary~\ref{cor:relations-vanish-u} (with $G=\{h,c_2\}$ and $d=22$) implies the injectivity of the cycle class map on the image of $\R^*(Y^m)$ in $\B^*(Y^m)$ for $m\leq45$. 

Finally, by Lemma~\ref{lem:GDCH-Y6} and Proposition~\ref{prop:generation6} we have equalities 
 $\GDCH^*_{B_6}(Y^m)=\R^*(Y^m)$ for $m\leq3$ and $\GDCH^i_{B_6}(Y^4)=\R^i(Y^4)$ for $i\leq12$, so that $\GDB^*_{B_6}(Y^m)$, respectively $\GDB^i_{B_6}(Y^4)$, is the image of $\R^*(Y^m)$, respectively $\R^i(Y^4)$, in $\B^*(Y^m)$.
\end{proof}

Proposition~\ref{prop:sixfold-alg-u} recovers \cite[Theorems~6.4 and~7.2]{BL}.

\subsection{Zero-cycles on the double EPW sextic and a motivic splitting}

\begin{rmk}\label{rmk:coarse} Let $M_{\mathrm{EPW}}$ denote the coarse moduli space of double EPW sextics (this is a GIT quotient
$M_{\mathrm{EPW}}=LG^\circ//PGL_6$, where $LG^\circ$ is an open subset of a Lagrangian Grassmannian \cite{O'G16}). Over a dense open subset $M_{\mathrm{EPW}}^\circ\subseteq M_{\mathrm{EPW}}$, it
 is possible to construct a tautological family of double EPW sextics (\cite[Section 8.1]{Ma22}, where this tautological family is denoted ${\mathcal F}\to{\mathcal M}$).
 It thus makes sense to speak about generically defined cycles with respect to $M^\circ_{\mathrm{EPW}}$, and Theorem \ref{thm:zero-cycles-X} implies that for any double EPW sextic $X$ parameterized by $M_{\mathrm{EPW}}^\circ$ one has   \[ \GDCH^4_{M^\circ_{\mathrm{EPW}}}(X) =\QQ.\]
 \end{rmk}

The analogue of the motivic splitting of Corollary~\ref{cor:sym2-splitting} exists for GM sixfolds, relative to the moduli stack ${\mathcal M}_{\mathrm{GM6}}$ of GM sixfolds (this moduli stack is constructed in \cite{DK3}):

\begin{cor}\label{cor:sym2-splitting6}
Let $Y=Y_b$ be a GM sixfold, and let $X=X_b$ be the associated double EPW
sextic. There is a split injection in $\Mrat$
\[\Phi \colon \
\h(Y^{(2)})\hookrightarrow \h(X)(-4)\oplus\bigoplus_j\h(Y)(a_j)\oplus\bigoplus_k\1(b_k)
\]
for suitable Tate twists $a_j,b_k$, such that both the injection and a left inverse are induced by
correspondences that are generically defined with respect to ${\mathcal M}_{\mathrm{GM6}}$.
\end{cor}

\begin{proof} Debarre--Kuznetsov \cite[Theorem 5.19]{DK2} have constructed an Abel--Jacobi isomorphism
   \[ \HH^2_{\prim}(X_b,\QQ)\xrightarrow{\cong}\HH^6_{\van}(Y_b,\QQ)(2)\ .\]
   The construction is via the universal $\sigma$-plane, and so this isomorphism is induced by a correspondence generically defined over ${\mathcal M}_{\mathrm{GM6}}$ (and the inverse isomorphism is also defined by such a correspondence) \cite[Theorem 5.19]{DK2} (cf. also \cite[Proposition 2.8]{BL}, where this construction is used relative over the parameter space $B_6$). The corollary is now proven just as Corollary \ref{cor:sym2-splitting}.
\end{proof}

\subsection{The MCK decomposition for GM sixfolds}

\begin{thm}\label{thm:MCK6}
The canonical Chow--K\"unneth decomposition of a GM sixfold is multiplicative.
\end{thm}

\begin{proof}
Let $Y=Y_b$, $b\in B_6$, be a GM sixfold. By Proposition~\ref{prop:sixfold-tau-u} it suffices to
prove that $\Upsilon_Y=0$ in $\CH^8(Y^2)$.

The class $\Upsilon_Y=\tau\cdot c_{2,1}$ is homologically trivial
(Proposition~\ref{prop:CK-basic-u}(iii)) and generically defined with respect to ${\mathcal M}_{\mathrm {GM6}}$.
It is readily seen that $\Upsilon_Y$ is invariant under the permutation of the two factors of $Y^2$.
Consequently $\Upsilon_Y$ lies in the image of the symmetrizing projector, that is, it defines a
homologically trivial and generically defined class $\overline\Upsilon_Y$ in $\CH^8$ of the motive
$\h(Y^{(2)})$.

Let us now apply Corollary~\ref{cor:sym2-splitting6}. Since
$\Phi$ and its left inverse are generically defined with respect to ${\mathcal M}_{\mathrm {GM6}}$, the components of
$\Phi_*\overline\Upsilon_Y$ are homologically trivial classes lying respectively in
$\GDCH^4_{{\mathcal M}_{\mathrm{GM6}}}(X_b)$, in the groups $\GDCH^{8+a_j}_{{\mathcal M}_{\mathrm{GM6}}}(Y)$ and in $\QQ$. The second groups inject into
cohomology by Proposition~\ref{prop:sixfold-tau-u}, and the third ones inject trivially.
Thus it only remains to prove that
  \begin{equation}\label{overGM6} \GDCH^4_{{\mathcal M}_{\mathrm{GM6}}}(X_b)=\QQ.\end{equation}
To this end, we claim that the map of coarse moduli spaces
$ \nu_6\colon M_{\mathrm{GM6}} \to M_{\mathrm{EPW}} $
  constructed in \cite[Section 6.1]{DK3}
  induces an equality
   \begin{equation}\label{basechange6}  \GDCH^\ast_{M_{\mathrm{EPW}}}(X)= \GDCH^\ast_{M_{\mathrm{GM6}}}(X)\ .\end{equation}
   In view of the inclusion $\GDCH^\ast_{{\mathcal M}_{\mathrm{GM6}}}(X) \subseteq \GDCH^\ast_{M_{\mathrm{GM6}}}(X)$
   and Remark \ref{rmk:coarse}, this claim implies equality \eqref{overGM6}.
 
 In order to prove the claim, we observe that the general fiber of the morphism $\nu_6$ is of the form $U/G$, where $U$ is the complement of an EPW sextic in $\PP^5$ and $G$ is a finite group \cite[Corollary 6.3]{DK3}.  
In view of the spread lemma (Lemma \ref{lem:base-change}(i)), to prove the claim it suffices to prove it after replacing $M_{\mathrm{EPW}}$ and $M_{\mathrm{GM6}}$ by some dense open subsets $M_{\mathrm{EPW}}^\circ$ resp. $M_{\mathrm{GM6}}^\circ$. Thus, choosing convenient open subsets we may assume that all fibers of $\nu_6$ are of the form $U/G$, and also (by generic flatness) that $\nu_6$ is a flat morphism. The claim then follows from Lemma \ref{lem:charles} below, applied to $\XX_{M_{\mathrm{GM6}}^\circ}\to\XX_{M_{\mathrm{EPW}}^\circ}$.
 \end{proof}
 
 \begin{lem}\label{lem:charles} Let $f\colon M\to N$ be a flat morphism among smooth quasi-projective varieties. Assume that every fiber $F$ of $f$ satisfies $\CH^0(F)=\QQ$ and $\CH^i(F)=0$ for $i>0$.
 Then $f^\ast\colon \CH^\ast(N)\to \CH^\ast(M)$ is surjective.
 \end{lem}
 
 \begin{proof} We prove this by induction on $\dim N$, the case $\dim N=0$ being trivial. 
 Let $N^{(1)}$ denote the codimension 1 subschemes of $N$, and let $\mu_N={\mathrm{Sp}} \C(N)$ be the generic point of~$N$.
 Flatness of $f$ implies there is a commutative diagram with exact rows
 \begin{equation*} \xymatrix{ {\displaystyle\bigoplus_{D\in N^{(1)}}} \CH^{i-1}(f^{-1}(D)) \ar[r] & \CH^i(M) \ar[r]  & \CH^i( M_{\mu(N)})  \ar[r] & 0  \\
  	{\displaystyle\bigoplus_{D\in N^{(1)}}}   \CH^{i-1}(D ) \ar[u]_{\oplus f_D^\ast} \ar[r] & \CH_i(N) \ar[u]_{f^\ast} \ar[r]  & \CH_i(\mu_N)   \ar[u]_{f^\ast_{\mu_N}}	 \ar[r] & 0
  }
\end{equation*} 
The right vertical arrow is surjective, and the left vertical arrow is surjective by the induction hypothesis. A diagram chase then gives surjectivity of the middle vertical arrow.
\end{proof}

\begin{cor}\label{cor:Y3-6}
Let $Y=Y_b$ be a GM sixfold, where $b\in B_6$. Then the cycle class map is injective on $\R^*(Y^m)$ for $m\leq45$, hence on $\GDCH^*_{B_6}(Y^m)$ for $m\leq3$ and on $\GDCH^i_{B_6}(Y^4)$ for $i\leq12$.
\end{cor}

\begin{proof}
By Theorem~\ref{thm:MCK6} and Proposition~\ref{prop:sixfold-tau-u}, $\Upsilon_Y=0$ in $\CH^8(Y\times Y)$, so that Corollary~\ref{cor:sixfold-CH-u} applies.
\end{proof}

\addcontentsline{toc}{section}{References}
\bibliographystyle{alpha} 
\bibliography{References}{}

\medskip \medskip

\noindent Lie Fu\\
\noindent{Universit\'e de Strasbourg, Institut de Recherche Mathématique Avancée (IRMA) and Institut Universitaire de France (IUF) } \\
\medskip \noindent{\texttt{lie.fu@math.unistra.fr}}

\noindent Robert Laterveer\\
\noindent{ Universit\'e de Strasbourg, Institut de Recherche Math\'ematique Avanc\'ee (IRMA) and CNRS      }\\
\medskip \noindent \texttt{laterv@math.unistra.fr}  

\noindent Carl Mazzanti\\
\noindent{Fakult\"at f\"ur Mathematik, Universit\"at Bielefeld, Germany}\\
\medskip \noindent \texttt{cmazzanti@math.uni-bielefeld.de}

\noindent Charles Vial\\
\noindent{Fakult\"at f\"ur Mathematik, Universit\"at Bielefeld, Germany}\\
\medskip \noindent \texttt{vial@math.uni-bielefeld.de}  
\end{document}